# On the Minimum Many-Valued Modal Logic over a Finite Residuated Lattice

FÉLIX BOU, FRANCESC ESTEVA and LLUÍS GODO, *Institut d'Investigació en Intel.ligència Artificial, IIIA - CSIC, Campus UAB, Bellaterra 08193, Spain. E-mail: fbou@iiia.csic.es; esteva@iiia.csic.es; godo@iiia.csic.es*

RICARDO OSCAR RODRÍGUEZ, *Dpto. de Computación, Fac. Ciencias Exactas y Naturales, Universidad de Buenos Aires, 1428 Buenos Aires, Argentina. E-mail: ricardo@dc.uba.ar*

## Abstract

This article deals with many-valued modal logics, based only on the necessity operator, over a residuated lattice. We focus on three basic classes, according to the accessibility relation, of Kripke frames: the full class of frames evaluated in the residuated lattice (and so defining the minimum modal logic), the ones only evaluated in the idempotent elements and the ones evaluated in 0 and 1. We show how to expand an axiomatization, with canonical truth-constants in the language, of a finite residuated lattice into one of the modal logic, for each one of the three basic classes of Kripke frames. We also provide axiomatizations for the case of a finite MV chain but this time without canonical truth-constants in the language.

## 1  Introduction

The basic idea of this article is to systematically investigate modal extensions of many-valued logics. Many-valued modal logics, under different forms and contexts, have appeared in the literature for different reasoning modeling purposes. For instance, Fitting introduces in [22, 23] a modal logic over logics valued on finite Heyting algebras, and provides a nice justification for studying such modal systems for dealing with opinions of experts with a dominance relation among them. But many papers offer technical contributions but without practical motivations. Although we will also mainly focus on developing a theoretical and general framework, let us say something about specific problems where the necessity of combining modal and many-valued semantics seems to be important. Our starting point is that the framework of classical logic is not enough to reason with vague concepts or with modal notions like belief, uncertainty, knowledge, obligations, time, etc. Residuated fuzzy (or many-valued) logics (in the sense of Hájek in [30]) appear as a suitable logical framework to formalise reasoning under vagueness, while a variety of modal logics address the logical formalization to reason about different modal notions as the ones mentioned above. Therefore, if one is interested in a logical account of both vagueness and some kind of modalities, one is led to study systems of many-valued or fuzzy modal logic.

A clear example of this are Fuzzy Description Logics (see e.g. [46] for an overview) which, as in the classical case, can be viewed both as fragments of t-norm fuzzy predicate logics (see e.g. Hájek's paper [31]) and as fuzzy (multi) modal logics. Other examples include generalizations of belief logics to deal with non-classical/fuzzy events [36, 29] or preliminary ideas on similarity-based fuzzy logics [27].

Having in mind these motivations, our main research goal at large is a systematic study of *many-valued modal logics* and their application, including the setting of fuzzy description logics. Although there is a rich literature about both modal [11, 4, 5] and many-valued logics [30, 13, 28], in our opinion there is no such systematic approach. This article may be seen as a first modest step in that direction, indeed, we restrict ourselves to investigate minimum many-valued modal logics for the necessity operator $\Box$ defined on top of logics of finite residuated lattices.

It is certainly true that several many-valued modal logics have been previously considered in the literature; but in most cases, with the two exceptions later cited, they are logics (we highlight [35, 30, 39] among them) that site at the top part of the hierarchy of modal logics, e.g., the modality is S5-like. In other words, they cannot be considered as the many-valued counterpart of the minimum classical modal logic K. Since in our opinion any systematic study of many-valued modal logics must start, like in the classical modal case, by characterizing the corresponding minimum many-valued modal logic we consider that this problem is the first question that someone who wants to study many-valued modal logics must answer. This is the reason why *the present article is focused on characterizing the minimum many-valued modal logic.* This does not mean that the authors do not consider important to study the full hierarchy of many-valued modal logics extending the minimum one, it is simply that this systematic study is left for future research.

Let us fix a residuated lattice $\mathbf{A}$. The intended meaning of its universe $A$ is the set of truth values. Thus, the classical modal logic framework corresponds to the case that $\mathbf{A}$ is the Boolean algebra of two elements. Due to technical reasons sometimes it will be convenient to add canonical constants to the residuated lattice $\mathbf{A}$. By adding *canonical constants* we mean to add one constant for every element in $A$ in such a way that these constants are semantically interpreted by its canonical interpretation. The residuated lattice obtained by adding these canonical constants will be denoted by $\mathbf{A^c}$.

The intuition behind the minimum many-valued modal logic over $\mathbf{A}$ (and analogously over $\mathbf{A^c}$) is that it has to be induced by the biggest class of Kripke frames. Therefore, the minimum logic must be semantically defined by using *Kripke frames* where the accessibility relation takes values in $A$. Thus, the accessibility relation in these Kripke frames is also many-valued. In other words, the accessibility relation may take values outside $\{0, 1\}$, and hence we cannot assume that it is a subset of $A \times A$.

As far as the authors are aware the only two cases in the literature that deal with a minimum modal logic over a residuated lattice $\mathbf{A}$, and hence with many-valued accessibility relations, are the ones respectively introduced in [22] and [10, 47]. The first one studies the case that $\mathbf{A}$ is a finite Heyting algebra (but adding canonical constants) and the second one studies the standard Gödel algebra (which is also a Heyting algebra). One particular feature about these two cases is that the normality

axiom

$$\Box(\varphi \to \psi) \to (\Box\varphi \to \Box\psi) \tag{K}$$

is valid in the minimum modal logic. However, there are a lot of residuated lattices where the normality axiom (K) is not valid in the minimum modal logic, e.g., finite MV chains. Thus, whenever (K) is not valid there is no axiomatization in the literature of the minimum modal logic. Even more, if (K) fails it is not so clear what could be a good candidate for an axiomatization.

The present article solves the problem of finding the minimum (local) modal logic, based only on the necessity operator $\Box$, for two cases: for finite residuated lattices with canonical constants and for finite MV chains without adding canonical constants. For these two cases we show *how to expand an axiomatization for the logic of the residuated lattice* **A** *into one of the minimum modal logic over* **A**. Thus, we assume in this article that we already know how to axiomatize the non-modal logic.

Besides the minimum (local) many-valued modal logic (i.e., the logic given by all Kripke frames) in this article we also consider the (local) many-valued modal logic given by two other classes of frames. These two other classes are the class of *crisp frames* where the accessibility relation takes values in $\{0, 1\}$, and the class of *idempotent frames* where the accessibility relation takes values in the set of idempotent elements of **A**. For every finite residuated lattice with canonical constants we also axiomatize in this article the logic given by idempotent frames, but for the logic given by crisp frames we have only succeeded when $\mathbf{A}^{\mathbf{c}}$ has a unique coatom. Our interest on the class of crisp frames is due to its obvious connection with the classical setting, and the interest on idempotent frames is based on the tight connection with the normality axiom (K) that will be shown later in the article.

In this article we stress that the non-modal logics considered are the ones defined over a particular residuated lattice and not over a class of residuated lattices[1] The reason to do so is a methodological one since a lot of our results are based on the role of canonical constants, which can only be considered in the framework of non-modal logics defined using only one residuated lattice. On the other hand, the restriction to finite residuate lattices is used in order to be sure that the non-modal logic is finitary.

It is worth pointing out that in this article logics are always (in the non-modal and in the modal setting) considered as consequence relations invariant under substitutions. So, the reader must understand in this article "completeness" as a synonym of what it is sometimes called "strong completeness". That is, we are talking about completeness even for infinite sets of assumptions. For the same reason we use the expression local and global modal logic, and not the terminology of "local and global consequence relations" which is more common in the modal tradition.

**Structure and contents of the article.** Section 2 contains the non-modal preliminaries; the definition of residuated lattice is stated and it is said what is the logic associated with a residuated lattice. For certain modal aspects it will be convenient to add canonical constants to the residuated lattice, and so in Section 2.2 we also consider the non-modal logic with canonical constants.

---

[1] It is common in the non-modal literature to consider equational classes of residuated lattices. However, this assumption does not seem reasonable to be done in the many-valued setting. For instance, in [8, p. 181] it is shown an example of two classes of residuated lattices generating the same equational class, but having different minimum many-valued modal logics.

In Section 3 we firstly introduce the Kripke semantics on the many-valued setting (with and without canonical constants). We notice that for the sake of simplicity in this article we only deal with the necessity modality $\Box$. The Kripke semantics is defined by

$$V(\Box\varphi, w) = \bigwedge\{R(w, w') \rightarrow V(\varphi, w') : w' \in W\},$$

and it is meaningful for the case that $\mathbf{A}$ is lattice complete. Although there are several possibilities to generalize the classical Kripke semantics to the many-valued realm we justify in Remark 3.5 why we think that this semantics is the most natural one; indeed, under this definition it holds that the semantics of the modal formula $\Box p$ corresponds to the first-order many-valued semantics of the formula $\forall y(Rxy \rightarrow Py)$. Then, the three basic classes of Kripke frames considered in this article (the full class of Kripke frames, idempotent Kripke frames and crisp Kripke frames) are introduced. In Section 3.1 the validity of some formulas, in the three basic classes of Kripke frames, is discussed. Among the validity results Corollary 3.13 is the more remarkable one. It says that the normality axiom ($\mathsf{K}$) is valid iff all elements in the residuated lattice are idempotent (i.e., $\mathbf{A}$ is a Heyting algebras). In particular this result shows a tight connection between ($\mathsf{K}$) and the class of idempotent frames. Besides the three basic classes of Kripke frames, we also introduce in this section another class (the class of Boolean Kripke frames) which is tightly connected to the class of crisp Kripke frames (see Theorem 3.16). Next, Section 3.2 introduces the local and the global many-valued modal logic defined by a class of Kripke frames. In this Section 3.2 we discuss several (meta)rules whose role is crucial in the modal setting, we show how to reduce the modal logic to the non-modal one (see Theorem 3.26) and we study the relationship between crisp and Boolean Kripke frames. In Section 3.3 we sum up the previous known results about the minimum many-valued modal logics. This *state of the art* is locate there, and not at the beginning of the article, in order to be able to state the already known results inside the notational framework introduced in this article.

Section 4 considers modal logics for the case that $\mathbf{A}$ is a finite residuated lattice and there are canonical constants. Firstly, in Section 4.1 we state several results illustrating the benefits of having canonical constants (these results work even for non finite residuated lattices). From these results it is clear that the presence of canonical constants simplifies proofs (but we must be really careful because some properties, like the Local Deduction Theorem, can be destroyed by adding canonical constants). Moreover, we also think that the presence of canonical constants helps to understand proofs. For example, the completeness proofs in Section 4 help to understand and clarify why the proofs in Section 5 work; that is, it is easier to understand the proofs in Section 5 after understanding the ones with canonical constants. In Sections 4.2, 4.3 and 4.4 we show, for the classes of Kripke frames and idempotent frames, how to expand an axiomatization of the non-modal logic of $\mathbf{A}$ with canonical constants into one of the local many-valued modal logic over $\mathbf{A}$ with canonical constants. In the case of crisp frames we only know how to do this expansion in case that $\mathbf{A}$ is a finite residuated lattice with a unique coatom. The proofs of these completeness results are based on the canonical frame construction. For the global many-valued modal logic we only know how to axiomatize the case of crisp frames.

In Section 5 we adapt the ideas from the previous section to the case that $\mathbf{A}$ is a finite $\mathsf{MV}$ chain and there are no canonical constants in the language. The trick is to realize that for finite $\mathsf{MV}$ chains we can somehow internalize the elements of the lattice using strongly characterizing formulas (see Remark 4.5). For the case of a finite MV chain we give axiomatizations of the local many-valued modal logic for each one of the basic classes of Kripke frames (in this case idempotent frames coincide with crisp ones). As in the previous section we have only succeed to axiomatize the global logic in the case of crisp frames. Although the logics defined by crisp frames over a finite MV chain were already axiomatized by HANSOUL and TEHEUX in [39] we have included them in this article because the proofs here given are different, and because we think that these new proofs are more intuitive.

In Section 6 we state what in our opinion are the main left open problems in this research field.

This article also contains two appendixes. Appendix A contains several remarks about the non-modal logics. Most of the results given in this appendix are well known in the literature. One exception is Theorem A.14. This new result provides a method to expand an axiomatization of the non-modal logic of **A** into one of the non-modal logic with canonical constants. Finally, Appendix B explains the non-modal companion framework as a simple method to discard the validity of some modal formulas (including the normality axiom).

**Notation.** Let us briefly explain the main conventions about notation that we assume in this article. Of course, in all these conventions sometimes subscripts and superscripts will be used. Algebras will be denoted by $\mathbf{A}, \mathbf{B}, \ldots$; their universes by $A, B, \ldots$; and their elements by $a, b, \ldots$. For the particular case of the algebra of formulas we will use $\mathbf{Fm}$, while the set of formulas will be denoted by $Fm$. The formulas, generated by variables $p, q, r, \ldots$, will be denoted by $\varphi, \psi, \gamma, \ldots$; and sets of formulas by $\Gamma, \Delta, \ldots$. In an algebraic disguise formulas are sometimes called terms; whenever we follow this approach $x, y, \ldots$ will denote term variables, $t, u, \ldots$ will denote terms and $t^{\mathbf{A}}, u^{\mathbf{A}}, \ldots$ will denote its interpretation in an algebra. The set of homomorphisms from algebra $\mathbf{A}$ into $\mathbf{B}$ is abbreviated to $Hom(\mathbf{A}, \mathbf{B})$. Operators $\mathbb{H}, \mathbb{S}, \mathbb{P}$ and $\mathbb{P}_U$ refer to the closure of a class of algebras under, respectively, homomorphic images, subalgebras, direct products and ultraproducts. And operators $\mathbb{Q}$ and $\mathbb{V}$ denote, respectively, the quasivariety and the variety generated by a class of algebras. Finally, we remark that while $f(x)$ denotes, as usual, the image under a map $f$ of an element $x$ in the domain, we will use $f[X]$ to denote the image set under the same map of a subset $X$ of the domain.

**Acknowledgements.** The collaborative work between the authors of this article was made possible by several research grants: CyT-UBA X484, research CONICET program PIP 5541, AT Consolider CSD2007-0022, LOMOREVI Eurocores Project FP006/FFI2008-03126-E/FILO, MULOG2 TIN2007-68005-C04-01 of the Spanish Ministry of Education and Science, including FEDER funds of the European Union, and 2009SGR-1433/1434 of the Catalan Government. The authors also wish to thank the anonymous referees for their helpful comments.

## 2   Non-modal preliminares

All algebras considered in this article are *residuated lattices*. By definition residuated lattices are algebras given in the propositional (algebraic) language $\langle \wedge, \vee, \odot, \rightarrow, 1, 0 \rangle$ of arity $(2, 2, 2, 2, 0, 0)$. In the future sections we will often refer to this language (perhaps expanded with constants) simply as the *non-modal* language. We will consider the connectives $\neg$ and $\leftrightarrow$ as they are usually defined: $\neg \varphi := \varphi \to 0$ and $\varphi \leftrightarrow \psi := (\varphi \to \psi) \odot (\psi \to \varphi)$. Other abbreviations (where $m \in \omega$) are $\varphi \oplus \psi := \neg(\neg \varphi \odot \neg \psi)$, $\varphi^m := \varphi \odot \overset{m}{\ldots} \odot \varphi$ and $m.\varphi := \varphi \oplus \overset{m}{\ldots} \oplus \varphi$. By definition an algebra $\mathbf{A} = \langle \wedge, \vee, \odot, \to, 1, 0 \rangle$ is a *residuated lattice* if and only if the reduct $\langle A, \wedge, 1, 0 \rangle$ is a bounded lattice with maximum 1 and minimum 0 (its order is denoted by $\leqslant$), the reduct $\langle A, \odot, 1 \rangle$ is a commutative monoid, and the *fusion* operation $\odot$ (sometimes also called the *intensional conjunction* or *strong conjunction*) is residuated, with $\to$ being its residuum; that is, for all $a, b, c \in A$,

$$a \odot b \leqslant c \iff b \leqslant a \to c. \tag{2.1}$$

In the literature these lattices are also well known under other names: e.g., integral, commutative residuated monoids and $\mathsf{FL_{ew}}$-algebras [40, 41, 49]. It is worth pointing out that the class of residuated lattices $\mathsf{RL}$ is a variety. Among its subvarieties are the ones widely studied in the context of fuzzy logics and of substructural logics. Some interesting subvarieties, which will appear later on, are: $\mathsf{HA}$, the subvariety of Heyting algebras, obtained from $\mathsf{RL}$ by adding idempotency of fusion $x \odot x \approx x$; $\mathsf{MTL}$, the subvariety of $\mathsf{RL}$ obtained by adding the prelinearity condition $(x \to y) \vee (y \to x) \approx 1$; $\mathsf{BL}$, the subvariety of $\mathsf{MTL}$ obtained by adding the divisibility condition $x \odot (x \to y) \approx x \wedge y$; $\mathsf{MV}$, the subvariety of $\mathsf{BL}$ obtained by making the negation involutive (i.e., $\neg\neg x \approx x$); $\mathsf{\Pi}$, the variety of product algebras, which is the subvariety of $\mathsf{BL}$ obtained by adding the cancellative property $\neg\neg x \to ((x \to (x \odot y)) \to (y \odot \neg\neg y)) \approx 1$; and the variety $\mathsf{G}$ of Gödel or Dummett algebras, obtained from $\mathsf{BL}$ by adding idempotency of fusion $x \odot x \approx x$. The main sources of information concerning these varieties (and related ones) are [13, 19, 25, 28, 30, 40, 49].

For modal purposes later we will require completeness in the residuated lattices. This *completeness* is the common lattice-theoretic property stating that all suprema and infima (even of infinite subsets of the domain) exist. It is well known that complete residuated lattices satisfy the law

$$x \odot \bigvee y_i \approx \bigvee (x \odot y_i) \tag{2.2}$$

In the non complete case the previous law also holds as far as we assume that the corresponding suprema exists. A remarkable consequence is that joins of idempotent elements are idempotent. Other facts obtained from (2.2) are

$$x \to (\bigwedge y_i) \approx \bigwedge (x \to y_i) \tag{2.3}$$

$$(\bigvee x_i) \to y \approx \bigwedge (x_i \to y) \tag{2.4}$$

However, the next two laws

$$x \to (\bigvee y_i) \approx \bigvee (x \to y_i) \tag{2.5}$$

$$(\bigwedge x_i) \to y \approx \bigvee (x_i \to y) \tag{2.6}$$

are false in general. Indeed, the class of complete residuated lattices satisfying (2.5) and (2.6) is precisely the class of complete $\mathsf{MTL}$ algebras (see [40]).

We notice that all finite lattices satisfy the completeness condition. Other examples of complete residuated lattices are *standard* BL-*algebras*. These are exactly the BL algebras over the real unit interval $[0, 1]$ (with the usual order), which are also characterized as the residuated lattices having a continuous t-norm as the $\odot$ operator. We remind that the three basic examples are:

- Łukasiewicz: $x \odot y := \max\{0, x + y - 1\}$ and $x \to y := \min\{1, 1 - x + y\}$.
- Gödel: $x \odot y := x \wedge y$ and $x \to y := 1$ if $x \leqslant y$ and $0$ if $y < x$.
- Product: $x \odot y := x \cdot y$ and $x \to y := 1$ if $x \leqslant y$ and $\frac{y}{x}$ if $y < x$.

We will denote these algebras respectively by $[\mathbf{0, 1}]_{\mathbf{\L}}$, $[\mathbf{0, 1}]_{\mathbf{G}}$ and $[\mathbf{0, 1}]_{\mathbf{\Pi}}$. Finally, we point out that for every $n \geqslant 2$ the set $\{\frac{m}{n-1} : m \in \{0, \ldots, n-1\}\}$ is the support of a subalgebra of $[\mathbf{0, 1}]_{\mathbf{\L}}$. This algebra is the only, up to isomorphism, MV chain with $n$ elements, and we will denote it by $\mathbf{L_n}$.

### 2.1   *The (non-modal) logic of a residuated lattice*

For every residuated lattice $\mathbf{A}$, there is a natural way to associate a logic. Following modern algebraic logic literature in this article we consider a logic as a consequence relation closed under substitutions. The main idea behind its definition is that for finitary deductions we are essentially capturing the quasivariety generated by $\mathbf{A}$. The logic associated with $\mathbf{A}$ will be denoted by $\mathbf{\Lambda(A)}$, and is obtained by putting, for all sets $\Gamma \cup \{\varphi\}$ of formulas,

$$\Gamma \vdash_{\mathbf{A}} \varphi \quad \iff \quad \forall h \in \mathrm{Hom}(\mathbf{F}m, \mathbf{A}), \text{ if } h[\Gamma] \subseteq \{1\} \text{ then } h(\varphi) = 1. \qquad (2.7)$$

Our aim in this article is to explain how to expand an axiomatization of the non-modal logic $\mathbf{\Lambda(A)}$ into one for the modal logic (that we will later introduce). Hence, throughout this article we assume that we have an axiomatization of $\mathbf{\Lambda(A)}$. In other words, our concerns are not focused on the non-modal part. However, the reader interested in the non-modal logic $\mathbf{\Lambda(A)}$ can take a look at Appendix A.1, where there is a summary of known results about this logic together with several references to already known axiomatizations for certain cases (which include finite and standard BL algebras).

### 2.2   *Adding canonical constants to the non-modal logic*

Next we are going to consider the expansion of the logic $\mathbf{\Lambda(A)}$ with a constant for every truth value (i.e., element of $\mathbf{A}$) such that each one of these constants will be interpreted as the corresponding truth value. In the rest of this section we explain the details of this approach.

For every residuated lattice $\mathbf{A}$ we consider its *expanded language by constants* as the expansion of the algebraic language of residuated lattices obtained by adding a new constant $\overline{a}$ for every element $a \in A$. The *canonical expansion* $\mathbf{A^c}$ (cf. [18]) of $\mathbf{A}$ is the algebra in this very expanded language whose reduct is $\mathbf{A}$ and such that the interpretation of $\overline{a}$ in $\mathbf{A^c}$ is canonical in the sense that it is $a$. We remind that the language for residuated lattices already have constants 0 and 1 (without an overline), and so indeed it would have been enough to introduce only new constants for the elements $a \in A \setminus \{0, 1\}$. The logic associated with $\mathbf{A^c}$ will be denoted by $\mathbf{\Lambda(A^c)}$, and is obtained by putting, for all sets $\Gamma \cup \{\varphi\}$ of formulas (possibly including constants from $\{\bar{a} : a \in A\}$),

$$\Gamma \vdash_{\mathbf{A^c}} \varphi \quad \Longleftrightarrow \quad \forall h \in \text{Hom}(\mathbf{F}m, \mathbf{A^c}), \text{ if } h[\Gamma] \subseteq \{1\} \text{ then } h(\varphi) = 1. \quad (2.8)$$

This logic is introduced following the same pattern that is used for $\mathbf{\Lambda(A)}$, and hence $\mathbf{\Lambda(A^c)}$ is clearly a conservative expansion of $\mathbf{\Lambda(A)}$. Analogously to the case without canonical constants, the purpose of this article is to explain how to expand an axiomatization of $\mathbf{\Lambda(A^c)}$ into one of the corresponding modal logic with canonical constants. Hence, throughout this article we assume that there is a fixed axiomatization of $\mathbf{\Lambda(A^c)}$.

In Appendix A.2 the reader can find several remarks about the logics $\vdash_{\mathbf{A^c}}$. In particular, for finite residuated lattices with a unique coatom it is given (see Corollary A.16) a method to convert an axiomatization of $\mathbf{\Lambda(A)}$ into one of $\mathbf{\Lambda(A^c)}$, which in opinion of the authors is a new result in the literature.

## 3    Introducing the minimum modal logic over A

This section is devoted to introduce the definition of the minimum modal logic associated with a residuated lattice $\mathbf{A}$ (or $\mathbf{A^c}$), and to state several properties about this logic. In this definition it will be required that $\mathbf{A}$ is complete. Our interest is on minimum logics (i.e., logics given by the full class of frames, and so not by reflexive, transitive, . . . ones) because we think that this is the first step in the way to understand the full landscape of modal logics over a residuated lattice $\mathbf{A}$ (or $\mathbf{A^c}$).

**DEFINITION 3.1**
The *modal language* is the expansion of the non-modal one (see Section 2) by a new unary connective: the *necessity operator* $\Box$. ⊣

Thus, the modal language contains the propositional language generated by the set $Var$ of propositional variables together with the connectives $\wedge, \vee, \odot, \rightarrow, 1, 0$ and $\Box$. Depending whether we have or not the canonical constants in the non-modal basis, these constants will be included or not in the modal language. As usual in the classical setting, for every $n \in \omega$ we consider $\Box^n$ as an abbreviation for $\Box \stackrel{n}{\ldots} \Box$. The *modal depth* of a modal formula is defined as the number of nested necessity operators. Thus, non-modal formulas are exactly those formulas with modal depth 0. Every $\Box$ occurrence also has associated a *modal degree*, which is defined as the modal depth of the subformula preceded by this $\Box$ occurrence. For example, the formula $\Box(p \rightarrow \Box q) \rightarrow (\Box p \rightarrow \Box q)$ has modal depth 2, and the modal degree of each one of $\Box$ occurrences is, from left to right, 1, 0, 0 and 0.

**DEFINITION 3.2**
An $\mathbf{A}$-*valued Kripke frame* (or simply a *Kripke frame*) is a pair $\mathfrak{F} = \langle W, R \rangle$ where $W$ is a non empty set (whose elements are called *worlds*) and $R$ is a binary relation valued in $A$ (i.e., $R : W \times W \longrightarrow A$) called *accessibility relation*. The Kripke frame is said to be *crisp* or *classical*[2] in case that the range of $R$ is included in $\{0, 1\}$, and it is

---

[2] In this article we prefer to use the word "crisp" rather than "classical". The reason to avoid this meaning of "classical" is because we prefer to use "classical" to mean that $\mathbf{A}$ is the Boolean algebra $\mathbf{2}$ with two elements.

*idempotent* if the range of $R$ is included in the set of idempotent elements of $\mathbf{A}$ (i.e., $R[W \times W] \subseteq \{a \in A : a \odot a = a\}$). The classes of Kripke frames, crisp Kripke frames and idempotent Kripke frames will be denoted, respectively, by $\mathsf{Fr}(\mathbf{A})$, $\mathsf{CFr}(\mathbf{A})$ and $\mathsf{IFr}(\mathbf{A})$ (or simply $\mathsf{Fr}$, $\mathsf{CFr}$ and $\mathsf{IFr}$ if there is no ambiguity).    ⊣

**Definition 3.3**
An $\mathbf{A}$-*valued Kripke model* (or simply a *Kripke model*) is a 3-tuple $\mathfrak{M} = \langle W, R, V \rangle$ where $\langle W, R \rangle$ is an $\mathbf{A}$-valued Kripke frame and $V$ is a map, called *valuation*, assigning to each variable in $Var$ and each world in $W$ an element of $A$ (i.e., $V : Var \times W \longrightarrow A$). In such a case we say that $\mathfrak{M}$ is based on the Kripke frame $\langle W, R \rangle$.    ⊣

From now on we assume throughout the article that the residuated lattice $\mathbf{A}$ is complete, i.e., all suprema and infima exist. If $\mathfrak{M} = \langle W, R, V \rangle$ is a Kripke model, the map $V$ can be uniquely extended to a map, also denoted by $V$, assigning to each modal formula and each world in $W$ an element of $A$ (i.e., $V : Fm \times W \longrightarrow A$) satisfying that:

- $V$ is an algebraic homomorphism, in its first component, for the connectives $\wedge, \vee, \odot, \rightarrow, 1$ and $0$,

- if the modal language contains canonical constants $\overline{a}$'s, then $V(\overline{a}, w) = a$.

- $V(\Box \varphi, w) = \bigwedge \{ R(w, w') \rightarrow V(\varphi, w') : w' \in W \}$.

We emphasize that lattice completeness[3] is what allow us to be sure that we can compute the value of formulas with $\Box$. Although the original map and its extension are different mappings there is no problem, since one is an extension of the other, to use the same notation $V$ for both.

**Remark 3.4** (Classical Setting)
For the case that $\mathbf{A}$ is the Boolean algebra of two elements all previous definitions correspond to the standard terminology in the field of classical modal logic (cf. [11, 4, 5]), and hence this approach can be understood as a generalization of the classical modal case. As far as the authors know the first one to talk about this way of extending the valuation $V$ into the many-valued modal realm was Fitting in [22].    ⊣

**Remark 3.5** (The corresponding first-order fragment)
In the classical setting it is well known [4, Section 2.4] that the semantics of the modal language can be seen as a fragment of first order classical logic (where propositional variables are seen as predicates): every modal formula expresses the same than its corresponding first order translation (which has one free variable), e.g., $\Box p$ expresses the same than $\forall y (Rxy \rightarrow Py)$. This translation is often called in the literature the *standard translation*. It is worth pointing out that this very translation also preserves the meaning when the modal formula is considered under the many-valued modal semantics, and the first order formula is considered under the first order many-valued semantics. In other words, the above way to extend $V$ from variables to modal formulas is the only way to do it as far as we want that the standard translation preserves meaning in the many-valued setting. Of course, the previous standard translation is not the only one that does the job in the classical setting (e.g., we could

---

[3] Another possibility not explored in this article would be to consider arbitrary residuated lattices and to restrict modal valuations to "safe" ones in the sense that the infima needed to compute $\Box$ do exist (cf. [30]).

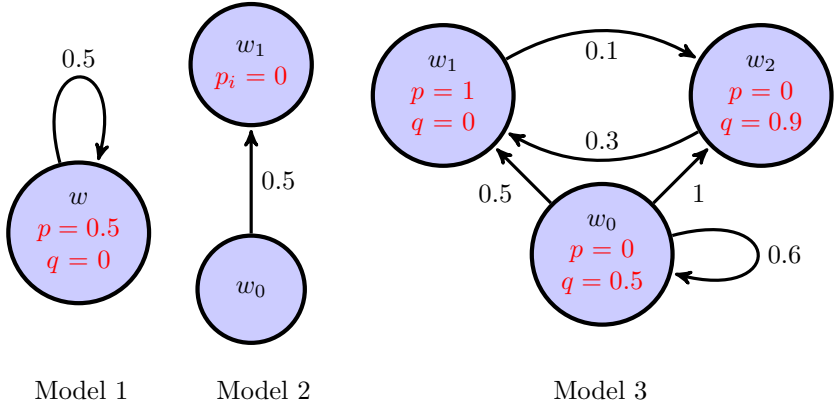

Fig. 1. Some examples of Kripke models

have translated $\Box p$ as $\forall y(\neg Rxy \lor Py)$) but the authors consider it the most natural one. ⊣

If we would like to introduce a *possibility operator* $\Diamond$, the right way to do it semantically (as far as we want that the standard translation preserves meaning) is to assume that $\Diamond$ is ruled by the condition

$$V(\Diamond \varphi, w) = \bigvee \{R(w, w') \odot V(\varphi, w') : w' \in W\}.$$

The reader can easily check that the following statements are equivalent:

- $V(\Diamond \varphi, w) = V(\neg \Box \neg \varphi, w)$ for every $V$, $\varphi$ and $w$,
- $\mathbf{A} \models x \odot \neg y \approx \neg(x \to y)$,
- $\mathbf{A} \models x \approx \neg \neg x$.

Thus, in general it is not convenient to consider $\Diamond$ as an abbreviation of $\neg \Box \neg$. The only case where this can be done without troubles is in the involutive case. We will not address the study of many-valued modal logics where both modalities $\Box$ and $\Diamond$ are present (except for the case that $\mathbf{A}$ is involutive, e.g. $\mathsf{MV}$ algebras, because then we get each modality from the other).

### Convention 3.6

In the classical framework it is usual to present Kripke frames in the form of diagrams by depicting worlds as points and drawing an arrow between each two worlds related by the accessibility relation. In the case of Kripke models what it is commonly done is writing inside (or next to) every point the variables that are valid in the world represented by this point. In the present many-valued framework we adopt a similar convention. Worlds are depicted as points (the name of the world will be written inside the point), and a labelled arrow connecting two worlds (perhaps they are just one) $w$ and $w'$ is drawn whenever $R(w, w') \neq 0$. The label of the arrow connecting $w$ and $w'$ is $R(w, w')$. These are the elements of a Kripke frame and we will use black color for them. The valuation will be given by writing inside each point, with red color, the value of the variables in this point. We shall often omit some variables in the description under the agreement that all variables not explicitly written have value 1. The reader can find the diagrams of several Kripke models in Figure 1.   ⊣

### 3.1   Some remarks about validity

**Definition 3.7**
We write $\mathfrak{M}, w \models^1 \varphi$ or simply $w \models^1 \varphi$, and say that $w$ *validates* $\varphi$, whenever $V(\varphi, w) = 1$. And we write $\mathfrak{M} \models^1 \varphi$ whenever $w \models^1 \varphi$ for every $w \in W$; then, we say $\varphi$ is *valid* in $\mathfrak{M}$. If $\mathfrak{F}$ is a frame, we say that $\varphi$ is *valid in the frame* $\mathfrak{F}$ when $\varphi$ is valid in any Kripke model based on $\mathfrak{F}$. We write it $\mathfrak{F} \models^1 \varphi$ for short. And if $\mathsf{K}$ is a class of frames then we write $\mathsf{K} \models^1 \varphi$ to mean that $\varphi$ is valid in all frames in this class.   ⊣

Similarly to $\models^1$ we could have introduced the relation $\models^+$ of being *positively valid* requiring that the valuation does not get the value 0 (i.e., it gets a positive value). In the classical modal setting there is a natural notion of satisfiability that is dual to validity, but in the many-valued modal setting we must be really careful (even in the involutive case) about this dual notion. Dually to the previous two notions we can define the notions of a formula being *satisfied* or *positively satisfied* in a Kripke model or frame. If the negation of **A** is involutive then it obviously holds that

- $\varphi$ is valid in a model $\mathfrak{M}$   iff   $\neg\varphi$ is not positively satisfied in $\mathfrak{M}$,
- $\varphi$ is positively valid in a model $\mathfrak{M}$   iff   $\neg\varphi$ is not satisfied in $\mathfrak{M}$.

Hence, the dual of validity, even in the involutive case, is not satisfiability; the dual is positive satisfiability. In the classical modal setting it really holds that the dual of validity is satisfiability, but this is only because in this particular case satisfiability coincides with positive satisfiability. And of course, without an involutive negation there is no connection between the previous notions.

**Definition 3.8**
A *semantic modal valuation*[4] is a map $v$ from the set of modal formulas into $A$ (i.e., $v : Fm \longrightarrow A$) such that $v = V(\bullet, w)$ for some Kripke model $\mathfrak{M} = \langle W, R, V \rangle$ and some world $w$ in $W$. Sometimes, if there is no ambiguity, we will use the same name $w$ to denote $V(\bullet, w)$.   ⊣

**Remark 3.9** (Properties on semantic modal valuations)
It is clear that semantic modal valuations are homomorphisms for all non-modal connectives, i.e., for $\wedge, \vee, \odot, \rightarrow, 1$ and $0$ (and the canonical constants when they are in the language); but this is not the case for $\Box$. Thus, for every non-modal formula $\varphi(p_1, \ldots, p_n)$ and every semantic modal valuation $v$ it holds that $v(\varphi) = \varphi^{\mathbf{A}}(v(p_1), \ldots, v(p_n))$. Another property that is clearly true is that if $\mathsf{Fr} \models^1 \varphi$, then $v(\varphi) = 1$ for every semantic modal valuation $v$.   ⊣

Next, and as a first approximation to modal logics, we discuss the validity of some formula schemata in the classes $\mathsf{Fr} \supseteq \mathsf{IFr} \supseteq \mathsf{CFr}$. Although they are indeed formula schema we will simply use the word "formula" (the reader must think that $\varphi, \psi, \ldots$ are somehow metavariables) for the rest of the article.

---

[4] We prefer to keep the adjective "semantic" because it is suggesting that these valuations are arising from Kripke models. We use this word to emphasize the difference between these semantic modal valuations and the points of the canonical models defined in Sections 4 and 5: these points may be seen as "syntactic" modal valuations.

The first remarkable fact is that the famous *normality* axiom ($\mathsf{K}$)

$$\Box(\varphi \to \psi) \to (\Box\varphi \to \Box\psi) \tag{K}$$

is, in general, not valid[5] in $\mathsf{Fr}$. A simple counterexample is to consider the algebra $\mathbf{A}$ as the $\mathsf{MV}$ chain with three elements and the Kripke frame $\mathfrak{F}$ as the first one depicted in Figure 1; then the semantic modal valuation $w$ satisfies that $w(\Box(p \to q)) = 1 = w(\Box p)$ and that $w(\Box q) = 0.5$. The same frame $\mathfrak{F}$ also shows that $(\Box\varphi \odot \Box\psi) \to \Box(\varphi \odot \psi)$ is not valid in general (think on the case that $\varphi = p$ and $\psi = p \to q$). In Appendix B, the method of the non-modal companion is introduced providing a really simple way to discard the validity of some modal formulas, among them ($\mathsf{K}$).

In the following proposition we give positive results stating the validity of some formulas. We will refer to the formula ($\mathsf{MD}$) as the *meet distributivity* axiom, and we remark that the formula (3.1) is one of the formulas considered by CAICEDO and RODRÍGUEZ in [10] (see Table 2).

PROPOSITION 3.10

1. Some valid formulas in $\mathsf{Fr}$ are

$$(\Box\varphi \wedge \Box\psi) \leftrightarrow \Box(\varphi \wedge \psi) \tag{MD}$$

$$\neg\neg\Box\varphi \to \Box\neg\neg\varphi. \tag{3.1}$$

$$\Box(\overline{a} \to \varphi) \leftrightarrow (\overline{a} \to \Box\varphi) \tag{Ax$_a$}$$

2. Some valid formulas in $\mathsf{IFr}$ are

$$\Box(\varphi \to \psi) \to (\Box\varphi \to \Box\psi) \tag{K}$$

$$(\Box\varphi \odot \Box\psi) \to \Box(\varphi \odot \psi) \tag{3.2}$$

3. Some valid formulas in $\mathsf{CFr}$ are

$$\Box\overline{0} \vee \neg\Box\overline{0} \tag{3.3}$$

$$\Box\overline{a} \vee (\Box\overline{a} \leftrightarrow \overline{a}) \tag{3.4}$$

$$\Box\overline{0} \vee (\Box\overline{a} \leftrightarrow \overline{a}) \tag{3.5}$$

PROOF. 1): The validity of ($\mathsf{MD}$) is an straightforward consequence of the fact that $x \to (y \wedge z) \approx (x \to y) \wedge (x \to z)$ holds in any residuated lattice, and the validity of ($\mathsf{Ax}_a$) easily follows from (2.3). Next we prove that $\mathsf{Fr} \models^1 \neg\neg\Box\varphi \to \Box\neg\neg\varphi$. Let us take any frame $\mathfrak{F}$ and a world $w$ in $\mathfrak{F}$, and let $V$ be a valuation on $\mathfrak{F}$. The proof is completed by showing that $V(\neg\neg\Box\varphi, w) \leqslant V(\Box\neg\neg\varphi, w)$. It suffices to prove that $V(\neg\neg\Box\varphi, w) \leqslant R(w, w') \to V(\neg\neg\varphi, w')$ for every world $w'$. Hence, let su consider an arbitrary world $w'$. We know that $R(w, w') \odot V(\neg\varphi, w') \odot V(\Box\varphi, w) \leqslant R(w, w') \odot V(\neg\varphi, w') \odot (R(w, w') \to V(\varphi, w')) \leqslant V(\neg\varphi, w') \odot V(\varphi, w') \leqslant 0$. We thus get $R(w, w') \odot V(\neg\varphi, w') \leqslant V(\neg\Box\varphi, w)$. This gives $R(w, w') \odot V(\neg\neg\Box\varphi, w) \leqslant V(\neg\neg\varphi, w')$; and hence $V(\neg\neg\Box\varphi, w) \leqslant R(w, w') \to V(\neg\neg\varphi, w')$.

---

[5] In the case that $\mathbf{A}$ is an $\mathsf{MTL}$ algebra it is easy to see that $\mathsf{Fr} \models^1 (\Box((\varphi \to \psi)^2) \to (\Box(\varphi^2) \to \Box\psi)$, where this formula resembles ($\mathsf{K}$). The reason why this formula is valid is essentially the fact that $\mathbf{A} \models (x^2 \wedge (x \to y)^2) \to y \approx 1$. Indeed, $\mathsf{Fr} \models^1 (\Box(\varphi^2) \wedge (\Box((\varphi \to \psi)^2)) \to \Box\psi$; this last formula being stronger than $\Box((\varphi \to \psi)^2) \to (\Box(\varphi^2) \to \Box\psi)$.

2): We only give the proof for the case of the normality axiom. The proof is a consequence of the validity of the quasi-equation

$$x \approx x \odot x \quad \Rightarrow \quad (x \to (y \to z)) \odot (x \to y) \leqslant x \to z$$

in all residuated lattices. To see $\mathsf{IFr} \models^1 (\mathsf{K})$ take any frame $\mathfrak{F}$ and a world $w$ in $\mathfrak{F}$, and let $V$ be a valuation on $\mathfrak{F}$. Then, for every world $w'$ it holds that $V(\Box(\varphi \to \psi), w) \odot V(\Box\varphi, w) \leqslant (R(w, w') \to V(\varphi \to \psi, w')) \odot (R(w, w') \to V(\varphi, w')) \leqslant R(w, w') \to V(\psi, w')$. Thus, $V(\Box(\varphi \to \psi), w) \odot V(\Box\varphi, w) \leqslant \bigwedge \{R(w, w') \to V(\psi, w') : w' \in W\} = V(\Box\psi, w)$. Hence, $V(\Box(\varphi \to \psi) \to (\Box\varphi \to \Box\psi), w) = 1$.

3): The validity of the first two formulas follows from the fact that for every valuation $V$ in a crisp frame and every world $w$ in this frame it holds that $V(\Box\bar{a}, w) \in \{a, 1\}$. Next we prove that $\mathsf{CFr} \models^1 \Box 0 \lor (\Box\bar{a} \leftrightarrow \bar{a})$. Let us take any valuation $V$ in a crisp Kripke frame and a world $w$ in this frame such that $V(\Box 0, w) \neq 1$. It is enough to prove that $V(\bar{a}, w) = V(\Box\bar{a}, w)$. Since $V(\Box 0, w) \neq 1$ we get that there is a world $w'$ such that $1 \neq R(w, w') \to V(0, w')$. As the frame is crisp we have that $R(w, w') = 1$. Therefore, $V(\bar{a}, w) = a = 1 \to a = \bigvee \{R(w, w'') : w'' \in W\} \to a = \bigwedge \{R(w, w'') \to a : w'' \in W\} = V(\Box\bar{a}, w)$. ∎

Before the previous proposition we saw that the formulas in the second item are not in general valid in $\mathsf{Fr}$. Analogously, the reader can easily check that the formulas in the third item are in general not valid in $\mathsf{IFr}$ (cf. Proposition 3.14). In the next proposition we use the notion of *(infinitely) distributive element* $a \in A$ defined as those elements satisfying that

$$\bigwedge \{a \lor x : x \in X\} \approx a \lor \bigwedge X \tag{3.6}$$

for every set $X \subseteq A$. It is obvious that if $\mathbf{A}$ is finite then the previous definition can be relaxed to the case that $X = \{x_1, x_2\}$. It is interesting to point out that if $\mathbf{A}$ is a finite algebra, then

- all Boolean elements[6] are distributive,
- all coatoms of $\mathbf{A}$ are distributive[7].

The first statement is a consequence of the fact that $p \lor \neg p \vdash_{\mathsf{RL}} (p \lor (q \land r)) \leftrightarrow ((p \lor q) \land (p \lor r))$. Here $\vdash_{\mathsf{RL}}$ refers to the non-modal logic associated with the class of all residuated lattices (see [25]). And the second one follows from the fact that if $a$ is a coatom, then the only possibility that (3.6) fails is when $a = a \lor (x_1 \land x_2)$ and $1 = (a \lor x_1) \land (a \lor x_2)$; which is a contradiction with the fact that $p \lor q, p \lor r \vdash_{\mathsf{RL}} p \lor (q \land r)$. In the proofs of these two statements it has been crucial the fact that $\vdash_{\mathsf{RL}}$ admits proofs by cases. Hence, in particular we know that Proposition 3.11 holds when $a$ is a coatom of a finite $\mathbf{A}$.

---

[6] We remind the reader that the set of Boolean elements is the set of $a$'s such that there is some $b \in A$ satisfying that $a \lor b = 1$ and $a \land b = 0$. In the setting of residuated lattices it is well known that the set of Boolean elements is exactly $\{a \in A : a \lor \neg a = 1\}$. Another characterization of this set is given by $\{a \in A : a \to b = \neg a \lor b$ for every $b \in A\}$. Indeed, using the fact that $\vdash_{\mathsf{RL}}$ admits proofs by cases it is easy to see that for every non-modal formula $\varphi(p, \bar{q})$, it holds that

– $\varphi^{\mathbf{A}}(p, \bar{q}) = 1$ for any residuated lattice $\mathbf{A}$ and any assignation such that $p^{\mathbf{A}} \in \{0, 1\}$, iff
– $\varphi^{\mathbf{A}}(p, \bar{q}) = 1$ for any residuated lattice $\mathbf{A}$ and any assignation such that $p^{\mathbf{A}}$ is a Boolean element.

[7] This is somehow suggesting the "picture" that the upper part of a residuated lattice is distributive in the previous sense.

**Proposition 3.11**
Let $a$ be a distributive element of $\mathbf{A}$. Then,

$$\Box(\overline{a} \vee \varphi) \to (\overline{a} \vee \Box\varphi). \tag{3.7}$$

is valid in $\mathsf{CFr}$.

PROOF. This is easily proved from the fact that $a$ is a distributive element. ∎

Next we consider the problem of definability for the main classes of frames considered in this article. First of all we prove that the class $\mathsf{IFr}$ of idempotent frames is modally definable by the normality axiom $(\mathsf{K})$. Indeed, the formula (3.2) and the formula $(\Box\varphi \odot \Box\varphi) \to \Box(\varphi \odot \varphi)$ are also defining the same class $\mathsf{IFr}$. We stress the fact that none of the formulas in Proposition 3.12 is using canonical constants.

**Proposition 3.12**
$\mathsf{IFr} = \{\mathfrak{F} : \mathfrak{F} \models^1 \Box(p \to q) \to (\Box p \to \Box q)\} = \{\mathfrak{F} : \mathfrak{F} \models^1 (\Box p \odot \Box q) \to \Box(p \odot q)\} = \{\mathfrak{F} : \mathfrak{F} \models^1 (\Box p \odot \Box p) \to \Box(p \odot p)\}$.

PROOF. The inclusion of $\mathsf{IFr}$ in the other three classes follow from Proposition 3.10. The rest of inclusions can be proved using the same idea, and so we will prove only one of them, namely that $\{\mathfrak{F} : \mathfrak{F} \models^1 (\Box p \odot \Box p) \to \Box(p \odot p)\} \subseteq \mathsf{IFr}$. Hence, let us consider a Kripke frame $\mathfrak{F} = \langle W, R \rangle$ and let us assume that there is an element $a \in A$ such that $a \odot a < a$ and $a = R(w, w')$ for certain worlds $w, w'$. We define a valuation $V$ by the conditions: (i) $V(p, w') = a$, and (ii) $V(p, w'') = 1$ for every world $w'' \neq w'$. Then, $V(\Box p, w) = 1$ and $V(\Box(p \odot p), w) \leqslant a \to a \odot a \neq 1$. Therefore, $(\Box p \odot \Box p) \to \Box(p \odot p)$ is not valid in $\mathfrak{F}$. ∎

**Corollary 3.13**
- Axiom $(\mathsf{K})$ is valid in $\mathsf{Fr}$    iff    $\mathbf{A}$ is a Heyting algebra    iff    $\mathsf{Fr} = \mathsf{IFr}$.
- $\{\varphi \in Fm : \mathsf{Fr} \models^1 \varphi\} = \{\varphi \in Fm : \mathsf{IFr} \models^1 \varphi\}$    iff    $\mathsf{Fr} = \mathsf{IFr}$.

PROOF. It is obvious that $\mathbf{A}$ is a Heyting algebra iff $\mathsf{Fr} = \mathsf{IFr}$. Then, by the previous proposition we obtain this corollary. ∎

Next we address the issue whether the class $\mathsf{CFr}$ of crisp Kripke frames is characterized by the set of its valid modal formulas. In general the answer is negative due to a result proved in [10]: over the standard Gödel algebra[8], the valid modal formulas (without canonical constants) in $\mathsf{IFr}$ are exactly the valid ones in $\mathsf{CFr}$. Hence, in general it may happen that $\mathsf{IFr} \neq \mathsf{CFr}$ while both classes share the same valid formulas. Although we have not got a characterization of the algebras $\mathbf{A}$ such that $\{\varphi \in Fm : \mathsf{IFr} \models^1 \varphi\} = \{\varphi \in Fm : \mathsf{CFr} \models^1 \varphi\}$ we can state some necessary conditions. The last item also gives us a characterization in a particular case.

**Proposition 3.14**
1. $\mathsf{IFr} \models^1 \Box 0 \vee \neg\Box 0$    iff    $\mathbf{A} \models x \approx x \odot x \Rightarrow \neg x \vee \neg\neg x \approx 1$.
2. If $\{\varphi \in Fm : \mathsf{IFr} \models^1 \varphi\} = \{\varphi \in Fm : \mathsf{CFr} \models^1 \varphi\}$, then $\mathbf{A} \models x \approx x \odot x \Rightarrow \neg x \vee \neg\neg x \approx 1$.
3. If there are canonical constants in the language and 1 is join irreducible, then $\{\varphi \in Fm : \mathsf{IFr} \models^1 \varphi\} = \{\varphi \in Fm : \mathsf{CFr} \models^1 \varphi\}$ iff $\mathsf{IFr} = \mathsf{CFr}$.

---

[8]A moment of reflection shows that this also holds for finite Gödel chains.

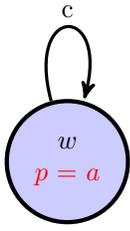

| $\odot$ | 0 | $a$ | $b$ | $c$ | $d$ | 1 |
|---|---|---|---|---|---|---|
| 0 | 0 | 0 | 0 | 0 | 0 | 0 |
| $a$ | 0 | $a$ | $a$ | $a$ | $a$ | $a$ |
| $b$ | 0 | $a$ | $a$ | $a$ | $b$ | $b$ |
| $c$ | 0 | $a$ | $a$ | $c$ | $c$ | $c$ |
| $d$ | 0 | $a$ | $b$ | $c$ | $d$ | $d$ |
| 1 | 0 | $a$ | $b$ | $c$ | $d$ | 1 |

| $\rightarrow$ | 0 | $a$ | $b$ | $c$ | $d$ | 1 |
|---|---|---|---|---|---|---|
| 0 | 1 | 1 | 1 | 1 | 1 | 1 |
| $a$ | 0 | 1 | 1 | 1 | 1 | 1 |
| $b$ | 0 | $c$ | 1 | 1 | 1 | 1 |
| $c$ | 0 | $b$ | $b$ | 1 | 1 | 1 |
| $d$ | 0 | $a$ | $b$ | $c$ | 1 | 1 |
| 1 | 0 | $a$ | $b$ | $c$ | $d$ | 1 |

An interesting MTL chain         A Kripke Model

Fɪɢ. 2. A counterexample to $\{\varphi : \mathsf{CFr} \models^1 \varphi\} = \{\varphi : \mathsf{IFr} \models^1 \varphi\}$

Proof. 1): Firstly let us assume that $\mathsf{IFr}(\mathbf{A}) \models^1 \Box 0 \vee \neg\Box 0$ and that $a \in A$ is an idempotent element. The proof finishes by considering the Kripke frame $\mathfrak{F} \in \mathsf{IFr}$ with only one point $w$ and such that $R(w, w) = a$.

To see the other direction we take any frame $\mathfrak{F}$ in $\mathsf{IFr}$ and a world $w$ in $\mathfrak{F}$, and let $V$ be a valuation on $\mathfrak{F}$. Then, $V(\Box 0, w) = \bigwedge\{R(w, w') \rightarrow V(0, w') : w' \in W\} = \bigwedge\{\neg R(w, w') : w' \in W\} = \neg\bigvee\{R(w, w') : w' \in W\}$. The element $a := \bigvee\{R(w, w') : w' \in W\}$ is idempotent because it is the join of a family of idempotent elements. We get that $V(\Box 0 \vee \neg\Box 0, w) = \neg a \vee \neg\neg a$, and thus it is 1 by assumption.

2): This is trivial from the previous item together with Proposition 3.10.

3): If $\mathsf{IFr} \neq \mathsf{CFr}$, then there is an idempotent element $a \in A \setminus \{0, 1\}$. By Proposition 3.10 it is enough to prove that $\mathsf{IFr} \not\models^1 \Box 0 \vee (\overline{a} \leftrightarrow \Box\overline{a})$; simply consider the Kripke frame $\mathfrak{F} \in \mathsf{IFr}$ with only one point $w$ and such that $R(w, w) = a$.    ∎

Next example shows that the second statement in the previous proposition is not an equivalence. In other words, it is not enough to consider the Stonean condition restricted to idempotent elements[9] in order to solve the previous open problem (even for the case without canonical constants), i.e., the condition $\mathsf{IFr} \models^1 \Box 0 \vee \neg\Box 0$ is not strong enough to guarantee that $\{\varphi : \mathsf{IFr} \models^1 \varphi\} = \{\varphi : \mathsf{CFr} \models^1 \varphi\}$.

Example 3.15
Let us consider $\mathbf{A}$ as the MTL chain given in Figure 2. It is clear that $\mathbf{A}$ satisfies $x \approx x \odot x \Rightarrow \neg x \vee \neg\neg x \approx 1$, i.e., $\mathsf{IFr} \models^1 \Box 0 \vee \neg\Box 0$. On the other hand, the reader can easily prove[10] that $\mathsf{CFr} \models^1 \Box(p \odot p) \leftrightarrow (\Box p \odot \Box p)$. Finally, the Kripke model depicted in Figure 2 witnesses that $\mathsf{IFr} \not\models^1 \Box(p \odot p) \leftrightarrow (\Box p \odot \Box p)$, and so $\{\varphi \in Fm : \mathsf{IFr} \models^1 \varphi\} \neq \{\varphi \in Fm : \mathsf{CFr} \models^1 \varphi\}$.    ⊣

In order to understand the set of modal formulas valid in the class of crisp frames it seems natural to introduce a new class of Kripke frames located between $\mathsf{IFr}$ and $\mathsf{CFr}$. Analogously to the classes $\mathsf{Fr}$, $\mathsf{IFr}$ and $\mathsf{CFr}$ we could have considered the class $\mathsf{BFr}$ of *Boolean frames* defined as the class of Kripke frames where the range of the accessibility relation is the set of Boolean elements (see footnote 6). Next we prove that the class $\mathsf{BFr}$ has the same valid formulas than the class $\mathsf{CFr}$ (even when canonical constants are allowed), and we also show that this new class $\mathsf{BFr}$ is always modally definable, using the formulas in Proposition 3.11, when there are canonical constants in the language. It is worth pointing out that if 1 is join irreducible then $\mathsf{BFr} = \mathsf{CFr}$, and hence without canonical constants in the language in general it is not possible to characterize $\mathsf{BFr}$ using modal formulas since the same counterexample cited above (and based on Gödel chains) applies here.

---

[9] We notice that this property holds in any BL chain.

[10] Indeed, for every finite chain $\mathbf{A}$ it holds that $\mathsf{CFr} \models^1 \Box(p \odot p) \leftrightarrow (\Box p \odot \Box p)$. However, it is possible to find infinite BL chains where $\Box(p \odot p) \leftrightarrow (\Box p \odot \Box p)$ is not valid; e.g., take $\mathbf{A}$ as the (complete) subalgebra of $\mathbf{Ł_3} \oplus [\mathbf{0}, \mathbf{1}]_{\mathbf{G}}$ obtained by deleting the idempotent element gluing both parts (cf. [19, p. 281]). It is worth pointing out that the recent result [37] characterizes the BL chains where $\Box(p \odot p) \leftrightarrow (\Box p \odot \Box p)$ is valid in $\mathsf{CFr}$ as the BL chains satisfying all standard first-order BL tautologies.

**Theorem 3.16**
Let $\mathbf{A}$ be a finite residuated lattice. Then, $\mathsf{BFr}$ and $\mathsf{CFr}$ have the same valid modal formulas (even allowing canonical constants in the language).

**Proof.** First of all we remark that by the finiteness assumption we know that $\mathbf{A} \cong \prod_{i<n} \mathbf{A}/\theta_i$ where the only Boolean elements of $\mathbf{A}/\theta_i$ are $\{0, 1\}$. This is a consequence of the characterization of directly indecomposable residuated lattices stated in [43, Proposition 1.5].

Since $\mathsf{CFr} \subseteq \mathsf{BFr}$ it suffices to prove that all modal formulas valid in $\mathsf{CFr}$ are also valid in $\mathsf{BFr}$. The proof is based on a method that converts a Boolean Kripke model over $\mathbf{A}$ into a family of $n$ crisp Kripke models $\mathbf{A}$. Let us assume that $\langle W, R, V \rangle$ is a Boolean Kripke model. Then, for every $i < n$ we consider the crisp Kripke model $\langle W_i, R_i, V_i \rangle$ defined by

- $W_i = W$,
- $R_i(w, w')$ is (i) $1^{\mathbf{A}}$ if $R(w, w') \in 1/\theta_1$, and (ii) $0^{\mathbf{A}}$ if $R(w, w') \in 0/\theta_1$,
- $V_i(p, w) = V(p, w)$.

The previous Kripke models are well defined because since $R(w, w')$ is a Boolean element we know that for every $i < n$, either $R(w, w') \in 1/\theta_1$ or $R(w, w') \in 0/\theta_1$. By construction of these Kripke models it is straightforward to check that for every modal formula $\varphi$ (perhaps including canonical constants), every world $w \in W$ and every $i < n$ it holds that

$$\pi_i(V(\varphi, w)) = \pi_i(V_i(\varphi, w)),$$

where $\pi_i$ is the $i$-th projection. Therefore, $V(\varphi, w) = \langle \pi_i(V_i(\varphi, w)) : i < n \rangle$ for every modal formula $\varphi$ and every world $w \in W$. Using this fact it is obvious all modal formulas valid in $\mathsf{CFr}$ are also valid in $\mathsf{BFr}$. ∎

**Proposition 3.17**
Let $\mathbf{A}$ be a finite residuated lattice, and let $\mathfrak{F}$ be a frame. The following statements are equivalent.

1. $\mathfrak{F}$ is a Boolean frame,
2. $\mathfrak{F} \models^1 \{\Box(\overline{a} \vee p) \to (\overline{a} \vee \Box p) : a$ is a distributive element of $\mathbf{A}\}$,
3. $\mathfrak{F} \models^1 \{\Box(\overline{k} \vee p) \to (\overline{k} \vee \Box p) : k$ is a coatom of $\mathbf{A}\}$.

**Proof.** $1 \Rightarrow 2$): This follows from Theorem 3.16 together with Proposition 3.11.

$2 \Rightarrow 3$): This is trivial since we know that all coatoms are distributive elements.

$3 \Rightarrow 1$): Let us assume that there is an element $a \in A$ such that $a$ is not a Boolean element and $a = R(w, w')$ for certain worlds $w, w'$. Hence, $a \vee \neg a \neq 1$, and so $a \vee \neg a \leqslant k$ for certain coatom $k$. We define a valuation $V$ by the conditions: (i) $V(p, w') = 0$, and (ii) $V(p, w'') = 1$ for every world $w'' \neq w'$. Then, $V(\Box(\bar{k} \vee p), w) = a \to (k \vee 0) = 1$ and $V(\bar{k} \vee \Box p, w) = k \vee (a \to 0) = k$. Therefore, $\Box(\bar{k} \vee p) \to (\bar{k} \vee \Box p)$ is not valid in $\mathfrak{F}$, which is a contradiction. ∎

It is worth pointing out that Proposition B.4 (and Corollary B.5) gives us another method to get valid formulas. We have left this result in Appendix B because although this proposition is formulated without appealing to the non-modal companion we can understand this result as explaining why for some modal formulas the non-modal companion method is enough to discard its validity. To finish this section we prove a similar result to Proposition B.4, but this time validity is restricted to a particularly interesting subclass of Kripke models, the class of modally witnessed models.

A Kripke model $\langle W, R, V \rangle$ is *modally witnessed*[11] when for every modal formula $\varphi$ and every world $w$, there is a world $w'$ (called the *witness*) such that

$$V(\Box\varphi, w) = R(w, w') \to V(\varphi, w').$$

It is obvious that this notion is the relaxed version of the witnessed first-order structures considered by Hájek in [32, 33] when replacing first order formulas with modal (propositional) ones. In case we restrict our attention to valid formulas in all modally witnessed Kripke models then we can get a similar result to Proposition B.4 without the requirement saying that $\varepsilon$ is expanding.

PROPOSITION 3.18
Let $\delta(p_1, \ldots, p_n)$ and $\varepsilon(p)$ be two non-modal formulas such that $\delta$ is non decreasing over $\mathbf{A^c}$ in any of its arguments. And let $\varphi_1, \ldots, \varphi_n, \varphi$ be non-modal formulas. If $\vdash_{\mathbf{A^c}} \delta(r \to \varphi_1, \ldots, r \to \varphi_n) \to \varepsilon(r \to \varphi)$ being $r$ a variable not appearing in $\{\varphi_1, \ldots, \varphi_n, \varphi\}$, then $\delta(\Box\varphi_1, \ldots, \Box\varphi_n) \to \varepsilon(\Box\varphi)$ is valid in all modally witnessed Kripke models.

PROOF. Take any modally witnessed Kripke model $\langle W, R, V \rangle$ and a world $w$ in $W$. We consider $w'$ as the witness for $V(\Box\varphi, w)$, i.e., $V(\Box\varphi, w) = R(w, w') \to V(\varphi, w')$. Then, $V(\delta(\Box\varphi_1, \ldots, \Box\varphi_n), w) = \delta^{\mathbf{A^c}}(V(\Box\varphi_1, w), \ldots, V(\Box\varphi_n, w)) \leqslant \delta^{\mathbf{A^c}}(R(w, w') \to V(\varphi_1, w'), \ldots, R(w, w') \to V(\varphi_n, w')) \leqslant \varepsilon^{\mathbf{A^c}}(R(w, w') \to V(\varphi, w')) = \varepsilon^{\mathbf{A^c}}(V(\Box\varphi, w)) = V(\varepsilon(\Box\varphi), w)$. ∎

## 3.2   *The local and global modal logics*

Now it is time to introduce the minimum modal logic we want to study in the present article. Our main interest is in the class $\mathsf{Fr}$ of all frames, but in order to be as general as possible the definitions are given for an arbitrary class $\mathsf{K}$ of frames. Thus, this definition also defines the modal logic associated with the other two basic classes of Kripke frames.

DEFINITION 3.19
In the following we assume $\mathbf{A}$ to be a given residuated lattice. Let $\mathsf{K}$ be a subclass of its Kripke frames. The set $\{\varphi \in Fm : \mathsf{K} \models^1 \varphi\}$ will be denoted by $\Lambda(\mathsf{K}, \mathbf{A})$. In

---

[11] We point out that in general finite Kripke models are not modally witnessed; but when $\mathbf{A}$ is a chain then all finite Krikpe models are modally witnessed. In the case that $\mathbf{A}$ is a chain it is worth pointing out a result of Hájek [31, Theorem 4] showing that there is a tight connection, as far as we only consider the modal language, between modally witnessed models and finite models.

case there are canonical constants we will use the notation $\Lambda(\mathsf{K}, \mathbf{A^c})$ to stress their presence. The *local (many-valued) modal logic* $\boldsymbol{\Lambda}(l, \mathsf{K}, \mathbf{A})$ is the logic obtained by defining, for all sets $\Gamma \cup \{\varphi\}$ of modal formulas (without canonical constants),

- $\Gamma \vdash_{\boldsymbol{\Lambda}(l, \mathsf{K}, \mathbf{A})} \varphi$ iff[12]
- for every semantic modal valuation $v$ arising from a Kripke frame in $\mathsf{K}$, if $v[\Gamma] \subseteq \{1\}$ then $v(\varphi) = 1$.

We will also use the notations $\boldsymbol{\Lambda}(l, \mathsf{K})$ or $\vdash^l_{\mathsf{K}(\mathbf{A})}$. And for the analogous definition, but including canonical constants, we will use the notation $\boldsymbol{\Lambda}(l, \mathsf{K}, \mathbf{A^c})$ or $\vdash^l_{\mathsf{K}(\mathbf{A^c})}$. On the other hand, the *global (many-valued) modal logic* $\boldsymbol{\Lambda}(g, \mathsf{K}, \mathbf{A})$ is the logic defined by

- $\Gamma \vdash_{\boldsymbol{\Lambda}(g, \mathsf{K}, \mathbf{A})} \varphi$ iff
- for every Kripke model $\mathfrak{M}$ arising from a Kripke frame in $\mathsf{K}$, if $\mathfrak{M} \models^1 \gamma$ for every $\gamma \in \Gamma$, then $\mathfrak{M} \models^1 \varphi$.

For the sake of simplicity we will write it simply $\boldsymbol{\Lambda}(g, \mathsf{K})$ or $\vdash^g_{\mathsf{K}(\mathbf{A})}$. If we follow the same definition scheme, but allowing canonical constants in the language, then we will use the notation $\boldsymbol{\Lambda}(g, \mathsf{K}, \mathbf{A^c})$ or $\vdash^g_{\mathsf{K}(\mathbf{A^c})}$. ⊣

By definition it is obvious that both logics, the local and the global, share the same set of theorems, namely $\Lambda(\mathsf{K}, \mathbf{A})$. Another obvious consequence of the definition is that $\boldsymbol{\Lambda}(l, \mathsf{K}, \mathbf{A}) \leqslant \boldsymbol{\Lambda}(g, \mathsf{K}, \mathbf{A})$, these two logics being conservative expansions of $\boldsymbol{\Lambda}(\mathbf{A})$. These last two sentences are also true with canonical constants in the language.

**Remark 3.20** (Global vs. Local)
The main advantage of the local modal logic is that in most common cases we can prove completeness using the canonical frame method (see Theorem 4.11). On the other hand, the global modal logic has other benefits. The first one is that $\boldsymbol{\Lambda}(g, \mathsf{K}, \mathbf{A})$ and $\boldsymbol{\Lambda}(g, \mathsf{K}, \mathbf{A^c})$ are always algebraizable in the Abstract Algebraic Logic framework (see [16, 24]). This is an easy consequence of the general theory using that $\varphi \to \psi \vdash^g_{\mathsf{K}} \Box\varphi \to \Box\psi$ together with the fact that the non-modal logics $\boldsymbol{\Lambda}(\mathbf{A})$ and $\boldsymbol{\Lambda}(\mathbf{A^c})$ are algebraizable with equivalence formulas $\{p \to q, q \to p\}$. Therefore, we can give an algebraic semantics (and hence, a truth functional semantics) for the global modal logic. Another benefit is that $\boldsymbol{\Lambda}(g, \mathsf{Fr}, \mathbf{A})$ is the fragment (given by the standard translation) of the first order logic associated with $\mathbf{A}$. Thus, we can transfer known results from first order logics to global modal logics. ⊣

In our opinion the developing of a many-valued modal logics hierarchy (adding properties like reflexivity, transitivity, etc. and following the same ideas than in the classical setting) cannot be successfully undertaken until we know how to axiomatize (whenever they are axiomatizable) the minimum modal logics $\boldsymbol{\Lambda}(l, \mathsf{Fr}, \mathbf{A})$ and $\boldsymbol{\Lambda}(g, \mathsf{Fr}, \mathbf{A})$. *What axioms and rules must be added to an axiomatization of* $\boldsymbol{\Lambda}(\mathbf{A})$ *to get an axiomatization of* $\boldsymbol{\Lambda}(l, \mathsf{Fr}, \mathbf{A})$*? And for* $\boldsymbol{\Lambda}(g, \mathsf{Fr}, \mathbf{A})$*?* Of course a first step to answer this question is to know how to obtain the set $\Lambda(\mathsf{Fr}, \mathbf{A})$ from an axiomatization of $\mathbf{\Lambda(A)}$.

---

[12] In the classical modal literature, it is not common to introduce the definition of the local modal logic using the notion of semantic modal valuation; alternative equivalent definitions are usually provided. However, we feel that this definition based on semantic modal valuations has the advantage of making more clear how to reduce the local modal logic to the non-modal logic (cf. Theorem 3.26). Comparing this definition of the local modal logic with the definition (2.7) of the non-modal logic it is obvious that in order to find this reduction one has to be able to characterize the semantic modal valuations as the non-modal homomorphisms that are sending a certain set of axioms to 1 (cf. Lemma 3.25, Theorem 4.11, etc).

This problem remains open in this general formulation (see Section 6 for a more detailed discussion). But in Sections 4 and 5 we solve some instances of the problem. We point out that, as far as the authors are aware, the results in these two sections present the first known axiomatizations for local many-valued modal logics where normality axiom ($\mathsf{K}$) fails.

First of all, we consider the problem of comparing the many-valued modal logics given by the main three classes of frames considered.

**PROPOSITION 3.21**

1. $\mathbf{\Lambda}(g, \mathsf{Fr}, \mathbf{A}) = \mathbf{\Lambda}(g, \mathsf{lFr}, \mathbf{A})$ iff $\mathbf{\Lambda}(g, \mathsf{Fr}, \mathbf{A^c}) = \mathbf{\Lambda}(g, \mathsf{lFr}, \mathbf{A^c})$ iff $\mathsf{Fr} = \mathsf{lFr}$.
2. $\mathbf{\Lambda}(l, \mathsf{Fr}, \mathbf{A}) = \mathbf{\Lambda}(l, \mathsf{lFr}, \mathbf{A})$ iff $\mathbf{\Lambda}(l, \mathsf{Fr}, \mathbf{A^c}) = \mathbf{\Lambda}(l, \mathsf{lFr}, \mathbf{A^c})$ iff $\mathsf{Fr} = \mathsf{lFr}$.
3. $\mathbf{\Lambda}(g, \mathsf{lFr}, \mathbf{A^c}) = \mathbf{\Lambda}(g, \mathsf{CFr}, \mathbf{A^c})$ iff $\mathsf{lFr} = \mathsf{CFr}$.
4. $\mathbf{\Lambda}(l, \mathsf{lFr}, \mathbf{A^c}) = \mathbf{\Lambda}(l, \mathsf{CFr}, \mathbf{A^c})$ iff $\mathsf{lFr} = \mathsf{CFr}$.

**PROOF.** The first two statements are a consequence of Corollary 3.13. Let us prove the last two. If $\mathsf{CFr} \neq \mathsf{lFr}$, then there is an idempotent element $a \in A \setminus \{0, 1\}$. The proof finishes by noticing that $\square \bar{a} \vdash^l_{\mathsf{CFr}(\mathbf{A^c})} \square 0$ (and so $\square \bar{a} \vdash^g_{\mathsf{CFr}(\mathbf{A^c})} \square 0$) while $\square \bar{a} \nvdash^g_{\mathsf{lFr}(\mathbf{A^c})} \square 0$ (and so $\square \bar{a} \nvdash^l_{\mathsf{lFr}(\mathbf{A^c})} \square 0$). ∎

We notice that in previous proposition it remains open to characterize when $\mathbf{\Lambda}(g, \mathsf{lFr}, \mathbf{A}) = \mathbf{\Lambda}(g, \mathsf{CFr}, \mathbf{A})$ or when $\mathbf{\Lambda}(l, \mathsf{lFr}, \mathbf{A}) = \mathbf{\Lambda}(l, \mathsf{CFr}, \mathbf{A})$. Using Theorem 3.26 it is easy to prove that this open question reduces to characterize the $\mathbf{A}$'s such that $\Lambda(\mathsf{lFr}, \mathbf{A}) = \Lambda(\mathsf{CFr}, \mathbf{A})$, but this other question is also open.

There are five (meta)rules that play a remarkable role in the modal setting. These rules are[13]

$$(\mathsf{N})\ \varphi \vdash \square\varphi \qquad\qquad (\mathsf{Mon})\ \varphi \to \psi \vdash \square\varphi \to \square\psi$$

$$(\mathsf{Pref})\ \frac{\gamma \vdash \varphi}{\square\gamma \vdash \square\varphi} \qquad\qquad (\mathsf{Pref}^*)\ \frac{\Gamma \vdash \varphi}{\square\Gamma \vdash \square\varphi}$$

$$(\mathsf{OrdPres}^*)\ \frac{r \to \Gamma \vdash r \to \varphi}{\square\Gamma \vdash \square\varphi}$$

where $r$ is a variable not appearing in $\Gamma \cup \{\varphi\}$. These five rules will be respectively called *Necessity, Monotonicity, Prefixing, Infinitary Prefixing* and *Infinitary Order Preserving*. We notice that the rule ($\mathsf{OrdPres}^*$) is not closed under substitutions; substitutions may destroy the requirement concerning the variable $r$.

**PROPOSITION 3.22**
Let $\mathsf{K}$ be a class of Kripke frames.

1. $\mathbf{\Lambda}(g, \mathsf{K}, \mathbf{A^c})$ is closed under ($\mathsf{N}$), that is, $\varphi \vdash^g_{\mathsf{K}(\mathbf{A^c})} \square\varphi$.

2. $\mathbf{\Lambda}(g, \mathsf{K}, \mathbf{A^c})$ is closed under ($\mathsf{Mon}$), that is, $\varphi \to \psi \vdash^g_{\mathsf{K}(\mathbf{A^c})} \square\varphi \to \square\psi$.

3. $\mathbf{\Lambda}(l, \mathsf{K}, \mathbf{A^c})$ is closed under ($\mathsf{OrdPres}^*$), that is, if $r$ is a variable not appearing in the set $\Gamma \cup \{\varphi\}$ of formulas and $\{r \to \gamma : \gamma \in \Gamma\} \vdash^l_{\mathsf{K}(\mathbf{A^c})} r \to \varphi$, then $\square\Gamma \vdash^l_{\mathsf{K}(\mathbf{A^c})} \square\varphi$.

---

[13] As expected we use $\square\Gamma$ and $r \to \Gamma$ to denote, respectively, the sets $\{\square\gamma : \gamma \in \Gamma\}$ and $\{r \to \gamma : \gamma \in \Gamma\}$.

4. If $\mathsf{K} \subseteq \mathsf{CFr}$ then both $\mathbf{\Lambda}(g, \mathsf{K}, \mathbf{A^c})$ and $\mathbf{\Lambda}(l, \mathsf{K}, \mathbf{A^c})$ are closed under $(\mathsf{Pref}^*)$.

5. If $\mathsf{K} \subseteq \mathsf{IFr}$ and $\mathbf{\Lambda}(l, \mathsf{K}, \mathbf{A^c})$ satifies the Local Deduction Theorem, then $\mathbf{\Lambda}(l, \mathsf{K}, \mathbf{A^c})$ is closed under $(\mathsf{Pref})$. That is, if $\gamma \vdash^l_{\mathsf{K(A^c)}} \varphi$, then $\Box\gamma \vdash^l_{\mathsf{K(A^c)}} \Box\varphi$.

Proof. The proofs of the first four statements are straightforward. Let us now prove the fifth one. We assume that $\vdash^l_{\mathsf{K(A^c)}} \gamma^n \to \varphi$ for some $n \in \omega$. Then, since the local modal logic and the global one have the same theorems by the second item we get that $\vdash^l_{\mathsf{IFr(A^c)}} \Box(\gamma^n) \to \Box\varphi$. On the other hand, it holds that $\Box\gamma \vdash^l_{\mathsf{IFr(A^c)}} (\Box\gamma)^n \vdash^l_{\mathsf{IFr(A^c)}} \Box(\gamma^n)$ (see Proposition 3.10). Using both things we get by Modus Ponens that $\Box\gamma \vdash^l_{\mathsf{K(A^c)}} \Box\varphi$. ∎

Example 3.23
It is worth pointing out that the fifth statement is false in general when we drop the hypothesis concerning the Local Deduction Theorem. For example, if we consider $\mathbf{A}$ as the weak nilpotent minimum algebra given in Figure 5 then it is clear that $\neg\neg p \vdash_{\mathsf{IFr}} p$ (see Example A.4). However, we get that $\Box\neg\neg p \not\vdash^l_{\mathsf{IFr}} \Box p$ by taking a Kripke model with a unique point $w$ and such that $R(w, w) = 0.75$ and $V(p, w) = 0.5$.  ⊣

Remark 3.24 (Classical Setting)
In the literature of classical modal logics it is very common to only consider the set $\Lambda(\mathsf{K}, \mathbf{2})$ as a primitive notion. The reason to do this is essentially that in this setting both logics, the local and the global one, can be obtained from the set $\Lambda(\mathsf{K}, \mathbf{2})$ because

- $\Gamma \vdash^l_{\mathsf{K(2)}} \varphi$   iff   $\Lambda(\mathsf{K}, \mathbf{2}) \cup \Gamma \vdash_{\mathbf{2}} \varphi$.
- $\Gamma \vdash^g_{\mathsf{K(2)}} \varphi$   iff   $\{\Box^n\gamma : n \in \omega, \gamma \in \Gamma\} \vdash_{\mathsf{K(2)}} \varphi$.

While we will be able to prove, when the canonical frame construction works, the generalization to $\mathbf{A}$ of the first item (cf. Theorem 4.11) this is not the case with the second item. Indeed, in this article it remains open whether $\mathbf{\Lambda}(g, \mathsf{K}, \mathbf{A})$ is the smallest consequence relation extending $\mathbf{\Lambda}(l, \mathsf{K}, \mathbf{A})$ that is closed under the Monotonicity rule $\varphi \to \psi \vdash \Box\varphi \to \Box\psi$. The same question with canonical constants is also open. We only know, because the same proof than in the classical setting works, that

- if $\mathsf{K} \subseteq \mathsf{CFr}$, then $\Gamma \vdash^g_{\mathsf{K(A^c)}} \varphi$ iff $\{\Box^n\gamma : n \in \omega, \gamma \in \Gamma\} \vdash^l_{\mathsf{K(A^c)}} \varphi$.

Therefore, when the class of frames is crisp the answer to the previous question is positive and indeed we get something stronger: the global modal logic is the smallest consequence relation extending the local one that is closed under the Necessity rule $\varphi \vdash \Box\varphi$.  ⊣

Now, using a semantic argument we explain why the modal logic can be reduced to the non-modal one. This explanation is based on a canonical Kripke model construction (cf. Definition 4.7, Definition 4.19, etc.), and it is worth pointing out that it works for arbitrary $\mathbf{A}$'s (not only finite ones) with and without canonical constants.

Lemma 3.25
Let $\mathsf{K} \in \{\mathsf{Fr}, \mathsf{IFr}, \mathsf{CFr}\}$. Then, for every map $h : Fm \to A$ the following statements are equivalent:

1. $h$ is a non-modal homomorphism such that $h[\Lambda(\mathsf{K}, \mathbf{A})] = \{1\}$,

2. $h$ is a semantic modal valuation arising from a Kripke frame in $\mathsf{K}$.

The same equivalence also holds when we replace $\mathbf{A}$ with $\mathbf{A^c}$.

PROOF. We will only prove the case without canonical constants, since the other one is analogous replacing $\mathbf{A}$ with $\mathbf{A^c}$. The only non-trivial direction is the downwards one. In order to prove this we are going to define a canonical Kripke model belonging to the class $\mathsf{K}$. First of all we define the set $B$ as (i) $A$ if $\mathsf{K} = \mathsf{Fr}$, (ii) $\{a \in A : a = a \odot a\}$ if $\mathsf{K} = \mathsf{IFr}$, and (iii) $\{0, 1\}$ if $\mathsf{K} = \mathsf{CFr}$. It is obvious that in all three cases the set $B$ is closed under arbitrary joins. Next we define the Kripke model $\langle W_{can}, R_{can}, V_{can} \rangle$ given by

- the set $W_{can}$ is the set of semantic modal valuations $v$ arising from frames in $\mathsf{K}$,
- the accessibility relation $R_{can}(v_1, v_2)$ is defined as the largest element in $B$ below $\bigwedge \{v_1(\Box \varphi) \to v_2(\varphi) : \varphi \in Fm\}$,
- the evaluation map is defined by $V_{can}(p, v) := v(p)$ for every variable $p$.

By definition this Kripke model belongs to $\mathsf{K}$. To finish the proof it is enough to prove that $V_{can}(\varphi, v) = v(\varphi)$ for every modal formula $\varphi$ and every $v \in W_{can}$. By induction it suffices to prove that $\bigwedge \{R_{can}(v, v') \to v'(\varphi) : v' \in W_{can}\} = v(\Box \varphi)$ for every modal formula $\varphi$ and every $v \in W_{can}$.

The inequality $\bigwedge \{R_{can}(v, v') \to v'(\varphi) : v' \in W_{can}\} \geqslant v(\Box \varphi)$ trivially follows from the definition of the accessibility relation $R_{can}$.

Let us now try to prove the other inequality. Since $v \in W_{can}$ we know that there is a Kripke model $\langle W, R, V \rangle$ in $\mathsf{K}$ and a world $w \in W$ such that $v = V(\bullet, w)$. It is obvious that all worlds $w' \in W$ can be seen as members of $W_{can}$, and using that $R(w, w') \leqslant \bigwedge \{w(\Box \psi) \to w'(\psi) : \psi \in Fm\}$ it follows that $R(w, w') \leqslant R_{can}(w, w')$. Thus, $v(\Box \varphi) = w(\Box \varphi) = \bigwedge \{R(w, w') \to w'(\varphi) : w' \in W\} \geqslant \bigwedge \{R_{can}(w, w') \to w'(\varphi) : w' \in W\} \geqslant \bigwedge \{R_{can}(w, v') \to v'(\varphi) : v' \in W_{can}\} = \bigwedge \{R_{can}(v, v') \to v'(\varphi) : v' \in W_{can}\}$, which finishes the proof. ∎

THEOREM 3.26
Let $\mathsf{K} \in \{\mathsf{Fr}, \mathsf{IFr}, \mathsf{CFr}\}$. Then, for every set $\Gamma \cup \{\varphi\}$ it holds that

$$\Gamma \vdash^l_{\mathsf{K}(\mathbf{A})} \varphi \quad \text{iff} \quad \Gamma \cup \Lambda(\mathsf{K}, \mathbf{A}) \vdash_{\mathbf{A}} \varphi.$$

Analogously it holds that

$$\Gamma \vdash^l_{\mathsf{K}(\mathbf{A^c})} \varphi \quad \text{iff} \quad \Gamma \cup \Lambda(\mathsf{K}, \mathbf{A^c}) \vdash_{\mathbf{A^c}} \varphi.$$

PROOF. This trivially follows from the previous lemma together with our definition of the local modal logic. ∎

To finish this section we compare the classes of Boolean and crisp frames. Previously (see Theorem 3.16) we proved that $\mathsf{BFr}$ and $\mathsf{CFr}$ share the same valid formulas, but unfortunately the construction there used does not apply to the consequence relations involved in the local and global modal logics. Next we will see that this is true for the local and global modal logics when there are no canonical constants, and false if there are canonical constants.

PROPOSITION 3.27
Let $\mathbf{A}$ be a finite residuated lattice with some non trivial Boolean element. Then, $\Lambda(\mathsf{BFr}, \mathbf{A^c}) = \Lambda(\mathsf{CFr}, \mathbf{A^c})$, but $\mathbf{\Lambda}(l, \mathsf{BFr}, \mathbf{A^c}) \neq \mathbf{\Lambda}(l, \mathsf{CFr}, \mathbf{A^c})$ and $\mathbf{\Lambda}(g, \mathsf{BFr}, \mathbf{A^c}) \neq \mathbf{\Lambda}(g, \mathsf{CFr}, \mathbf{A^c})$.

---

- the set of axioms is the smallest set closed under substitutions containing
  - an axiomatic basis for $\mathbf{\Lambda(L_n)}$ (see [13]),
  - $\Box(\varphi \to \psi) \to (\Box\varphi \to \Box\psi)$, $\Box(\varphi \oplus \varphi) \leftrightarrow \Box\varphi \oplus \Box\varphi$ and $\Box(\varphi \odot \varphi) \leftrightarrow \Box\varphi \odot \Box\varphi$,
- the Modus Ponens rule and the Necessity rule $\varphi \vdash \Box\varphi$.

---

TABLE 1. Axiomatization of the logic $\mathbf{\Lambda}(g, \mathsf{CFr}, \mathsf{L_n})$

PROOF. The first part is Theorem 3.16. Let us now consider a Boolean element $a \notin \{0,1\}$. Now it is easy to see that

$$\Box 0 \to \overline{a} \vdash^l_{\mathsf{CFr}(\mathbf{A^c})} \neg\Box 0 \qquad \text{and} \qquad \Box 0 \to \overline{a} \vdash^g_{\mathsf{CFr}(\mathbf{A^c})} \neg\Box 0.$$

On the other hand, the fact that

$$\Box 0 \to \overline{a} \not\vdash^l_{\mathsf{BFr}(\mathbf{A^c})} \neg\Box 0 \qquad \text{and} \qquad \Box 0 \to \overline{a} \not\vdash^g_{\mathsf{BFr}(\mathbf{A^c})} \neg\Box 0$$

can be proved using the Boolean Kripke model with only point $w$ such that $R(w,w) = \neg a$. ∎

PROPOSITION 3.28
Let $\mathbf{A}$ be a finite residuated lattice. Then, $\Lambda(\mathsf{BFr}, \mathbf{A}) = \Lambda(\mathsf{CFr}, \mathbf{A})$, $\mathbf{\Lambda}(l, \mathsf{BFr}, \mathbf{A}) = \mathbf{\Lambda}(l, \mathsf{CFr}, \mathbf{A})$ and $\mathbf{\Lambda}(g, \mathsf{BFr}, \mathbf{A}) = \mathbf{\Lambda}(g, \mathsf{CFr}, \mathbf{A})$.

PROOF. The first part is Theorem 3.16. The inequalities $\mathbf{\Lambda}(l, \mathsf{BFr}, \mathbf{A}) \leqslant \mathbf{\Lambda}(l, \mathsf{CFr}, \mathbf{A})$ and $\mathbf{\Lambda}(g, \mathsf{BFr}, \mathbf{A}) \leqslant \mathbf{\Lambda}(g, \mathsf{CFr}, \mathbf{A})$ are trivial. The inequality $\mathbf{\Lambda}(l, \mathsf{BFr}, \mathbf{A}) \geqslant \mathbf{\Lambda}(l, \mathsf{CFr}, \mathbf{A})$ is a consequence of Theorem 3.26. Finally, the inequality $\mathbf{\Lambda}(g, \mathsf{BFr}, \mathbf{A}) \geqslant \mathbf{\Lambda}(g, \mathsf{CFr}, \mathbf{A})$ follows from the reduction (see Remark 3.24) of the glocal modal logic to the local one. ∎

### 3.3    *Previous related works*

FITTING axiomatizes in [22, 23] somehow the set $\Lambda(\mathsf{Fr}, \mathbf{A^c})$ (and also the set $\Lambda(\mathsf{CFr}, \mathbf{A^c})$) for the case that $\mathbf{A}$ is a finite Heyting algebra. The main difference with our approach is that he is using sequents and not formulas. FERMÜLLER and LANGSTEINER study in [21] the set[14] $\Lambda(\mathsf{CFr}, \mathbf{A^c})$ under the paradigm of tableaux calculus.

If we restrict our attention to closer approaches to ours there are rather few results in the literature focussing in the minimum modal logic $\mathbf{\Lambda}(g, \mathsf{Fr}, \mathbf{A})$ (and the same for $\mathbf{\Lambda}(g, \mathsf{Fr}, \mathbf{A})$). For example, in [39] HANSOUL and TEHEUX proved that Table 1 axiomatizes $\mathbf{\Lambda}(g, \mathsf{CFr}, \mathsf{L_n})$, but they do not consider the logic $\mathbf{\Lambda}(g, \mathsf{Fr}, \mathsf{L_n})$. Two of the few papers where the authors consider the full class of Kripke frames and not only crisp ones are [10, 47]. There it is shown that $\mathbf{\Lambda}(g, \mathsf{Fr}, [\mathbf{0, 1}]_{\mathbf{G}}) = \mathbf{\Lambda}(g, \mathsf{CFr}, [\mathbf{0, 1}]_{\mathbf{G}})$ and that Table 2 is an axiomatization for this logic. From the results in [39, 10] it is easy to see that $\mathbf{\Lambda}(l, \mathsf{CFr}, \mathsf{L_n})$ and $\mathbf{\Lambda}(l, \mathsf{Fr}, [\mathbf{0, 1}]_{\mathbf{G}})$ are axiomatized by its set of theorems (not only its axioms) and the unique rule of Modus Ponens. Therefore, it is clear that they satisfy the Local Deduction Theorem.

---

[14] Indeed, they study a more general notion of modality that they call distributed, being the necessity operator $\Box$ just a particular case.

---

- the set of axioms is the smallest set closed under substituions containing
  - an axiomatic basis for $\boldsymbol{\Lambda}\big([\mathbf{0},\mathbf{1}]_{\mathbf{G}}\big)$ (see [30]),
  - $\Box(\varphi \rightarrow \psi) \rightarrow (\Box\varphi \rightarrow \Box\psi)$ and $\neg\neg\Box\varphi \rightarrow \Box\neg\neg\varphi$,
- the Modus Ponens rule and the Necessity rule $\varphi \vdash \Box\varphi$.

---

TABLE 2. Axiomatization of the logic $\boldsymbol{\Lambda}\big(g,\mathsf{Fr},[\mathbf{0},\mathbf{1}]_{\mathbf{G}}\big)$

The logic $\boldsymbol{\Lambda}\big(g,\mathsf{CFr},[\mathbf{0},\mathbf{1}]_{\mathbf{L}}\big)$ is also axiomatized in [39], but adding to Table 1 the axioms $\Box(\varphi \oplus \varphi^n) \leftrightarrow ((\Box\varphi) \oplus (\Box\varphi)^n)$ (for every $n \in \omega$) and the infinitary rule[15]

$$(\mathsf{InfGreat})\ \ \frac{\varphi \oplus \varphi, \varphi \oplus \varphi^2, \ldots, \varphi \oplus \varphi^n, \ldots}{\varphi}$$

We have decided to denote this rule by the name *Infinitely Great* because it is somehow describing in the assumptions the infinitely great elements of an MV algebra, i.e., the elements of the radical filter (cf. [13, Definition 3.6.3]). Thus, the MV models of this rule must satisfy that the radical filter is trivial, and hence they have to be semisimple MV algebras. In [39] it is proved that the previous addition gives us an strong complete axiomatization for $\boldsymbol{\Lambda}\big(g,\mathsf{CFr},[\mathbf{0},\mathbf{1}]_{\mathbf{L}}\big)$. It is unknown whether this infinitary rule is admissible, i.e, whether we get the same set of theorems if the infinitary rule is deleted. What it is known is that $\Lambda(\mathsf{CFr},[\mathbf{0},\mathbf{1}]_{\mathbf{L}}) = \bigcap_{n\in\omega} \Lambda(\mathsf{CFr},\mathbf{L_n})$ because the same proof of [30, Theorem 5.4.30] works (it is enough to check that the model construction of this proof preserves crispness).

## 4   Completeness of the modal logic when there are canonical constants and A is finite

In this section we first analyze the advantages of having canonical constants in the language and then we provide a method to expand a complete axiomatization of $\boldsymbol{\Lambda}(\mathbf{A^c})$, whenever $\mathbf{A}$ is finite, into complete axiomatizations for the modal logics $\boldsymbol{\Lambda}(l,\mathsf{Fr},\mathbf{A^c})$, $\boldsymbol{\Lambda}(l,\mathsf{IFr},\mathbf{A^c})$, $\boldsymbol{\Lambda}(l,\mathsf{CFr},\mathbf{A^c})$ and $\boldsymbol{\Lambda}(g,\mathsf{CFr},\mathbf{A^c})$. For the last two cases we also need the assumption that $\mathbf{A}$ has a unique coatom. The proofs are based on the canonical model construction that is so well known in the classical setting. Unfortunately, axiomatizations for $\boldsymbol{\Lambda}(g,\mathsf{Fr},\mathbf{A^c})$ and $\boldsymbol{\Lambda}(g,\mathsf{IFr},\mathbf{A^c})$ are at present unknown, and we think that they cannot be settled using the canonical model construction, and hence some other argument must be used.

### 4.1   Benefits of having canonical constants

In the completeness proofs given in Sections 4.2, 4.3 and 4.4 we prove completeness with two assumptions. One assumption is finiteness. This assumption guarantees that $\vdash_{\mathbf{A^c}}$ is finitary[16]. The other assumption is the addition of canonical constants.

---

[15] In the non-modal case several infinitary rules have been considered with the property of giving strong completeness for the standard Łukasiewicz chain. We point out that in our opinion this infinitary rule (InfGreat) is one of the simpler ones because it only uses one formula (cf. [48]).

[16] In order to generalize the proofs here given to infinite residuated lattices, we suggest in Section 6 to use infinitary proof systems together with strongly complete axiomatizations of the non-modal logics.

The addition of canonical constants, except for the Łukasiewicz case (see Remark 4.5 and Section 5), seems really hard to overtake.

Of course, the main obvious benefit of having canonical constants in the language is the increasing of expressive power because canonical constants allow us to express certain rules inside our formal language. For example, for every finite $X \subseteq A$ the rule

$$\frac{((\overline{a_1} \to \varphi_1) \odot \ldots \odot (\overline{a_n} \to \varphi_n)) \to (\overline{a} \to \varphi) \text{ for every } a_1, \ldots, a_n, a \in X}{(\Box\varphi_1 \odot \ldots \odot \Box\varphi_n) \to \Box\varphi}$$

is valid (cf. Corollary B.5) in the global consequence associated with any frame such that $R$ only takes values inside $X$ (i.e., $R : W \times W \longrightarrow X$). This rule gives us a way to somehow rewrite Corollary B.5 using a rule that is closed under substitutions.

However, there are other benefits of having canonical constants that are more hidden. Next we illustrate with two propositions[17] the difference between having or not canonical constants in the language. These propositions suggest that in the presence of canonical constants the behaviour of a semantic modal valuation $v$ is somehow determined by the set $v^{-1}[\{1\}]$ of formulas. This yields a different behaviour in modal logics, depending whether there are canonical constants or not, since by definition the local modal logic only takes into account those semantic modal valuations that send the involved formulas to 1.

We remark that the main advantage of the second representation in the next proposition is that the set $\{v_2(\varphi) : \varphi \in Fm, v_1(\Box\varphi) = 1\}$ is closed under finite meets.

**Proposition 4.1**
Let us assume that the canonical constants are in the modal language. If $a \in A$, and $v_1$ and $v_2$ are two semantic modal valuations, then the following statements are equivalent:

- $a \leqslant \bigwedge \{v_1(\Box\varphi) \to v_2(\varphi) : \varphi \in Fm\}$,
- $a \leqslant \bigwedge \{v_2(\varphi) : \varphi \in Fm, v_1(\Box\varphi) = 1\}$.

**Proof.** One direction is trivial. To prove the other, let us assume that $a \leqslant \bigwedge\{v_2(\varphi) : \varphi \in Fm, v_1(\Box\varphi) = 1\}$ and let us consider a modal formula $\varphi$. We are left with the task of checking that $a \leqslant v_1(\Box\varphi) \to v_2(\varphi)$. If we prove that $v_1(\Box\varphi) \leqslant a \to v_2(\varphi)$, the assertion follows. The proof is finished by showing[18] that for every $b \in A$, if $b \leqslant v_1(\Box\varphi)$ then $b \leqslant a \to v_2(\varphi)$. For every $b \in A$, if $b \leqslant v_1(\Box\varphi)$ then $1 = b \to v_1(\Box\varphi) = v_1(\overline{b} \to \Box\varphi) = v_1(\Box(\overline{b} \to \varphi))$; hence $a \leqslant v_2(\overline{b} \to \varphi) = b \to v_2(\varphi)$; and so $b \leqslant a \to v_2(\varphi)$. ∎

**Example 4.2**
It is not difficult to find counterexamples to the previous proposition in case that we do not have canonical constants. For instance, let us consider the semantic modal

---

[17] Another example of this different behaviour between having or not canonical constants can be obtained comparing the results in Section 4.3 with Lemma 5.10.

[18] This way of conducting the proof could seem tricky since we can simply consider $b := v_1(\Box\varphi)$ (this can be done because we added one canonical constant for every element in $A$). However, the advantage of following this tricky approach is that it also works as far as we only have canonical constants for a subset $C$ satisfying for every $a \in A$ the condition $a = \sup\{c \in C : c \leqslant a\} = \inf\{c \in C : a \leqslant c\}$. This condition guarantees that the following three statements are equivalent: (i) $a \leqslant b$, (ii) for every $c \in C$, if $c \leqslant a$ then $c \leqslant b$, and (iii) for every $c \in C$, if $b \leqslant c$ then $a \leqslant c$. When we say that the proof also works for a subset $C$ satisfying the previous condition we mean that it is enough to replace "$b \in A$" with "$b \in C$" in the present proof to get a proof for the case that there are only canonical constants for elements in $C$. In all future proofs we will also consider the same tricky way in order to be as general as possible. It is worth pointing out that if $A$ is the real unit interval $[0, 1]$ then a subset $C$ satisfying this condition is the subset of its rational numbers.

valuations $w_0$ and $w_1$ given in the second diagram of Figure 1 over the MV chain $\mathbf{L_3}$ of three values (cf. [39, Lemma 5.4]). It is obvious that $w_1(\varphi) \in \{0, 1\}$ for every modal formula $\varphi$ without canonical constants. Thus,

- $1 \not\leqslant w_0(\Box 0) \to w_1(0)$, while
- $1 \leqslant \{w_1(\varphi) : \varphi \in Fm, w_0(\Box \varphi) = 1\}$. ⊣

PROPOSITION 4.3
Let us assume that the canonical constants are in the modal language. If $v_1$ and $v_2$ are two semantic modal valuations, then the following statements are equivalent:

- $v_1 = v_2$ (i.e., for every modal formula $\varphi$, it holds that $v_1(\varphi) = v_2(\varphi)$),
- $\{\varphi \in Fm : v_1(\varphi) = 1\} = \{\varphi \in Fm : v_2(\varphi) = 1\}$.

PROOF. For the non trivial implication, let us assume that $\{\varphi \in Fm : v_1(\varphi) = 1\} = \{\varphi \in Fm : v_2(\varphi) = 1\}$ and let us consider a modal formula $\varphi$. We are left with the task of checking that $v_1(\varphi) = v_2(\varphi)$. By symmetry it is enough to prove that $v_1(\varphi) \leqslant v_2(\varphi)$. If we prove that for every $b \in A$, if $b \leqslant v_1(\varphi)$ then $b \leqslant v_2(\varphi)$, then the assertion follows. Let $b$ be an element of $A$ such that $b \leqslant v_1(\varphi)$. Then, $1 = b \to v_1(\varphi) = v_1(\overline{b} \to \varphi)$. This clearly forces $1 = v_2(\overline{b} \to \varphi) = b \to v_2(\varphi)$. We thus get $b \leqslant v_2(\varphi)$. ∎

EXAMPLE 4.4
Again it is easy to find counterexamples when there are no canonical constants. For instance, using that in non-modal interpretations over the standard Gödel algebra $[\mathbf{0}, \mathbf{1}]_{\mathbf{G}}$ the set of non-modal formulas evaluated to 1 depends only on the relative ordering (see [3, Remark 2.6] for details) it is easy to find a counterexample for $[\mathbf{0}, \mathbf{1}]_{\mathbf{G}}$. ⊣

REMARK 4.5 ((Strongly) Characterizing formulas)
A particular case where Proposition 4.3 holds even without canonical constants is the case that $\mathbf{A}$ is a subalgebra of $[\mathbf{0}, \mathbf{1}]_{\mathbf{L}}$ (i.e., $\mathbf{A}$ is a simple MV algebra). The reason why this holds is that for every rational number $\alpha \in [0, 1]$ there is a non-modal formula $\eta_\alpha(p)$, from now on called the *characterizing formula of the interval* $[\alpha, 1]$, with only one variable such that for every $a \in A$, it holds that

$$\eta_\alpha^{\mathbf{A}}(a) = 1 \quad \text{iff} \quad a \in [\alpha, 1]. \tag{4.1}$$

The existence of these formulas is an easy consequence of McNaughton's Theorem (see [13, Theorem 9.1.5]). How can we use these formulas to prove Proposition 4.3? It suffices to realize[19] that for every semantic modal valuation $v$ over $\mathbf{A}$ and every modal formula $\varphi$, it holds that $v(\varphi) = \bigvee\{\alpha \in [0, 1] \cap \mathbb{Q} : \alpha \leqslant v(\varphi)\} = \bigvee\{\alpha \in [0, 1] \cap \mathbb{Q} : \eta_\alpha^{\mathbf{A}}(v(\varphi)) = 1\} = \bigvee\{\alpha \in [0, 1] \cap \mathbb{Q} : v(\eta_\alpha(\varphi)) = 1\} = \bigvee\{\alpha \in [0, 1] \cap \mathbb{Q} : \eta_\alpha(\varphi) \in v^{-1}[\{1\}]\}$.

In case that $\mathbf{A}$ is a finite subalgebra of $[\mathbf{0}, \mathbf{1}]_{\mathbf{L}}$ (i.e., $\mathbf{A}$ is $\mathbf{L_n}$ for some $n \in \omega$), then every interval $[\alpha, 1]$, with $\alpha$ a rational number, also has a *strongly characterizing formula* $\eta_\alpha(p)$ in the sense that besides (4.1) it also holds that

$$\mathbf{A} \models \eta_\alpha(p) \vee \neg\eta_\alpha(p) \approx 1. \tag{4.2}$$

---

[19] This idea was somehow noticed by HANSOUL and TEHEUX in [39, Proposition 5.5].

---

- the set of axioms is the smallest set closed under substitutions containing
  - the axiomatic basis for $\mathbf{\Lambda(A^c)}$,
  - $\Box 1$, $(\Box\varphi \wedge \Box\psi) \to \Box(\varphi \wedge \psi)$ and $\Box(\overline{a} \to \varphi) \leftrightarrow (\overline{a} \to \Box\varphi)$,
- the rules of a basis for $\mathbf{\Lambda(A^c)}$ and the Monotonicity rule $\varphi \to \psi \vdash \Box\varphi \to \Box\psi$.

---

TABLE 3. Axiomatization of the set $\Lambda(\mathsf{Fr}, \mathbf{A^c})$ when $\mathbf{A}$ is finite

This last condition is simply saying that for every $a \in A$, it holds that $\eta_\alpha^{\mathbf{A}}(a) \in \{0, 1\}$. In other words, $\eta_\alpha(p)$ is a strongly characterizing formula of the interval $[\alpha, 1]$ iff $\eta_\alpha^{\mathbf{A}}$ is the characteristic function of the interval $[\alpha, 1]$. Of course, by continuity reasons there are no strongly characterizing formulas for the case $[\mathbf{0}, \mathbf{1}]_{\mathbf{L}}$. A weaker property than (4.2) but one that can be accomplished for the case $[\mathbf{0}, \mathbf{1}]_{\mathbf{L}}$ is that the unary map $\eta_\alpha^{\mathbf{A}}$ is non decreasing. $\dashv$

## 4.2   Completeness of $\mathbf{\Lambda}\big(l, \mathsf{Fr}, \mathbf{A^c}\big)$ when $\mathbf{A}$ is finite

The aim of this section is to prove that the axiomatization given in Table 3 is characterizing the set $\Lambda(\mathsf{Fr}, \mathbf{A^c})$ in case that $\mathbf{A}$ is finite. From this result we will be able to give an axiomatization for $\mathbf{\Lambda}(l, \mathsf{Fr}, \mathbf{A^c})$. Throughout the rest of Section 4 we assume that we have fixed an axiomatic (axioms and rules) basis for $\mathbf{\Lambda(A^c)}$. Indeed, the rules in Table 3 include the axioms and rules of this axiomatic basis for $\mathbf{\Lambda(A^c)}$. It is obvious that all axioms and rules given in this table are sound (even in case that $\mathbf{A}$ is infinite).

DEFINITION 4.6 (cf. Table 3)
A *(many-valued) modal logic set* over a residuated lattice $\mathbf{A^c}$ is any set $\mathsf{L}$ of modal formulas closed under substitutions such that

- $\mathsf{L}$ contains an axiomatic basis for $\mathbf{\Lambda(A^c)}$,
- $\mathsf{L}$ contains the formulas of the form $\Box 1$, $(\Box\varphi \wedge \Box\psi) \to \Box(\varphi \wedge \psi)$ and $\Box(\overline{a} \to \varphi) \leftrightarrow (\overline{a} \to \Box\varphi)$,
- $\mathsf{L}$ is closed under the rules of a basis for $\mathbf{\Lambda(A^c)}$,
- $\mathsf{L}$ is closed under the rule (Mon).

$\mathsf{L}$ is said to be *consistent* in case that $L$ is not the set of all formulas. $\dashv$

By definition it is obvious that modal logic sets are closed under intersections. Hence, there is a minimum modal logic set. The minimum modal logic set is exactly the set described in Table 3. It is also clear that modal logic sets are closed under $\vdash_{\mathbf{A^c}}$. From now on whenever we refer to $\vdash_{\mathbf{A^c}}$ (i.e., $\mathbf{\Lambda(A^c)}$) we are considering all subformulas starting with a $\Box$ as propositional variables. Thus, by (2.8) it holds that $\Gamma \vdash_{\mathbf{A^c}} \varphi$ iff for every non-modal homomorphism $h$ from the algebra of modal formulas into $\mathbf{A^c}$, if $h[\Gamma] \subseteq \{1\}$ then $h(\varphi) = 1$. Using the closure under $\vdash_{\mathbf{A^c}}$ it is obvious that $\mathsf{L}$ is closed under Modus Ponens, and also that $\mathsf{L}$ is consistent iff $0 \notin \mathsf{L}$.

DEFINITION 4.7
The *canonical Kripke model* $\mathfrak{M}_{\mathbf{can}}^{\mathbf{c}}(\mathsf{L})$ associated with a consistent modal logic set $\mathsf{L}$ is the Kripke model $\langle W_{can}^c, R_{can}^c, V_{can}^c \rangle$ where

- the set $W^c_{can}$ is the set[20] of non-modal homomorphisms $v : \mathbf{Fm} \longrightarrow \mathbf{A^c}$ (we point out the language of the algebra of formulas includes the necessity modality) such that $v[\mathsf{L}] = \{1\}$,

- the accessibility relation $R^c_{can}$ is defined by

$$R^c_{can}(v_1, v_2) := \bigwedge \{v_1(\Box\varphi) \to v_2(\varphi) : \varphi \in Fm\},$$

- the evaluation map is defined by $V^c_{can}(p, v) := v(p)$ for every variable $p$.      ⊣

In the previous definition the superscript refers to the fact that there are canonical constants. We notice that it is the consistency of $\mathsf{L}$ what implies that the set $W^c_{can}$ just defined is non empty. This is so because $0 \notin \mathsf{L}$, and hence $\mathsf{L} \not\vdash_{\mathbf{A^c}} 0$, what gives a non-modal homomorphism $h$ from the set of modal formulas into $\mathbf{A^c}$ such that $h[\mathsf{L}] = \{1\}$.

By the same argument that was given for Proposition 4.1 we can prove that $R^c_{can}(v_1, v_2) = \bigwedge\{v_2(\varphi) : \varphi \in Fm, v_1(\Box\varphi) = 1\}$. Using that this last set is closed under finite meets we get that, since $\mathbf{A}$ is finite,

$$R^c_{can}(v_1, v_2) = \min\{v_2(\varphi) : \varphi \in Fm, v_1(\Box\varphi) = 1\}. \tag{4.3}$$

The proof of the following Truth Lemma could be simplified using this last equality, but we prefer to avoid its use in order to give a proof that can be adopted in the next section where we do not have canonical constants in the language. We also remark that the intuition behind the first claim inside the proof of Truth Lemma is given by the rule (OrdPres*) where the new variable $r$ is instantiated by every element of the residuated lattice.

LEMMA 4.8 (Truth Lemma)
Let $\mathsf{L}$ be a consistent modal logic set over a finite residuated lattice $\mathbf{A^c}$. The canonical Kripke model $\mathfrak{M}^c_{\mathfrak{can}}(\mathsf{L})$ satisfies $V^c_{can}(\varphi, v) = v(\varphi)$ for every formula $\varphi$ and every world $v$.

PROOF. The proof is done by induction on the formula. The only non trivial case is when this formula starts with a $\Box$. Let $v$ be a world of the canonical Kripke model. We have to prove that $V^c_{can}(\Box\varphi, v) = v(\Box\varphi)$ under the induction hypothesis that $V^c_{can}(\varphi, v') = v'(\varphi)$ for every world $v'$. Hence, it suffices to prove that

$$\bigwedge\{R^c_{can}(v, v') \to v'(\varphi) : v' \in W^c_{can}\} = v(\Box\varphi).$$

The fact that the inequality $\geqslant$ holds is a trivial consequence of the definition of the accessibility relation $R^c_{can}$. In order to check the other inequality we will prove that for every $a \in A$, if $a \leqslant \bigwedge\{R^c_{can}(v, v') \to v'(\varphi) : v' \in W^c_{can}\}$ then $a \leqslant v(\Box\varphi)$. Thus, let us consider an element $a$ such that $a \leqslant \bigwedge\{R^c_{can}(v, v') \to v'(\varphi) : v' \in W^c_{can}\}$, i.e., $R^c_{can}(v, v') \leqslant a \to v'(\varphi)$ for every world $v'$.

---

[20] Another way to present the points of the canonical model is to define the elements of $W^c_{can}$ as the sets $\Sigma$ of modal formulas (including canonical constants) that are closed under $\vdash_{\mathbf{A^c}}$ and such that for every modal formula $\varphi$ there is a unique element $a \in A$ such that $\overline{a} \leftrightarrow \varphi \in \Sigma$. It is clear that there is a 1-to-1 correspondence $h \longmapsto h^{-1}[\{1\}]$ between both presentations of points in the canonical model. The advantage of the presentation we have taken is that it trivially generalizes to the case that there are no canonical constants in the language.

**Claim 4.9**

For every element $b \in A$, it holds that

$$\mathsf{L} \cup \{\overline{b} \to (\overline{d} \to \psi) : \psi \in Fm, d \in A, d \leqslant v(\Box\psi)\} \vdash_{\mathbf{A^c}} \overline{b} \to (\overline{a} \to \varphi).$$

PROOF OF CLAIM. To prove this claim we take a non-modal homomorphism $h$ that sends all hypothesis to 1. Hence, for every formula $\psi$ and every element $d \in A$, if $d \leqslant v(\Box\psi)$ then $d \leqslant b \to h(\psi)$. That is, for every formula $\psi$, it holds that $v(\Box\psi) \leqslant b \to h(\psi)$. Hence, $b \leqslant \bigwedge\{v(\Box\psi) \to h(\psi) : \psi \in Fm\}$. Using that $h$ is a point of the canonical Kripke model, we get from the previous inequality that $b \leqslant R_{can}^c(v, h)$. By the assumed property of the element $a$ it follows that $b \leqslant a \to h(\varphi)$. Hence, $h(\overline{b} \to (\overline{a} \to \varphi)) = 1$.    Q.E.D. (Claim)

**Claim 4.10**

There is an $m \in \omega$, formulas $\psi_1, \ldots, \psi_m$ and elements $d_1, \ldots, d_m$ such that $d_i \leqslant v(\Box\psi_i)$ for every $i \in \{1, \ldots, m\}$ and

$$\mathsf{L} \vdash_{\mathbf{A^c}} (\bigwedge_{1 \leqslant i \leqslant m} (\overline{d_i} \to \psi_i)) \to (\overline{a} \to \varphi).$$

PROOF OF CLAIM. Using that $\mathbf{A}$ is finite, and so $\vdash_{\mathbf{A^c}}$ is finitary, it is obvious by the previous claim that for every $b \in A$ there is an $m_b \in \omega$, formulas $\psi_1^b, \ldots, \psi_{m_b}^b$ and elements $d_1^b, \ldots, d_{m_b}^b$ such that (i) $d_i^b \leqslant v(\Box\psi_i^b)$ for every $i \in \{1, \ldots, m_b\}$, and (ii) $\mathsf{L} \cup \{\overline{b} \to (\overline{d_i^b} \to \psi_i^b) : 1 \leqslant i \leqslant m_b\} \vdash_{\mathbf{A^c}} \overline{b} \to (\overline{a} \to \varphi)$. Let us assume that $A = \{b_1, \ldots, b_n\}$. Considering the finite sequence $\psi_1^{b_1}, \ldots, \psi_{m_{b_1}}^{b_1}, \psi_1^{b_2}, \ldots, \psi_{m_{b_2}}^{b_2}, \ldots, \psi_1^{b_n}, \ldots, \psi_{m_{b_n}}^{b_n}$ of formulas, and the finite sequence $d_1^{b_1}, \ldots, d_{m_{b_1}}^{b_1}, d_1^{b_2}, \ldots, d_{m_{b_2}}^{b_2}, \ldots, d_1^{b_n}, \ldots, d_{m_{b_n}}^{b_n}$ of elements we get that there is an $m \in \omega$, formulas $\psi_1, \ldots, \psi_m$ and elements $d_1, \ldots, d_m$ such that (i) $d_i \leqslant v(\Box\psi_i)$ for every $i \in \{1, \ldots, m\}$, and (ii) for every element $b \in A$ it holds that

$$\mathsf{L} \cup \{\overline{b} \to (\overline{d_i} \to \psi_i) : 1 \leqslant i \leqslant m\} \vdash_{\mathbf{A^c}} \overline{b} \to (\overline{a} \to \varphi).$$

Now we prove that these formulas and elements satisfy the desired property stated in the claim. To prove this we consider a non-modal homomorphism $h$ such that $h[\mathsf{L}] = \{1\}$. The fact that $h((\bigwedge_{1 \leqslant i \leqslant m}(\overline{d_i} \to \psi_i)) \to (\overline{a} \to \varphi)) = 1$ follows from the property (ii) taking $b$ as $h(\bigwedge_{1 \leqslant i \leqslant m}(\overline{d_i} \to \psi_i))$.    Q.E.D. (Claim)

Using[21] that $\mathsf{L}$ is closed under $\vdash_{\mathbf{A^c}}$, we get that $((\overline{d_1} \to \psi_1) \wedge \ldots \wedge (\overline{d_m} \to \psi_m)) \to (\overline{a} \to \varphi)$ belongs to $\mathsf{L}$. Therefore, $\Box((\overline{d_1} \to \psi_1) \wedge \ldots \wedge (\overline{d_m} \to \psi_m)) \to \Box(\overline{a} \to \varphi) \in \mathsf{L}$ by the Monotonicity rule. Thus, $(\Box(\overline{d_1} \to \psi_1) \wedge \ldots \wedge \Box(\overline{d_m} \to \psi_m)) \to \Box(\overline{a} \to \varphi) \in \mathsf{L}$ by the meet distributivity axiom. And so, by the fact that $\mathsf{L}$ is a modal logic set we get that $((\overline{d_1} \to \Box\psi_1) \wedge \ldots \wedge (\overline{d_m} \to \Box\psi_m)) \to (\overline{a} \to \Box\varphi)$ belongs to $\mathsf{L}$. Finally, using that $v(((\overline{d_1} \to \Box\psi_1) \wedge \ldots \wedge (\overline{d_m} \to \Box\psi_m)) \to (\overline{a} \to \Box\varphi)) = 1 \to (a \to v(\Box\varphi))$ and that $v[\mathsf{L}] = \{1\}$ we obtain that $a \to v(\Box\varphi) = 1$, i.e., $a \leqslant v(\Box\varphi)$. This finishes the proof. ∎

Now we can easily prove that the set $\Lambda(\mathsf{Fr}, \mathbf{A^c})$ is the one axiomatized in Table 3, and also give a presentation of $\boldsymbol{\Lambda}(l, \mathsf{Fr}, \mathbf{A^c})$. We point out that we do not know if $\boldsymbol{\Lambda}(g, \mathsf{Fr}, \mathbf{A^c})$ is precisely the consequence relation given by Table 3.

---

[21] It is worth pointing out that up to now we have only used in this proof the definition of $R_{can}^c$. In other words, we have not used any of the properties encapsulated in the fact that $\mathsf{L}$ is a modal logic set. This will allow us to use the same proof until the previous second claim for the case that we do not have canonical constants in the language (see the proof of Lemma 5.6).

Theorem 4.11 (Axiomatization of $\mathbf{\Lambda}(l, \mathsf{Fr}, \mathbf{A^c})$)

Let $\mathsf{L}$ be the smallest modal logic set over a finite residuated lattice $\mathbf{A^c}$. Then,

1. $v$ is a semantic modal valuation iff $v$ is a point of $\mathfrak{M}_{\mathbf{can}}^{\mathfrak{c}}(\mathsf{L})$.

2. $\Gamma \vdash_{\mathsf{Fr}(\mathbf{A^c})}^{l} \varphi$   iff   $\mathsf{L} \cup \Gamma \vdash_{\mathbf{A^c}} \varphi$,   for every set $\Gamma \cup \{\varphi\}$ of formulas.

3. $\mathsf{L} = \Lambda(\mathsf{Fr}, \mathbf{A^c})$.

4. $\mathbf{\Lambda}(l, \mathsf{Fr}, \mathbf{A^c})$ is axiomatized by (i) $\mathsf{L}$ as a set of axioms, and (ii) the rules of the basis for $\mathbf{\Lambda}(\mathbf{A^c})$.

Proof. 1): This first statement is a trivial consequence of the Truth Lemma.

2): It is clear that we have the following chain of equivalences

- $\mathsf{L} \cup \Gamma \vdash_{\mathbf{A^c}} \varphi$, iff
- for every non-modal homomorphism $h : \mathbf{Fm} \longrightarrow \mathbf{A^c}$, if $h[\mathsf{L}] = \{1\}$ and $h[\Gamma] \subseteq \{1\}$ then $h(\varphi) = 1$, iff
- for every $v \in W_{can}^c$, if $v[\Gamma] \subseteq \{1\}$ then $v(\varphi) = 1$, iff
- for every semantic modal valuation $v$, if $v[\Gamma] \subseteq \{1\}$ then $v(\varphi) = 1$, iff
- $\Gamma \vdash_{\mathsf{Fr}(\mathbf{A^c})}^{l} \varphi$.

3): This is a consequence of the previous item.

4): This is again a consequence of the second item.  ■

## 4.3   Completeness of $\mathbf{\Lambda}(l, \mathsf{IFr}, \mathbf{A^c})$ when $\mathbf{A}$ is finite

The purpose of this section is to prove that $\Lambda(\mathsf{IFr}, \mathbf{A^c})$, whenever $\mathbf{A}$ is a finite residuated lattice, is the smallest modal logic set over $\mathbf{A^c}$ that contains at least one of the following three kinds of formula (schema)

$$\Box(\varphi \rightarrow \psi) \rightarrow (\Box\varphi \rightarrow \Box\psi) \quad (\Box\varphi \odot \Box\psi) \rightarrow \Box(\varphi \odot \psi) \quad (\Box\varphi \odot \Box\varphi) \rightarrow \Box(\varphi \odot \varphi)$$

We remark that the first formula is (K) and that in Proposition 3.12 we proved that each of them characterizes $\mathsf{IFr}$. We will see that if any of these three formulas is in a modal logic set then all of them also belong to it. We will start studying the case of the last of these formulas (it seems the weakest one), and later we will show the connection between the three formulas.

Lemma 4.12

Let $\mathsf{L}$ be the smallest modal logic set over a finite residuated lattice $\mathbf{A^c}$ containing $(\Box\varphi \odot \Box\varphi) \rightarrow \Box(\varphi \odot \varphi)$. Then, the frame of $\mathfrak{M}_{\mathbf{can}}^{\mathfrak{c}}(\mathsf{L})$ is an idempotent frame.

Proof. Let $a$ be $R_{can}^c(v_1, v_2)$ for some points $v_1$ and $v_2$ of the canonical Kripke model. By (4.3) we know that there is a formula $\varphi$ such that $v_1(\Box\varphi) = 1$ and $v_2(\varphi) = a$. Then, $v_1(\Box\varphi \odot \Box\varphi) = 1 \odot 1 = 1$. Using that $(\Box\varphi \odot \Box\varphi) \rightarrow \Box(\varphi \odot \varphi) \in \mathsf{L}$ we get that $1 = v_1(\Box(\varphi \odot \varphi))$. On the other hand, by the Truth Lemma we know that $v_1(\Box(\varphi \odot \varphi)) \leqslant R_{can}^c(v_1, v_2) \rightarrow (v_2(\varphi) \odot v_2(\varphi)) = a \rightarrow (a \odot a)$. Therefore, $1 \leqslant a \rightarrow (a \odot a)$, i.e., $a = a \odot a$.  ■

Theorem 4.13 (Axiomatization of $\mathbf{\Lambda}(l, \mathsf{IFr}, \mathbf{A^c})$)

Let $\mathsf{L}$ be the smallest modal logic set over a finite residuated lattice $\mathbf{A^c}$ containing $(\Box\varphi \odot \Box\varphi) \rightarrow \Box(\varphi \odot \varphi)$. Then,

1. $v$ is a semantic modal valuation arising from an idempotent frame iff $v$ is a point of $\mathfrak{M}_{\mathbf{can}}^{\mathfrak{c}}(\mathsf{L})$.

2. $\Gamma \vdash^l_{\mathsf{IFr}(\mathbf{A^c})} \varphi$    iff    $\mathsf{L} \cup \Gamma \vdash_{\mathbf{A^c}} \varphi$,    for every set $\Gamma \cup \{\varphi\}$ of formulas.

3. $\mathsf{L} = \Lambda(\mathsf{IFr}, \mathbf{A^c})$.

4. $\mathbf{\Lambda}(l, \mathsf{IFr}, \mathbf{A^c})$ is axiomatized by (i) $\mathsf{L}$ as a set of axioms, and (ii) the rules of a basis for $\mathbf{\Lambda}(\mathbf{A^c})$.

PROOF. The proof is analogous to the one given for Theorem 4.11 (but now using Lemma 4.12). ∎

COROLLARY 4.14
Let $\mathsf{L}$ be a modal logic set over a finite residuated lattice $\mathbf{A^c}$. Then,

1. $\Box(\varphi \to \psi) \to (\Box\varphi \to \Box\psi) \in \mathsf{L}$, iff
2. $(\Box\varphi \odot \Box\psi) \to \Box(\varphi \odot \psi) \in \mathsf{L}$, iff
3. $(\Box\varphi \odot \Box\varphi) \to \Box(\varphi \odot \varphi) \in \mathsf{L}$.

PROOF. $1 \Rightarrow 2$): Applying two times the axiom (K) to the fact that $\varphi \to (\psi \to (\varphi \odot \psi)) \in \mathsf{L}$ we get that $\Box\varphi \to (\Box\psi \to \Box(\varphi \odot \psi)) \in \mathsf{L}$. Therefore, $(\Box\varphi \odot \Box\psi) \to \Box(\varphi \odot \psi) \in \mathsf{L}$.

$2 \Rightarrow 3$): This is trivial.

$3 \Rightarrow 1$): Let $\mathsf{L}'$ be the smallest modal logic set containing $(\Box\varphi \odot \Box\varphi) \to \Box(\varphi \odot \varphi)$. It suffices to prove that $\Box(\varphi \to \psi) \to (\Box\varphi \to \Box\psi) \in \mathsf{L}'$, but this is a trivial consequence of $\mathsf{L}' = \Lambda(\mathsf{IFr}, \mathbf{A^c})$ (see Theorem 4.13) together with Proposition 3.10. ∎

REMARK 4.15 (Idempotent Canonical Kripke model)
There is an alternative way to prove the results given in this Section 4.3. This other method consists on introducing an *idempotent canonical Kripke model* $\mathfrak{M}_{\mathbf{ican}}^{\mathfrak{c}}(\mathsf{L})$ as in Definition 4.7 except for the accessibility relation where we take the largest idempotent below $\bigwedge\{v_1(\Box\varphi) \to v_2(\varphi) : \varphi \in Fm\}$. This idempotent always exists because idempotent elements are closed under arbitrary joins. Then, the same strategy than in the proof of Lemma 4.8 allows us to conclude the Truth Lemma, but this time the two intermediate claims are (i) for every idempotent element $b \in A$, it holds that $\mathsf{L} \cup \{\overline{b} \to (\overline{d} \to \psi) : \psi \in Fm, d \in A, d \leqslant v(\Box\psi)\} \vdash_{\mathbf{A^c}} \overline{b} \to (\overline{a} \to \varphi)$ , and (ii) there is an $m \in \omega$, formulas $\psi_1, \ldots, \psi_m$ and elements $d_1, \ldots, d_m$ such that $d_i \leqslant v(\Box\psi_i)$ for every $i \in \{1, \ldots, m\}$ and $\mathsf{L} \vdash_{\mathbf{A^c}} ((\overline{d_1} \to \psi_1) \wedge \ldots \wedge (\overline{d_m} \to \psi_m))^n \to (\overline{a} \to \varphi)$ where $n$ is the cardinal of $A$. We have not adopted this method in this section because it was not necessary, but in the next section we will have to consider a similar method in order to get completeness for the class of crisp frames. ⊣

## 4.4  Completeness of $\mathbf{\Lambda}(l, \mathsf{CFr}, \mathbf{A^c})$ and $\mathbf{\Lambda}(g, \mathsf{CFr}, \mathbf{A^c})$ when $\mathbf{A}$ is finite and has a unique coatom

We have been looking for an axiomatization of $\mathbf{\Lambda}(l, \mathsf{CFr}, \mathbf{A^c})$ and $\mathbf{\Lambda}(g, \mathsf{CFr}, \mathbf{A^c})$ in case that $A$ is finite, but we have only succeeded when $\mathbf{A}$ has a unique coatom $k$. In this section we present the solution to this problem. The main difficulty to generalize our proof to the general case is that we need a method to reduce the consequence relation $\vdash_{\mathbf{A^c}}$ to its set of theorems, and we only know how to do it when there is a

> - the set of axioms is the smallest set closed under substitutions containing
>   - the axiomatic basis for $\mathbf{\Lambda(A^c)}$,
>   - $\Box 1$, $(\Box\varphi \wedge \Box\psi) \to \Box(\varphi \wedge \psi)$, $\Box(\overline{a} \to \varphi) \leftrightarrow (\overline{a} \to \Box\varphi)$ and $\Box(\overline{k} \vee \varphi) \to (\overline{k} \vee \Box\varphi)$,
> - the rules are
>   - those of the basis for $\mathbf{\Lambda(A^c)}$,
>   - the Monotonicity rule $\varphi \to \psi \vdash \Box\varphi \to \Box\psi$.

TABLE 4: Axiomatization of $\mathbf{\Lambda}(g, \mathsf{CFr}, \mathbf{A^c})$ when $\mathbf{A}$ is finite and $k$ is its unique coatom

unique coatom $k$. For this particular case we know that this can be done using the trivial equivalence

$$\gamma_1, \ldots, \gamma_n \vdash_{\mathbf{A^c}} \varphi \qquad \text{iff} \qquad \vdash_{\mathbf{A^c}} (\gamma_1 \wedge \ldots \wedge \gamma_n) \to (\varphi \vee \overline{k}), \qquad (4.4)$$

For the rest of this section we assume that $\mathbf{A}$ has a unique coatom (i.e., penultimate element) $k$. Under this assumption it is obvious that the only Boolean elements in $\mathbf{A}$ are the trivial ones 0 and 1; so indeed $\mathsf{BFr} = \mathsf{CFr}$. We already know that $\Box(\overline{k} \vee \varphi) \to (\overline{k} \vee \Box\varphi) \in \Lambda(\mathsf{CFr}, \mathbf{A^c})$. In this section we will prove that $\Lambda(\mathsf{CFr}, \mathbf{A^c})$ is precisely the smallest modal logic set $\mathsf{L}$ such that $\Box(\overline{k} \vee \varphi) \to (\overline{k} \vee \Box\varphi) \in \mathsf{L}$.

**DEFINITION 4.16** (cf. Table 4)
A *crisp (many-valued) modal logic set* over $\mathbf{A^c}$ (where $k$ is its unique coatom) is any modal logic set over $\mathbf{A^c}$ such that $\mathsf{L}$ contains the formulas of the form $\Box(\overline{k} \vee \varphi) \to (\overline{k} \vee \Box\varphi)$.

Thus, the statement claimed above can be rewritten as saying that $\Lambda(\mathsf{CFr}, \mathbf{A^c})$ is the smallest crisp modal logic set. A first idea to prove this statement is to show that the frame of $\mathfrak{M}^{\mathfrak{c}}_{\mathsf{can}}(\mathsf{L})$, where $\mathsf{L}$ is the smallest crisp modal logic set, is a crisp frame. Unfortunately this method does not work in general as the following proposition points out. It is worth stressing that, as a consequence of Theorem 4.22, we could replace in this proposition the set $\mathsf{L}$ with $\Lambda(\mathsf{CFr}, \mathbf{A^c})$.

**PROPOSITION 4.17**
Let $\mathbf{A}$ be the ordinal sum $\mathbf{A_1} \oplus \mathbf{A_2}$ of two finite $\mathsf{MTL}$ chains such that $\mathbf{A_1}$ and $\mathbf{A_2}$ are non trivial (i.e., $\min\{|A_1|, |A_2|\} \geqslant 2$). And let $\mathsf{L}$ be the smallest crisp modal logic set over $\mathbf{A^c}$. Then, the frame of $\mathfrak{M}^{\mathfrak{c}}_{\mathsf{can}}(\mathsf{L})$ is not crisp.

PROOF. Let us consider $a \in A$ as the idempotent element separating the components $A_1$ and $A_2$. Since $\mathbf{A_1}$ and $\mathbf{A_2}$ are non trivial it is obvious that $a \notin \{0^{\mathbf{A}}, 1^{\mathbf{A}}\}$, and also that for every $b \in A$, $a \odot b = a \wedge b$. Let $w_0$, $w_1$ and $w_2$ be the semantic modal valuations given in Figure 3. Next we are going to see that these three valuations belong to the corresponding canonical Kripke model (i.e., $w_i[\mathsf{L}] = \{1\}$) and that $R^c_{can}(w_0, w_2) = a$. This is enough in order to conclude that $\mathfrak{M}^{\mathfrak{c}}_{\mathsf{can}}(\mathsf{L})$ is not crisp.

In order to prove that $R^c_{can}(w_0, w_2) = a$ we have to check that $a = \bigwedge\{w_0(\Box\varphi) \to w_1(\varphi) : \varphi \in Fm\}$. This is trivial because $a = R(w_0, w_2) \leqslant \bigwedge\{w_0(\Box\varphi) \to w_2(\varphi) : \varphi \in Fm\} \leqslant w_0(\Box p) \to w_2(p) = 1 \to a = a$.

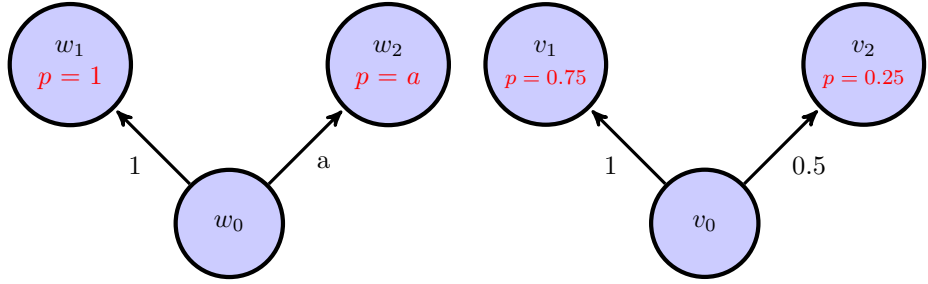

Now it is time to check that $w_i[\mathsf{L}] = \{1\}$. The only non trivial part is to check that $w_0[\{\Box(\overline{k} \vee \varphi) \rightarrow (\overline{k} \vee \Box \varphi) : \varphi \in Fm\}] = \{1\}$. In order to check this we first prove the following claim.

**Claim 4.18**
For every modal formula $\varphi$ (possibly including canonical constants), it holds that

- $a \leqslant w_1(\varphi)$  iff  $a \leqslant w_2(\varphi)$,
- if $a \not\leqslant w_1(\varphi)$ then $w_1(\varphi) = w_2(\varphi)$.

**Proof of Claim.** This claim can be straightforwardly proved by a simultaneous induction using the fact that $\mathbf{A}$ is an ordinal sum of the elements below $a$ with the elements above $a$.                                                                    Q.E.D. (Claim)

Now we are ready to prove that for every modal formula $\varphi$, it holds that $1 = w_0(\Box(\overline{k} \vee \varphi) \rightarrow (\overline{k} \vee \Box \varphi))$. Since $k$ is a coatom it is enough to prove that if $1 = w_0(\Box(\overline{k} \vee \varphi))$ then $1 = w_0(\Box \varphi)$. The fact that $1 = w_0(\Box(\overline{k} \vee \varphi)) = (1 \rightarrow w_1(\overline{k} \vee \varphi)) \wedge (a \rightarrow w_2(\overline{k} \vee \varphi))$ implies that $1 \leqslant w_1(\varphi)$. Hence using the previous claim we get that $a \leqslant w_2(\varphi)$. Therefore, $w_0(\Box \varphi) = (1 \rightarrow w_1(\varphi)) \wedge (a \rightarrow w_2(\varphi)) = 1$. ∎

Although we have just seen that the canonical Kripke model defined before does not help us to do the job of obtaining a completeness proof, next we will use an alternative canonical Kripke model which indeed does the job. This canonical Kripke model is forced to be crisp by definition. The idea behind this other canonical Kripke model is very similar to the one stated in Remark 4.15.

**Definition 4.19**
The *crisp canonical Kripke model* $\mathfrak{M}^{\mathfrak{c}}_{\mathfrak{ccan}}(\mathsf{L})$ associated with a crisp modal logic set $\mathsf{L}$ is the Kripke model $\langle W^c_{ccan}, R^c_{ccan}, V^c_{ccan} \rangle$ where

- the set $W^c_{ccan}$ is the set of non-modal homomorphisms $v : \mathbf{Fm} \longrightarrow \mathbf{A^c}$ (we point out the language of the algebra of formulas includes the necessity modality) such that $v[\mathsf{L}] = \{1\}$,

- the accessibility relation $R^c_{ccan}$ is defined by

$$R^c_{ccan}(v_1, v_2) := \begin{cases} 1 & \text{if } v_1(\Box \varphi) \leqslant v_2(\varphi) \text{ for every } \varphi \\ 0 & \text{otherwise} \end{cases}$$

- the evaluation map is defined by $V^c_{ccan}(p, v) := v(p)$ for every variable $p$.                ⊣

Analogously to what we pointed out after Definition 4.7 it is not difficult to see that if $\mathbf{A}$ is finite then

$$R_{ccan}^c(v_1, v_2) = \begin{cases} 1 & \text{if } \forall \varphi (v_1(\Box \varphi) = 1 \Rightarrow v_2(\varphi) = 1) \\ 0 & \text{otherwise} \end{cases}$$

**Lemma 4.20** (Truth Lemma)
Let $\mathsf{L}$ be a crisp modal logic set over $\mathbf{A^c}$ (where $k$ is its unique coatom). The crisp canonical Kripke model $\mathfrak{M}_{ccan}^c(\mathsf{L})$ satisfies that $V_{ccan}^c(\varphi, v) = v(\varphi)$ for every formula $\varphi$ and every world $v$.

**Proof.** The beginning of this proof is the same one that was given for Lemma 4.8 (replacing the canonical Kripke model for $\mathsf{L}$ with the crisp one for the same $\mathsf{L}$) until Claim 4.9. Now we continue the proof from this point.

**Claim 4.21**
There are $m \in \omega$, formulas $\psi_1, \ldots, \psi_m$ and elements $d_1, \ldots, d_m$ such that $d_i \leqslant v(\Box \psi_i)$ for every $i \in \{1, \ldots, m\}$ and

$$\mathsf{L} \cup \{\overline{d_1} \to \psi_1, \ldots, \overline{d_m} \to \psi_m\} \vdash_{\mathbf{A^c}} \overline{a} \to \varphi.$$

**Proof of Claim.** By the same proof given for Claim 4.9 we know[22] that $\mathsf{L} \cup \{\overline{d} \to \psi : \psi \in Fm, d \in A, d \leqslant v(\Box \psi)\} \vdash_{\mathbf{A^c}} \overline{a} \to \varphi$. Using the fact that $\vdash_{\mathbf{A^c}}$ is finitary the proof is finished. Q.E.D. (Claim)

Using the fact that $k$ is the unique coatom (cf. (4.4)) it follows that $\mathsf{L} \vdash_{\mathbf{A^c}} ((\overline{d_1} \to \psi_1) \land \ldots \land (\overline{d_m} \to \psi_m)) \to (\overline{k} \lor (\overline{a} \to \varphi))$. Using that $\mathsf{L}$ is closed under $\vdash_{\mathbf{A^c}}$, we get that $((\overline{d_1} \to \psi_1) \land \ldots \land (\overline{d_m} \to \psi_m)) \to (\overline{k} \lor (\overline{a} \to \varphi)) \in \mathsf{L}$. Therefore, using the properties of crisp modal logic sets we get that $((\overline{d_1} \to \Box \psi_1) \land \ldots \land (\overline{d_m} \to \Box \psi_m)) \to (\overline{k} \lor (\overline{a} \to \Box \varphi)) \in \mathsf{L}$. Finally, using that $v((\overline{d_1} \to \Box \psi_1) \land \ldots \land (\overline{d_m} \to \Box \psi_m)) = 1$ and that $v[\mathsf{L}] = \{1\}$ we obtain that $k \lor (a \to v(\Box \varphi)) = 1$. Hence, $a \leqslant v(\Box \varphi)$. This finishes the proof. ∎

**Theorem 4.22** (Axiomatization of $\mathbf{\Lambda}(l, \mathsf{CFr}, \mathbf{A^c})$ and $\mathbf{\Lambda}(g, \mathsf{CFr}, \mathbf{A^c})$)
Let $\mathsf{L}$ be the smallest crisp modal logic set over $\mathbf{A^c}$ (where $k$ is its unique coatom). Then,

1. $v$ is a semantic modal valuation arising from a crisp frame iff $v$ is a point of $\mathfrak{M}_{ccan}^c(\mathsf{L})$.
2. $\Gamma \vdash_{\mathsf{CFr}(\mathbf{A^c})}^l \varphi$   iff   $\mathsf{L} \cup \Gamma \vdash_{\mathbf{A^c}} \varphi$,   for every set $\Gamma \cup \{\varphi\}$ of formulas.
3. $\mathsf{L} = \Lambda(\mathsf{CFr}, \mathbf{A^c})$.
4. $\mathbf{\Lambda}(g, \mathsf{CFr}, \mathbf{A^c})$ is axiomatized by the axioms and rules given in Table 4.

**Proof.** The first three items are proved like in Theorem 4.11 (but this time using the crisp canonical Kripke model). For the last one we notice that the statement

$$\Gamma \vdash_{\mathsf{CFr}(\mathbf{A^c})}^g \varphi \quad \text{iff} \quad \{\Box^n \gamma : n \in \omega, \gamma \in \Gamma\} \vdash_{\mathsf{CFr}(\mathbf{A^c})}^l \varphi$$

---

[22] It is worth noting that this statement can be rewritten as saying that for every element $b \in \{0, 1\}$, it holds that

$$\mathsf{L} \cup \{\overline{b} \to (\overline{a} \to \psi) : \psi \in Fm, d \in A, d \leqslant v(\Box \psi)\} \vdash_{\mathbf{A^c}} \overline{b} \to (\overline{a} \to \varphi).$$

This last version shows the connection with crisp frames where the accessibility relation only takes values in $\{0, 1\}$.

is easily proved considering the submodel generated by a world (see for instance [4]). Hence, $\mathbf{\Lambda}(g, \mathsf{CFr}, \mathbf{A^c})$ is the smallest logic extending $\mathbf{\Lambda}(l, \mathsf{CFr}, \mathbf{A^c})$ that is closed under the Necessity rule. Thus, using the second item we get that $\mathbf{\Lambda}(g, \mathsf{CFr}, \mathbf{A^c})$ is exactly the consequence relation given in Table 4. ∎

Using the previous theorem it is obvious that $\Box(\varphi \to \psi) \to (\Box\varphi \to \Box\psi)$ is a theorem of the calculus given in Table 4. Therefore, if in the calculus given in Table 4 we would have replaced the Monotonicity rule with the normality axiom (K) together with the Necessity rule, then we would get another complete axiomatization of $\mathbf{\Lambda}(g, \mathsf{CFr}, \mathbf{A^c})$.

To finish this section we give a semantic argument that justifies why in order to obtain a semantic completeness for crisp frames it is enough to add the formulas $\Box(\overline{k} \vee \varphi) \to (\overline{k} \vee \Box\varphi)$. This result says that when these formulas are valid then we can replace a Kripke model with a crisp "submodel" (the one given by taking the same worlds but replacing the accessibility relation with its associated crisp one). A particular case of this construction is the replacement of the canonical Kripke model used in previous sections with the crisp canonical Krikpe model introduced in this section.

LEMMA 4.23

Let us assume that $\mathbf{A}$ is finite (with $k$ its unique coatom) and let us assume that there are canonical constants in the language. If $\mathfrak{M} = \langle W, R, V \rangle$ is a Kripke model and $w \in W$, then the following statements are equivalent:

1. $\mathfrak{M}, w \models^1 \{\Box(\overline{k} \vee \varphi) \to (\overline{k} \vee \Box\varphi) : \varphi \in Fm\}$,
2. $V(\Box\varphi, w) = \bigwedge\{V(\varphi, w') : R(w, w') = 1, w' \in W\}$ for every modal formula $\varphi$,
3. there is a crisp Kripke model $\langle W', R', V' \rangle$ and a world $w' \in W'$ such that $V(\bullet, w) = V'(\bullet, w')$.

PROOF. $1 \Rightarrow 2$ : Let us consider a modal formula $\varphi$. We define $a = \bigwedge\{V(\varphi, w') : R(w, w') = 1, w' \in W\}$. It is trivial that $V(\Box\varphi, w) \leqslant a$, and our aim is to prove that $a \leqslant V(\Box\varphi, w)$ also holds. We know by the definition of $a$ that $V(\Box(\overline{k} \vee (\overline{a} \to \varphi)), w) = 1$. Then, using the assumption for the particular case of the formula $\overline{a} \to \varphi$ we get that $1 \leqslant k \vee (a \to V(\Box\varphi, w))$. Hence, $a \leqslant V(\Box\varphi, w)$.

$2 \Rightarrow 3$ : It is obvious that we can take the Kripke model $\langle W, R', V \rangle$ where $R'(w, w')$ is defined as (i) 1 if $R(w, w') = 1$, and (ii) 0 if $R(w, w') \neq 1$.

$3 \Rightarrow 1$ : This is trivial. ∎

It is worth pointing out that in the previous proof it is crucial that the set $Fm$ of modal formulas allows to use canonical constants. Indeed, if there are no canonical constants in the language then it is not difficult to find counterexamples to the equivalence of the first two statements in the previous lemma. For example, let us consider the Gödel chain with universe $\{0, 0.25, 0.5, 0.75, 1\}$ and the semantic modal valuation $v_0$ given in Figure 3. Using that $v_0(\Box(\overline{0.5} \vee p) \to (\overline{0.5} \vee \Box p)) = 0.5 \neq 1$ it is obvious that $v_0$ cannot be obtained as a semantic modal valuation arising from a crisp Kripke frame. On the other hand, by induction it is obvious that for every modal formula $\varphi$ without canonical constants, it holds that

- $v_1(\varphi) \in \{0, 0.75, 1\}$,
- $v_1(\varphi) = 1$ iff $v_2(\varphi) = 1$,

---

- the set of axioms is the smallest set closed under substitutions containing
  - an axiomatic basis for $\mathbf{\Lambda(L_n)}$ (see [13, Section 8.5]),
  - $\Box 1$ and $(\Box\varphi \wedge \Box\psi) \rightarrow \Box(\varphi \wedge \psi)$,
- the Modus Ponens rule, the Monotonicity rule and for every $a \in L_n \setminus \{0\}$ the rule

$$(\mathsf{R}_a) \ \frac{(\eta_{a_2 \odot b}(\varphi_2) \wedge \eta_{a_3 \odot b}(\varphi_3) \wedge \ldots \wedge \eta_{a_n \odot b}(\varphi_n)) \rightarrow \eta_{a \odot b}(\varphi) \ \text{for every } b \in L_n, \ b > \neg a}{(\eta_{a_2}(\Box\varphi_2) \wedge \eta_{a_3}(\Box\varphi_3) \ldots \wedge \eta_{a_n}(\Box\varphi_n)) \rightarrow \eta_a(\Box\varphi)}$$

where $a_2 = \frac{1}{n-1}$, $a_3 = \frac{2}{n-1}$, $\ldots$, $a_{n-1} = \frac{n-2}{n-1}$ and $a_n = 1$.

TABLE 5. Axiomatization of the set $\Lambda(\mathsf{Fr}, \mathbf{L_n})$

- $v_1(\varphi) = 0.75$ iff $v_2(\varphi) = 0.25$,
- $v_1(\varphi) = 0$ iff $v_2(\varphi) = 0$.

Using this it is very easy to check that $v_0(\Box(\overline{0.75} \vee \varphi) \rightarrow (\overline{0.75} \vee \Box\varphi)) = 1$ for every modal formula $\varphi$ without canonical constants.

## 5 Completeness of the modal logic given by a finite MV chain

In this section our aim is to study the set $\Lambda(\mathsf{Fr}, \mathbf{L_n})$ where $L_n$ is the finite MV chain with $n$ elements. From now on we assume that $n \geqslant 3$ is fixed. We stress that there are no canonical constants in the language, and that the possibility operator $\diamond$ is definable since MV algebras are involutive (i.e., $\diamond\varphi$ is an abbreviation for $\neg\Box\neg\varphi$). The completeness proofs that we will give are based on the ones of previous sections (those proofs were with canonical constants in the language) together with the fact that on $\mathbf{L_n}$ there are strongly characterizing formulas $\eta_a(p)$ (see Remark 4.5). From now on whenever we talk about $\eta_a$ we assume, except if it is explicitly said something different, that $\eta_a$ is a strongly characterizing formula for the interval $[a, 1]$. Besides focussing on $\Lambda(\mathsf{Fr}, \mathbf{L_n})$ we also consider in this section the set $\Lambda(\mathsf{CFr}, \mathbf{L_n})$, which clearly coincides with $\Lambda(\mathsf{IFr}, \mathbf{L_n})$ because idempotent frames over $\mathbf{L_n}$ coincide with crisp frames.

### 5.1 Completeness of $\mathbf{\Lambda}\big(l, \mathsf{Fr}, L_n\big)$

We will show that Table 5 gives us an axiomatization of the set $\Lambda(\mathsf{Fr}, \mathbf{L_n})$. The idea behind the rules $(\mathsf{R}_a)$ given in this table is that they are a way to rewrite the rules (where $a_2, a_3, \ldots, a_n$ are like in Table 5)

$$\frac{((\overline{a_2} \rightarrow \varphi_2) \wedge \ldots \wedge (\overline{a_n} \rightarrow \varphi_n)) \rightarrow (\overline{a} \rightarrow \varphi)}{((\overline{a_2} \rightarrow \Box\varphi_2) \wedge \ldots \wedge (\overline{a_n} \rightarrow \Box\varphi_n)) \rightarrow (\overline{a} \rightarrow \Box\varphi)}$$

without using canonical constants. First of all, let us check that all axioms and rules in Table 5 are sound.

LEMMA 5.1
The logic $\mathbf{\Lambda}\big(g, \mathsf{Fr}, \mathbf{L_n}\big)$ is closed under all axioms and rules given in Table 5.

Proof. It suffices to prove that it is closed under the rules $(\mathsf{R}_a)$. To see this, let $a \in A$ and let $V$ be a valuation from a Kripke model. We assume that $V((\eta_{a_2 \odot b}(\varphi_2) \wedge \eta_{a_3 \odot b}(\varphi_3) \wedge \ldots \wedge \eta_{a_n \odot b}(\varphi_n)) \rightarrow \eta_{a \odot b}(\varphi), w) = 1$ for every world $w$ and every $b > \neg a$. We have to prove that $V((\eta_{a_2}(\Box \varphi_2) \wedge \eta_{a_3}(\Box \varphi_3) \ldots \wedge \eta_{a_n}(\Box \varphi_n)) \rightarrow \eta_a(\Box \varphi), w) = 1$ for an arbitrary world $w$. Using that the formulas $\eta$ are strongly characterizing it suffices to prove that if $V(\eta_{a_i}(\Box \varphi_i), w) = 1$ for every $i \in \{2, \ldots, n\}$, then $V(\eta_a(\Box \varphi), w) = 1$. That is, we have to prove that if $a_i \leqslant V(\Box \varphi_i, w)$ for every $i \in \{2, \ldots, n\}$, then $a \leqslant V(\Box \varphi, w)$. Hence, let us assume that $a_i \leqslant V(\Box \varphi_i, w)$ for every $i \in \{2, \ldots, n\}$. We must prove that $a \leqslant V(\Box \varphi, w)$, i.e., that $a \leqslant R(w, w') \rightarrow V(\varphi, w')$ for every world $w'$. Thus, we consider a world $w'$ and we define $b := R(w, w')$. If $b \leqslant \neg a$ then it is obvious that $a \leqslant b \rightarrow 0 \leqslant R(w, w') \rightarrow V(\varphi, w')$. Let us now consider the case that $b > \neg a$. The fact that $a_i \leqslant V(\Box \varphi_i, w)$ tells us that $a_i \odot b \leqslant V(\varphi_i, w')$ for every $i \in \{2, \ldots, n\}$. Therefore, $V(\eta_{a_i \odot b}(\varphi_i), w') = 1$ for every $i \in \{2, \ldots, n\}$. Using the assumption about the upper part of the rule $(\mathsf{R}_a)$ we get that $V(\eta_{a \odot b}(\varphi), w') = 1$. Thus, $a \leqslant b \rightarrow V(\varphi, w') = R(w, w') \rightarrow V(\varphi, w')$. This finishes the proof. ∎

Before giving the proof based on the canonical model construction, as a matter of example next we analyze in detail the rules in Table 5 for the case of the MV chain with three points.

Example 5.2 (Case of Ł$_{\mathbf{3}}$)
For $n = 3$ the rules $(\mathsf{R}_a)$ considered are

$$(\mathsf{R}_{0.5}) \quad \frac{(\eta_{0.5}(\varphi_2) \wedge \eta_1(\varphi_3)) \rightarrow \eta_{0.5}(\varphi)}{(\eta_{0.5}(\Box \varphi_2) \wedge \eta_1(\Box \varphi_3)) \rightarrow \eta_{0.5}(\Box \varphi)}$$

$$(\mathsf{R}_1) \quad \frac{(\eta_0(\varphi_2) \wedge \eta_{0.5}(\varphi_3)) \rightarrow \eta_{0.5}(\varphi) \qquad (\eta_{0.5}(\varphi_2) \wedge \eta_1(\varphi_3)) \rightarrow \eta_1(\varphi)}{(\eta_{0.5}(\Box \varphi_2) \wedge \eta_1(\Box \varphi_3)) \rightarrow \eta_1(\Box \varphi)}$$

Using that $\eta_0(p) := 1$, $\eta_{0.5}(p) := p \oplus p$ and $\eta_1(p) := p \odot p$, we get that the previous rules can be rewritten as

$$(\mathsf{R}_{0.5}) \quad \frac{((\varphi_2 \oplus \varphi_2) \wedge (\varphi_3 \odot \varphi_3)) \rightarrow (\varphi \oplus \varphi)}{((\Box \varphi_2 \oplus \Box \varphi_2) \wedge (\Box \varphi_3 \odot \Box \varphi_3)) \rightarrow (\Box \varphi \oplus \Box \varphi)}$$

$$(\mathsf{R}_1) \quad \frac{(1 \wedge (\varphi_3 \oplus \varphi_3)) \rightarrow (\varphi \oplus \varphi) \qquad ((\varphi_2 \oplus \varphi_2) \wedge (\varphi_3 \odot \varphi_3)) \rightarrow (\varphi \odot \varphi)}{((\Box \varphi_2 \oplus \Box \varphi_2) \wedge (\Box \varphi_3 \odot \Box \varphi_3)) \rightarrow (\Box \varphi \odot \Box \varphi)}$$

Some theorems that can be obtained using rule $(\mathsf{R}_{0.5})$ are

$$((\Box p \oplus \Box p) \wedge (\Box q \odot \Box q)) \rightarrow (\Box(p \odot q) \oplus \Box(p \odot q))$$

$$((\Box p \oplus \Box p) \wedge (\Box(q \oplus q \oplus (\neg p \odot \neg p)) \odot \Box(q \oplus q \oplus (\neg p \odot \neg p)))) \rightarrow (\Box q \oplus \Box q)$$

$$((\Box(q \oplus \neg p) \oplus \Box(q \oplus \neg p)) \wedge (\Box p \odot \Box p)) \rightarrow (\Box q \oplus \Box q)$$

$$((\Box 1 \oplus \Box 1) \wedge (\Box 1 \odot \Box 1)) \rightarrow (\Box(p \vee \neg p) \oplus \Box(p \vee \neg p)),$$

and some theorems using rule $(\mathsf{R}_1)$ are

$$((\Box 0 \oplus \Box 0) \wedge (\Box 1 \odot \Box 1)) \rightarrow (\Box(p \vee \neg p) \odot \Box(p \vee \neg p))$$

$$((\Box(p \odot p) \oplus \Box(p \odot p)) \wedge (\Box(p \oplus p) \odot \Box(p \oplus p))) \rightarrow (\Box p \odot \Box p)$$

$$((\Box p \oplus \Box p) \wedge (\Box(p \vee \neg p) \odot \Box(p \vee \neg p))) \to (\Box(p \vee q \vee \neg q) \odot \Box(p \vee q \vee \neg q)).$$

Using standard matrix arguments it is possible to prove that the rules $(\mathsf{R}_{0.5})$ and $(\mathsf{R}_1)$ are independent of all the other axioms and rules in Table 5. For example, the matrix $\langle \mathbf{L_3} \times \mathbf{L_2}, \Box, \{(1,1)\} \rangle$ where

$$\Box x := \begin{cases} (0,0), & \text{if } \pi_1(x) = 0 \text{ (i.e., } x \in \{(0,0),(0,1)\}) \\ (1,1), & \text{if } \pi_1(x) \neq 0 \end{cases}$$

is a model of all axioms and rules except[23] for $(\mathsf{R}_{0.5})$. And the matrix $\langle \mathbf{L_3} \times \mathbf{L_2}, \Box', \{(1,1)\} \rangle$ where

$$\Box' x := \begin{cases} (0.5,1), & \text{if } \pi_1(x) \neq 1 \\ (1,1), & \text{if } \pi_1(x) = 1 \text{ (i.e., } x \in \{(1,0),(1,1)\}) \end{cases}$$

is a model of all axioms and rules except[24] for $(\mathsf{R}_1)$. $\dashv$

**DEFINITION 5.3** (cf. Table 5)
A *(many-valued) modal logic set* over $\mathbf{L_n}$ is any set $\mathsf{L}$ of modal formulas closed under substitutions such that

- $\mathsf{L}$ contains an axiomatic basis for $\boldsymbol{\Lambda}(\mathbf{L_n})$ (see [13, Section 8.5]),
- $\mathsf{L}$ contains the formulas of the form $\Box 1$ and $(\Box\varphi \wedge \Box\psi) \to \Box(\varphi \wedge \psi)$,
- $\mathsf{L}$ is closed under Modus Ponens and the Monotonicity rule,
- for every $a \in \mathbf{L_n} \setminus \{0\}$, $\mathsf{L}$ is closed under the rule $(\mathsf{R}_a)$. $\dashv$

It is trivial that the minimum modal logic set exists, and it is exactly the set described in Table 5. And it is obvious that all modal logic sets are closed under $\vdash_{\mathbf{L_n}}$. Next we show an slight generalization of the rule $(\mathsf{R}_a)$ that it is also preserved by modal logic sets.

**LEMMA 5.4**
Let $\mathsf{L}$ be a modal logic set over $\mathbf{L_n}$. For every $a_1, \ldots, a_m, a \in \mathbf{L_n}$ it holds that $\mathsf{L}$ is closed under the rule

$$(\mathsf{R}_a^{a_1 \ldots a_m}) \frac{(\eta_{a_1 \odot b}(\varphi_1) \wedge \ldots \wedge \eta_{a_m \odot b}(\varphi_m)) \to \eta_{a \odot b}(\varphi) \text{ for every } b \in \mathbf{L_n},\, b > \neg a}{(\eta_{a_1}(\Box\varphi_1) \wedge \ldots \wedge \eta_{a_m}(\Box\varphi_m)) \to \eta_a(\Box\varphi)}$$

**PROOF.** This is an straightforward consequence of the fact that for every $a \in \mathbf{L_n}$ it holds that $\vdash_{\mathbf{L_n}} \eta_a(\varphi_1 \wedge \varphi_2) \leftrightarrow (\eta_a(\varphi_1) \wedge \eta_a(\varphi_2))$, and so this formula belongs to $\mathsf{L}$. ∎

**DEFINITION 5.5**
The *canonical Kripke model* $\mathfrak{M}_{\mathtt{can}}(\mathsf{L})$ associated with a modal logic set $\mathsf{L}$ over $\mathbf{L_n}$ is the Kripke model $\langle W_{can}, R_{can}, V_{can} \rangle$ defined as in Definition 4.7 except for the fact that now there are no canonical constants in the formulas. $\dashv$

The proof of the following Truth Lemma follows the same pattern than the one given for Lemma 4.8. We also notice that at some intermediate steps of the proof we consider the logic $\vdash_{\mathbf{L_n^c}}$ where there are canonical constants.

---

[23] Hint: To see that this matrix is not a model of $(\mathsf{R}_{0.5})$ it is enough to check that the first of the above stated theorems derivable from $(\mathsf{R}_{0.5})$ is not valid in this matrix. A witness of this non validity is shown for instance taking an evaluation where $p$ and $q$ are evaluated as $(0.5, 0)$.

[24] Hint: To see that this matrix is not a model of $(\mathsf{R}_1)$ we can consider the first of its stated theorems. This theorem fails when we evaluate $p$ as $(0.5, 0)$.

LEMMA 5.6 (Truth Lemma)
Let $\mathsf{L}$ be a modal logic set over $\mathbf{L_n}$. The canonical Kripke model $\mathfrak{M}_{\mathbf{can}}(\mathsf{L})$ satisfies that $V_{can}(\varphi, v) = v(\varphi)$ for every formula $\varphi$ and every world $v$.

PROOF. The beginning of the this proof is the same one that was given for Lemma 4.8, except for the fact that now there are no canonical constants in the formulas, until Claim 4.10 (see Footnote 21). Now we continue the proof from this point. It is worth pointing out that in the first of these claims we are considering the conservative expansion $\vdash_{\mathbf{L_n^c}}$ with canonical constants.

CLAIM 5.7
There is an $m \in \omega$, formulas $\psi_1, \ldots, \psi_m$ and elements $d_1, \ldots, d_m$ such that $d_i \leqslant v(\Box\psi_i)$ for every $i \in \{1, \ldots, m\}$ and

$$\mathsf{L} \vdash_{\mathbf{L_n^c}} (\bigwedge_{1 \leqslant i \leqslant m} (\overline{d_i} \to \psi_i)) \to (\overline{a} \to \varphi).$$

PROOF OF CLAIM. This was proved inside Lemma 4.8.    Q.E.D. (Claim)

CLAIM 5.8
For every $b \in \mathbf{L}_n$, it holds that

$$\mathsf{L} \vdash_{\mathbf{L_n}} (\bigwedge_{1 \leqslant i \leqslant m} \eta_{d_i \odot b}(\psi_i)) \to \eta_{a \odot b}(\varphi).$$

PROOF OF CLAIM. Let $b$ be an element of $\mathbf{L}_n$. To prove this claim we take a non-modal homomorphism $h$ such that $h[\mathsf{L}] = \{1\}$. We have to prove that $h((\bigwedge_{1 \leqslant i \leqslant m} \eta_{d_i \odot b}(\psi_i)) \to \eta_{a \odot b}(\varphi)) = 1$. Since strongly characterizing formulas only take values in $\{0, 1\}$ it is enough to prove that if $h(\eta_{d_1 \odot b}(\psi_1)) = \ldots = h(\eta_{d_m \odot b}(\psi_m)) = 1$ then $h(\eta_{a \odot b}(\varphi)) = 1$. Hence, let us assume that $h(\eta_{d_i \odot b}(\psi_i)) = 1$ for every $i \in \{1, \ldots, m\}$. Then, $d_i \odot b \leqslant h(\psi_i)$ for every $i \in \{1, \ldots, m\}$. Therefore, $b \leqslant h((\overline{d_1} \to \psi_1) \wedge \ldots \wedge (\overline{d_m} \to \psi_m))$. By the previous claim we get that $b \leqslant h(\overline{a} \to \varphi)$, i.e., $a \odot b \leqslant h(\varphi)$. Thus, $h(\eta_{a \odot b}(\varphi)) = 1$, and so this claim is proved.    Q.E.D. (Claim)

Using that $\mathsf{L}$ is closed under $\vdash_{\mathbf{L_n}}$, we get that for every $b \in \mathbf{L}_n$, it holds that $(\bigwedge_{1 \leqslant i \leqslant m} \eta_{d_i \odot b}(\psi_i)) \to \eta_{a \odot b}(\varphi)$ belongs to $\mathsf{L}$. Therefore, $(\eta_{d_1}(\Box\psi_1) \wedge \ldots \wedge \eta_{d_m}(\Box\psi_m)) \to \eta_a(\Box\varphi) \in \mathsf{L}$ by the rule $(\mathsf{R}_a^{d_1 \ldots d_m})$. Finally, using that $v((\eta_{d_1}(\Box\psi_1) \wedge \ldots \wedge \eta_{d_m}(\Box\psi_m)) \to \eta_a(\Box\varphi)) = 1 \to v(\eta_a(\Box\varphi))$ and that $v[\mathsf{L}] = \{1\}$ we obtain that $v(\eta_a(\Box\varphi)) = 1$, i.e., $a \leqslant v(\Box\varphi)$. This finishes the proof.    ∎

THEOREM 5.9 (Axiomatization of $\mathbf{\Lambda}(l, \mathsf{Fr}, \mathbf{L_n})$)
Let $\mathsf{L}$ be the smallest modal logic set over $\mathbf{L_n}$. Then,

1. $v$ is a semantic modal valuation iff $v$ is a point of $\mathfrak{M}_{\mathbf{can}}(\mathsf{L})$.

2. $\Gamma \vdash^l_{\mathsf{Fr}(\mathbf{L_n})} \varphi$   iff   $\mathsf{L} \cup \Gamma \vdash_{\mathbf{L_n}} \varphi$,   for every set $\Gamma \cup \{\varphi\}$ of formulas.

3. $\mathsf{L} = \Lambda(\mathsf{Fr}, \mathbf{L_n})$.

4. $\mathbf{\Lambda}(l, \mathsf{Fr}, \mathbf{L_n})$ is axiomatized by (i) $\mathsf{L}$ as a set of axioms, and (ii) Modus Ponens rule as its unique rule.

PROOF. These items are proved like in Theorem 4.11.    ∎

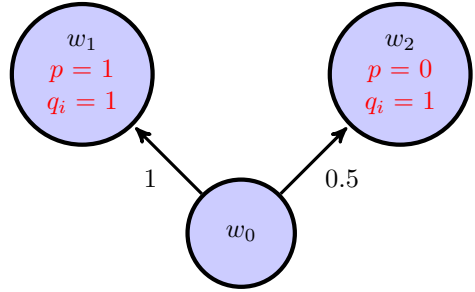

Fig. 4. A non crisp Kripke model over $\mathbf{L_n}$ where (K) is valid

There are still several important open problems around the modal logics over $\mathbf{L_n}$. The main one perhaps is to axiomatize $\mathbf{\Lambda}(g, \mathsf{Fr}, \mathbf{L_n})$. A related question is *whether the global modal logic is the smallest consequence relation extending the local one that is closed under the Monotonicity rule.* If the answer to this question is positive then we would get that $\mathbf{\Lambda}(g, \mathsf{Fr}, \mathbf{L_n})$ is axiomatized by the axioms and rules in Table 5 but restricting the use of the rules $(\mathsf{R}_a)$ only to the case that the assumption is a theorem.

## 5.2 Completeness of $\mathbf{\Lambda}(l, \mathsf{CFr}, L_\mathbf{n})$ and $\mathbf{\Lambda}(g, \mathsf{CFr}, L_\mathbf{n})$

To finish the study of the modal logics over $\mathbf{L_n}$ given by all three basic classes of frames now we devote our attention to the case of crisp frames (we do not need to pay attention to idempotent frames because in $\mathbf{L_n}$ the only idempotent elements are precisely those in $\{0,1\}$). Although the logic $\mathbf{\Lambda}(g, \mathsf{CFr}, \mathbf{L_n})$ was axiomatized in [39] by the calculus in Table 1, we will deal with this issue in the present section in order to show how this result can be obtained as an application of what has been done in Section 5.1. This approach gives us a different proof from the one in [39], and in our opinion it gives more insight.

In case that there are canonical constants we saw in Section 4.3 that it is enough to add the normality axiom (K) to $\Lambda(\mathsf{Fr}, \mathsf{L}_\mathbf{n}^\mathsf{c})$ in order to axiomatize the modal logic of crisp frames (remember that in this context they coincide with idempotent ones). However, the following result shows that this is not the case without canonical constants in the language.

### Lemma 5.10
Let $\mathsf{L}$ be the smallest modal logic set over $\mathbf{L_n}$ such that it contains the normality axiom (K). Then, $\mathsf{L} \subsetneq \Lambda(\mathsf{CFr}, \mathbf{L_n})$.

Proof. The inclusion is obvious. To see it is strict it is enough to prove that $(\Box p \oplus \Box p) \leftrightarrow \Box(p \oplus p) \notin \mathsf{L}$ because this formula belongs to $\Lambda(\mathsf{CFr}, \mathbf{L_n})$. Let us consider the Kripke model $\mathfrak{M}$ given in Figure 4. Next we see that the formula scheme (K) is valid[25] in $\mathfrak{M}$.

### Claim 5.11
For every $\varphi$ and $\psi$, it holds that $\mathfrak{M} \models^1 \Box(\varphi \to \psi) \to (\Box\varphi \to \Box\psi)$.

---

[25] Of course (cf. Theorems 4.13 and 4.14) if we allow canonical constants in the language then $\mathfrak{M} \not\models^1$ (K). For example, $w_0(\Box((\overline{0.5} \oplus p) \to p) \to (\Box(\overline{0.5} \oplus p) \to \Box p)) = 1 \to (1 \to 0.5) \neq 1$.

PROOF OF CLAIM. This is clear for points $w_1$ and $w_2$, and so it suffices to prove that $\mathfrak{M}, w_0 \models^1 \Box(\varphi \to \psi) \to (\Box\varphi \to \Box\psi)$. The idea behind this proof is that the evaluation of an arbitrary formula (remember that there are no canonical constants) in points $w_1$ and $w_2$ always belongs to $\{0, 1\}$. Thus, we get that $w_1(\varphi \odot (\varphi \to \psi)) = w_1(\varphi \wedge (\varphi \to \psi))$ and that $w_2(\varphi \odot (\varphi \to \psi)) = w_2(\varphi \wedge (\varphi \to \psi))$. Therefore, $w_0(\Box\varphi \odot \Box(\varphi \to \psi)) \leqslant w_0(\Box\varphi \wedge \Box(\varphi \to \psi)) = w_0(\Box(\varphi \wedge (\varphi \to \psi))) = (1 \to w_1(\varphi \wedge (\varphi \to \psi))) \wedge (0.5 \to w_2(\varphi \wedge (\varphi \to \psi))) = (1 \to w_1(\varphi \odot (\varphi \to \psi))) \wedge (0.5 \to w_2(\varphi \odot (\varphi \to \psi))) = w_0(\Box(\varphi \odot (\varphi \to \psi))) \leqslant w_0(\Box\psi)$. This finishes the proof of this claim.                                                                         Q.E.D. (Claim)

Hence, we know that $\{\varphi : \mathfrak{M} \models^1 \varphi\}$ contains the scheme (K) and that it is closed under the Modus Ponens rule, the Monotonicity rule and the rules ($\mathsf{R}_a$) (see Lemma 5.1). Although $\{\varphi : \mathfrak{M} \models^1 \varphi\}$ is not a modal logic set (because it is not closed under substitution since $\Box q_0$ is in this set while $\Box p$ not) we can easily see from the previous sentence that $\mathsf{L} \subseteq \{\varphi : \mathfrak{M} \models^1 \varphi\}$. Thus, using that $\mathfrak{M}, w_0 \not\models^1 (\Box p \oplus \Box p) \leftrightarrow \Box(p \oplus p)$ we get that $(\Box p \oplus \Box p) \leftrightarrow \Box(p \oplus p) \notin \mathsf{L}$. ∎

Next we see in the following three results that what we need to add to $\Lambda(\mathsf{Fr}, \mathbf{L_n})$ are the formulas $\eta_a(\Box\varphi) \to \Box\eta_a(\varphi)$ where $a \in \mathsf{L}_n \setminus \{0\}$. The first result will give us the soundness (remember that $\eta_a(p)$ is non decreasing) and the second one the completeness. Lemma 5.12 is already stated in [39], but without a proof.

LEMMA 5.12 ([39, Proposition 2.2])
Let us assume that $\delta(p)$ is a non-modal formula and is non decreasing over $\mathbf{L_n}$. Then, $\delta(\Box p) \leftrightarrow \Box\delta(p)$ is valid on the class of crisp frames.

PROOF. To prove this lemma let us consider a valuation $V$ on a crisp frame and a world $w$ in the same frame. We notice that the fact that $\delta(p)$ is non decreasing guarantees that, by continuity, it commutes with arbitrary meets. Then,

$$V(\delta(\Box p), w) =$$
$$\delta^{\mathbf{L_n}}(V(\Box p, w)) =$$
$$\delta^{\mathbf{L_n}}\left(\bigwedge\{V(p, w') : w' \in W, R(w, w') = 1\}\right) =$$
$$\bigwedge\{\delta^{\mathbf{L_n}}(V(p, w')) : w' \in W, R(w, w') = 1\} =$$
$$\bigwedge\{V(\delta(p), w') : w' \in W, R(w, w') = 1\} =$$
$$V(\Box\delta(p), w).$$                                                                                      ∎

LEMMA 5.13
Let $\mathsf{L}$ be the smallest modal logic set over $\mathbf{L_n}$ containing $\eta_a(\Box\varphi) \to \Box\eta_a(\varphi)$ for every $a \in \mathsf{L}_n \setminus \{0\}$. Then, the frame of $\mathfrak{M_{can}}(\mathsf{L})$ is a crisp frame.

PROOF. Let $a$ be $R_{can}(v_1, v_2)$ for some points $v_1$ and $v_2$ of the canonical Kripke model. We have to prove that $a \in \{0, 1\}$. Let us assume that $a \neq 0$. Using that $\mathbf{L_n}$ is a finite chain it is clear from the definition of $R_{can}(v_1, v_2)$ that there is some formula $\varphi$ such that $a = v_1(\Box\varphi) \to v_2(\varphi)$. Let $b$ be $v_1(\Box\varphi)$. Without loss of generality we can assume that $b \neq 0$ (because if $b = 0$ then $a = b \to v_2(\varphi) = 1$). Then, $v_1(\eta_b(\Box\varphi)) = 1$. Hence, using that $\eta_b(\Box\varphi) \to \Box\eta_b(\varphi) \in \mathsf{L}$ we get that $v_1(\Box\eta_b(\varphi)) = 1$. By the Truth Lemma we obtain that $1 = v_1(\Box\eta_b(\varphi)) \leqslant R_{can}(v_1, v_2) \to v_2(\eta_b(\varphi)) = a \to v_2(\eta_b(\varphi))$.

Since $\eta_b$ only takes values in $\{0,1\}$ we get that $v_2(\eta_b(\varphi)) = 1$, i.e., $b \leqslant v_2(\varphi)$, i.e., $v_1(\Box\varphi) \leqslant v_2(\varphi)$. Therefore, $a = v_1(\Box\varphi) \leqslant v_2(\varphi) = 1$. ∎

**Theorem 5.14** (Axiomatization of $\mathbf{\Lambda}(l, \mathsf{CFr}, \mathbf{L_n})$ and $\mathbf{\Lambda}(g, \mathsf{CFr}, \mathbf{L_n})$)
Let $\mathsf{L}$ be the smallest modal logic set over $\mathbf{L_n}$ containing $\eta_a(\Box\varphi) \to \Box\eta_a(\varphi)$ for every $a \in \mathrm{L}_n \setminus \{0\}$. Then,

1. $v$ is a semantic modal valuation arising from a crisp frame iff $v$ is a point of $\mathfrak{M}_{\mathsf{can}}(\mathsf{L})$.

2. $\Gamma \vdash^l_{\mathsf{CFr}(\mathbf{L_n})} \varphi$   iff   $\mathsf{L} \cup \Gamma \vdash_{\mathbf{L_n}} \varphi$,   for every set $\Gamma \cup \{\varphi\}$ of formulas.

3. $\mathsf{L} = \Lambda(\mathsf{CFr}, \mathbf{L_n})$.

4. $\mathbf{\Lambda}(l, \mathsf{CFr}, \mathbf{L_n})$ is axiomatized by (i) $\mathsf{L}$ as a set of axioms, and (ii) Modus Ponens rule as its unique rule.

5. $\mathbf{\Lambda}(g, \mathsf{CFr}, \mathbf{L_n})$ is axiomatized by (i) $\mathsf{L}$ as a set of axioms, and (ii) Modus Ponens rule and Necessity rule.

**Proof.** The proof is like the ones given for Theorems 4.11 and 4.22. ∎

To finish this section we prove a lemma showing the connection of our previous axiomatization with the one found by Hansoul and Teheux in [39]. Lemma 5.16 is already proved in the paper [39] and it says that it is enough to add the axioms $\Box(\varphi \odot \varphi) \leftrightarrow (\Box\varphi \odot \Box\varphi)$ and $\Box(\varphi \oplus \varphi) \leftrightarrow (\Box\varphi \oplus \Box\varphi)$ in order to get the formulas $\eta_a(\Box\varphi) \to \Box\eta_a(\varphi)$. It is worth pointing out that this lemma talks about the smallest set and not about the smallest modal logic set.

**Remark 5.15** (Trick used in the proof of Lemma 5.16)
Let us consider $\tau_1(p) := p \odot p$ and $\tau_2(p) := p \oplus p$. It is known (see [39, Definition 5.3] and [50]) that for every $a \in [0,1]$ if $a$ is[26] a finite sum of negative powers of 2, then in $[\mathbf{0}, \mathbf{1}]_{\underline{\mathbb{L}}}$ there is a characterizing formula $\eta_a(p)$ for the interval $[a, 1]$ such that $\eta_a$ is obtained as a composition of $\tau_1$ and $\tau_2$. From here, composing several times $\tau_1$ with the previous characterizing formula, it is obvious that in $\mathbf{L_n}$ for every $a \in \mathrm{L}_n \setminus \{0, 1\}$ there is a strongly characterizing formula $\eta_a(p)$ for the interval $[a, 1]$ such that $\eta_a$ is obtained as a composition of $\tau_1$ and $\tau_2$. And so, it is clear that the same holds for every $a \in \mathrm{L}_n \setminus \{0\}$. ⊣

**Lemma 5.16** (cf. [39, Proposition 6.3])
1. $\Lambda(\mathsf{CFr}, \mathbf{L_n})$ is the smallest set closed under the axioms, but replacing (K) with (MD), and rules given in Table 1.
2. $\Lambda(\mathsf{CFr}, \mathbf{L_n})$ is the smallest set closed under the axioms and rules given in Table 1.

**Proof.** 1): Let $\mathsf{L}$ be the smallest set closed under the axioms, but replacing (K) with (MD), and rules given in Table 1. By Theorem 5.14 we know that $\Lambda(\mathsf{CFr}, \mathbf{L_n})$ is the smallest modal logic set over $\mathbf{L_n}$ containing $\eta_a(\Box\varphi) \to \Box\eta_a(\varphi)$ for every $a \in \mathrm{L}_n \setminus \{0\}$. Therefore, it suffices to prove that $\mathsf{L}$ is a modal logic set and contains the formulas $\eta_a(\Box\varphi) \to \Box\eta_a(\varphi)$ for every $a \in \mathrm{L}_n \setminus \{0\}$.

First of all we prove that $\mathsf{L}$ contains the formulas $\eta_a(\Box\varphi) \to \Box\eta_a(\varphi)$. Let $a \in \mathrm{L}_n \setminus \{0\}$. By Remark 5.15 we can assume that $\eta_a(p)$ is obtained as a composition of $\tau_1(p)$ and $\tau_2(p)$. Using that $\tau_1(\Box p) \leftrightarrow \Box\tau_1(p) \in \mathsf{L}$ and that $\tau_2(\Box p) \leftrightarrow \Box\tau_2(p) \in \mathsf{L}$ it is obvious that $\eta_a(\Box\varphi) \rightarrow \Box\eta_a(\varphi) \in \mathsf{L}$.

---

[26] It is worth pointing out that the set $C$ of elements that are finite sums of negative powers of 2 satisfies the condition given in Footnote 18.

Finally we must check that $\mathsf{L}$ is a modal logic set. The only non trivial part is to see that $\mathsf{L}$ is closed under the rules $(\mathsf{R}_a)$, but the following derivation sketch

$$\frac{\dfrac{\dfrac{\dfrac{\dfrac{(\eta_{a_2 \odot b}(\varphi_2) \wedge \eta_{a_3 \odot b}(\varphi_3) \wedge \ldots \wedge \eta_{a_n \odot b}(\varphi_n)) \rightarrow \eta_{a \odot b}(\varphi) \text{ for every } b > \neg a}{(\eta_{a_2}(\varphi_2) \wedge \eta_{a_3}(\varphi_3) \wedge \ldots \wedge \eta_{a_n}(\varphi_n)) \rightarrow \eta_a(\varphi)}}{\Box(\eta_{a_2}(\varphi_2) \wedge \eta_{a_3}(\varphi_3) \wedge \ldots \wedge \eta_{a_n}(\varphi_n)) \rightarrow \Box\eta_a(\varphi)}}{(\Box\eta_{a_2}(\varphi_2) \wedge \Box\eta_{a_3}(\varphi_3) \wedge \ldots \wedge \Box\eta_{a_n}(\varphi_n)) \rightarrow \Box\eta_a(\varphi)}}{(\eta_{a_2}(\Box\varphi_2) \wedge \eta_{a_3}(\Box\varphi_3) \ldots \wedge \eta_{a_n}(\Box\varphi_n)) \rightarrow \eta_a(\Box\varphi)}}$$

shows that the rules $(\mathsf{R}_a)$ are derivable from the theorems $\eta_a(\Box\varphi) \leftrightarrow \Box\eta_a(\varphi)$ using the Monotonicity rule and the meet distributivity axiom.

2): The proof is the same as in the previous item, just realizing that in the previous derivation sketch we can replace in the intermediate steps meet with fusion (because strongly idempotent formulas are only evaluated into $\{0, 1\}$). ∎

We stress that the first item of the previous lemma gives us a strengthing of the axiomatization provided in [39] since $(\mathsf{MD})$ holds in all Kripke frames and not oly in the crisp ones.

## 6 Concluding Remarks

This article is among the first to address the study of minimum modal logics in the context of residuated lattices, and there are still a lot of open questions. In this section we collect what in our opinion are the main open problems concerning the framework discussed in this article.

Before giving the list of open problems we want to stress that several important frameworks different than the one considered in this article have not yet been studied. Among them we want to mention at least three. The first one is the development of a general theory where $\Box$ and $\Diamond$ are simultaneously in the language. The second one concerns the theory of many-valued modal logics over classes of residuated lattices (not just one residuated lattice), which can be naturally defined as the intersection of the modal logics over each one of the members of the class. And the last one is about comparing many-valued modal logics (without canonical constants in the language) given by two different residuated lattices. One straightforward result of this kind based on direct powers is that

$$\mathbf{\Lambda}(g, \mathsf{Fr}, \mathbf{A^I}) = \mathbf{\Lambda}(g, \mathsf{Fr}, \mathbf{A}) \text{ and } \mathbf{\Lambda}(l, \mathsf{Fr}, \mathbf{A^I}) = \mathbf{\Lambda}(l, \mathsf{Fr}, \mathbf{A}).$$

The reason why this holds is that any frame over $\mathbf{A^I}$ can be seen as a family, indexed by $I$, of frames. Indeed, the same argument also shows that the first order logic given by $\mathbf{A^I}$ is also the first order logic given by $\mathbf{A}$. The authors consider that the development of a model theory for first-order many-valued logics, which is still missing in the literature, will also benefit the realm of many-valued modal logics. Another property comparing two different residuated latices is that if $\mathbf{A}$ is completely embeddable into $\mathbf{B}$, then $\mathbf{\Lambda}(l, \mathsf{Fr}, \mathbf{B}) \leqslant \mathbf{\Lambda}(l, \mathsf{Fr}, \mathbf{A})$ and $\mathbf{\Lambda}(g, \mathsf{Fr}, \mathbf{B}) \leqslant \mathbf{\Lambda}(g, \mathsf{Fr}, \mathbf{A})$. It is not enough to have an embedding (see Footnote 10), it has to preserve arbitrary meets and joins.

Now we list the main open problems concerning our framework. In the formulation of these problems **A** is an arbitrary complete residuated lattice (maybe non finite).

**Problem 1.** How can we expand the consequence relation $\mathbf{\Lambda}(\mathbf{A})$ in order to get $\mathbf{\Lambda}(g, \mathsf{Fr}, \mathbf{A})$? And for $\mathbf{\Lambda}(l, \mathsf{Fr}, \mathbf{A})$? And also the same problems with canonical constants in the language.

**Problem 2.** Is there any **A** such that the set of theorems of $\vdash_{\mathbf{A}}$ is recursively axiomatizable while the set $\Lambda(\mathsf{Fr}, \mathbf{A})$ is not?

**Problem 3.** What is the computational complexity of the set $\Lambda(\mathsf{Fr}, \mathbf{A})$? And the same question for the other classes of frames.

**Problem 4.** Is $\mathbf{\Lambda}(g, \mathsf{Fr}, \mathbf{A})$ the smallest consequence relation extending $\mathbf{\Lambda}(l, \mathsf{Fr}, \mathbf{A})$ that is closed under the Monotonicity rule? And the same problem with canonical constants.

**Problem 5.** How can we characterize those **A**'s such that $\Lambda(\mathsf{IFr}, \mathbf{A}) = \Lambda(\mathsf{CFr}, \mathbf{A})$?

**Problem 6.** How can we axiomatize the set $\Lambda(\mathsf{Fr}, [\mathbf{0}, \mathbf{1}]_{\mathrm{L}})$ and its related modal logics? And what about $\Lambda(\mathsf{CFr}, [\mathbf{0}, \mathbf{1}]_{\mathrm{L}})$?

**Problem 7.** Is there some axiomatization of the set $\Lambda(\mathsf{Fr}, \mathbf{L_n})$ simpler than the one given in Table 5?

To finish the article we make some comments about the Problems 1, 3, 6 and 7. From an intuitive point of view it seems that Problem 1 could be rewritten as wondering how we can expand an axiomatization of $\vdash_{\mathbf{A}}$ in order to get an axiomatization of $\mathbf{\Lambda}(g, \mathsf{Fr}, \mathbf{A})$. However, this last formulation has some ambiguity (for finite algebras, like the ones studies in this article, this ambiguity disappears). The ambiguity comes from the fact that if $\vdash_{\mathbf{A}}$ is not finitary then we know that there are no strongly complete axiomatizations (at least in the context of finitary formal proof systems). Hence, we can think on this problem in two different ways. A first problem is the previous question relaxing the meaning of an axiomatization to sets with finite assumptions. And a second problem is to change the paradigm and to use infinitary proof systems where we still have chances to get strongly complete axiomatizations for non finitary consequence relations. Both questions are also interesting.

Next we focus on Problem 3. As far as the authors are aware the only knwon result about this topic is the computational complexity of $\Lambda(\mathsf{Fr}, [\mathbf{0}, \mathbf{1}]_{\mathbf{G}})$. In [47] it is shown that this set is PSPACE-complete. On the other hand, it is easy to prove that if **A** is a finite residuated lattice then the sets $\Lambda(\mathsf{Fr}, \mathbf{A})$ and $\Lambda(\mathsf{Fr}, \mathbf{A^c})$ are decidable. This can be easily proved using a filtration method[27]. It is also worth pointing out that in [31] it is proved that $\Lambda(\mathsf{Fr}, [\mathbf{0}, \mathbf{1}]_{\mathrm{L}})$ is decidable.

Let us give some ideas about Problem 6. One of the few things that it is known about the set $\Lambda(\mathsf{Fr}, [\mathbf{0}, \mathbf{1}]_{\mathrm{L}})$ is that it coincides with $\bigcap_n \Lambda(\mathsf{Fr}, \mathbf{L_n})$. This result is a consequence of [30, Theorem 5.4.30]. Unfortunately, the proof given in Section 5.1 for the case of finite MV chains does not seem easy to be adapted to the case of the standard Łukasiewicz algebra. At least there are two clear difficulties when one tries to generalize the proof. The first one concerns about writing valid rules that play the same role than $(\mathsf{R_a})$: the problem is that in the standard Łukasiewicz algebra there are no strongly characterizing formulas. And the second difficulty comes from the fact that the non-modal logic is not finitary (this is used to get the first claim inside the proof of Lemma 5.6).

---

[27] Indeed, the filtration method can also be used to prove that the logic defined using simultaneously $\Box$ and $\Diamond$ is also decidable

At this moment it is not clear for us how to overcome the first difficulty. But about the second one we think that one possible way is to use infinitary proof systems. This method has been used in [39] to show an strong axiomatization of $\Lambda(g, \mathsf{CFr}, [\mathbf{0}, \mathbf{1}]_{\mathrm{L}})$ (see Section 3.3), but even in the classical modal setting the method has been successfully considered (see [44]). In our opinion the addition of infinitary rules may help to obtain a completeness proof based on a canonical model construction. Although in the literature there are some reasons to believe that the canonical model construction should not help to settle the completeness for the case of the standard Łukasiewicz algebra (see [26, 45]), these limitations may not apply to the case of using a canonical model construction based on infinitary proof systems. We think, roughly speaking, that infinitary proof systems may help to overtake these limitations.

Finally, we focus on Problem 7 about the search of a more transparent axiomatization for $\Lambda(\mathsf{Fr}, \mathbf{L_n})$. That is, *can we simplify the axiomatization given in Table 5?* The authors can imagine at least two different strategies that may help to simplify the axiomatization. A first one would be based on understanding better which non-modal formulas are always evaluated in $[a, 1] \cap \mathrm{L}_n$ for a certain value $a \in \mathrm{L}_n$. We feel that the understanding[28] of these non-modal formulas may help to solve this open question due to the connection with the rules[29] ($\mathsf{R}_a$). Next we illustrate this connection with an straightforward remark: if $\varphi$ is always evaluated in $[a, 1]$ and $\eta_a$ is a (maybe not strongly) characterizing formula of the interval $[a, 1]$ then the formula $\eta_a(\Box\varphi)$ belongs to $\Lambda(\mathsf{Fr}, \mathbf{L_n})$ (indeed it also belongs to $\Lambda(\mathsf{Fr}, [\mathbf{0}, \mathbf{1}]_{\mathrm{L}})$). For example, using $a = 0.5$ we can get that $\Box(p \vee \neg p) \oplus \Box(p \vee \neg p)$ is a theorem. However, nowadays it is still unclear to the authors whether this research line could really help to simplify Table 5.

The other strategy is to visualize a many-valued accessibility relation $R$ as a non-increasing sequence $\langle R_a : a \in \mathrm{L}_n \setminus \{0\} \rangle$ of crisp accessibility relations. By being non-increasing we refer to the fact that if $a \leqslant a'$ then $R_{a'} \subseteq R_a$. It is clear that every many-valued accessibility relation $R$ can be canonically represented in this way using $R_a := \{(w, w') : R(w, w') \geqslant a\}$; and on the other hand, every non-increasing sequence $\langle R_a : a \in \mathrm{L}_n \setminus \{0\} \rangle$ of crisp accessibility relations is canonically represented by the many-valued accessibility relation

$$R(w, w') := \begin{cases} \min\{a \in \mathrm{L}_n \setminus \{0\} : (w, w') \in R_a\}, & \text{if } (w, w') \in \bigcup_a R_a \\ 0, & \text{otherwise.} \end{cases}$$

This suggests that $\Lambda(\mathsf{Fr}, \mathbf{L_n})$ has canonically associated a multi-modal $\{\Box_a : a \in \mathrm{L}_n \setminus \{0\}\}$ logic where each modality $\Box_a$ is semantically interpreted by the crisp accessibility relation $R_a$. The reason why this perspective could possibly help to simplify the axiomatization in Table 5 is that for the case of crisp frames there are really simple axiomatizations of the many-valued modal logic over $\mathbf{L_n}$ due to Hansoul and Teheux (see Lemma 5.16). Indeed, some previous attempts to axiomatize $\Lambda(\mathsf{Fr}, \mathbf{L_n})$ done by some of the authors were based on this idea (see [8]).

---

[28] This research has been started by one of the authors in the recent manuscript [7].

[29] We are mainly thinking in the informal version stated just before Lemma 5.1.

| | 0 | 0.25 | 0.5 | 0.75 | 1 |
|---|---|---|---|---|---|
| $n$ | 1 | 0.75 | 0.25 | 0.25 | 0 |

$$x \odot y := \left\{ \begin{array}{ll} 0, & \text{if } x \leqslant n(y) \\ x \wedge y & \text{if } n(y) < x \end{array} \right. \qquad x \to y := \left\{ \begin{array}{ll} 1, & \text{if } x \leqslant y \\ n(x) \vee y & \text{if } y < x \end{array} \right.$$

FIG. 5. A weak nilpotent minimum algebra

# A   Some remarks about the non-modal logics

## A.1   *The (non-modal) logic of a residuated lattice*

The purpose of this appendix is to summarize several known results about the non modal logic $\mathbf{\Lambda}(\mathbf{A})$, which we remind is defined by

$$\Gamma \vdash_{\mathbf{A}} \varphi \quad \Longleftrightarrow \quad \forall h \in \mathrm{Hom}(\mathbf{Fm}, \mathbf{A}), \text{ if } h[\Gamma] \subseteq \{1\} \text{ then } h(\varphi) = 1. \tag{2.7}$$

The fact that this logic is introduced using only one residuated lattice supposes a difference with the standard framework where logics are defined using a variety of residuated lattices (cf. [30, 25]). For example, it may be the case that the Local Deduction Theorem fails for the logic $\vdash_{\mathbf{A}}$ (see Example A.4). Hence, the reader must be careful about not confusing this framework with the standard one, since there are several differences.

PROPOSITION A.1 ([16, 24])
The logic $\mathbf{\Lambda}(\mathbf{A})$ is algebraizable (in the sense of Blok and Pigozzi).

PROPOSITION A.2 ([9, 16])

1. The logics $\mathbf{\Lambda}(\mathbf{A})$ and $\mathbf{\Lambda}(\mathbf{B})$ share the same theorems iff $\mathbf{A}$ and $\mathbf{B}$ generate the same variety.

2. The logics $\mathbf{\Lambda}(\mathbf{A})$ and $\mathbf{\Lambda}(\mathbf{B})$ share the same finitary deductions iff $\mathbf{A}$ and $\mathbf{B}$ generate the same quasivariety.

3. The logics $\mathbf{\Lambda}(\mathbf{A})$ and $\mathbf{\Lambda}(\mathbf{B})$ coincide iff the closure under the operators $\mathbb{I}, \mathbb{S}, \mathbb{P}_{\sigma-f}$ (where $\mathbb{P}_{\sigma-f}$ denotes the operator of reduced products over countably complete filters) of the classes $\{\mathbf{A}\}$ and $\{\mathbf{B}\}$ are the same class.

PROPOSITION A.3

1. If $\mathbf{A}$ is finite, then $\mathbf{\Lambda}(\mathbf{A})$ is finitary. That is, if $\Gamma \vdash_{\mathbf{A}} \varphi$ then there is a finite $\Gamma' \subseteq \Gamma$ such that $\Gamma' \vdash_{\mathbf{A}} \varphi$.

2. If $\mathbf{A}$ is such that 1 is join irreducible (i.e., if $a \vee b = 1$ then either $a = 1$ or $b = 1$), then $\mathbf{\Lambda}(\mathbf{A})$ admits proofs by cases. That is, $\Gamma, \varphi \vee \psi \vdash_{\mathbf{A}} \delta$ iff it holds that $\Gamma, \varphi \vdash_{\mathbf{A}} \delta$ and that $\Gamma, \psi \vdash_{\mathbf{A}} \delta$.

3. If $\mathbb{Q}(\mathbf{A})$ is a variety, then[30] the Local Deduction Theorem holds. That is[31], $\gamma_1, \ldots, \gamma_n, \gamma_{n+1} \vdash_{\mathbf{A}} \varphi$ iff there is some $m \in \omega$ such that $\gamma_1, \ldots, \gamma_n \vdash_{\mathbf{A}} \gamma_{n+1}^m \to \varphi$.

EXAMPLE A.4
We point out that the hypothesis of the third item in Proposition A.3 cannot be deleted. For instance, the weak nilpotent minimum algebra (see [19]) over the domain $\{l/4 : 0 \leqslant l \leqslant 4, l \in \omega\}$ (with the order induced by the reals) that is given in Figure 5 is an example witnessing the failure of the Local Deduction Theorem. It fails because $\neg\neg p \vdash_{\mathbf{A}} p$, while evaluating $p$ as 0.5 we get that for every $m \in \omega$, it holds that $\not\vdash_{\mathbf{A}} (\neg\neg p)^m \to p$. ⊣

---

[30] Indeed, these two conditions are equivalent. The other direction follows from the fact that the Local Deduction Theorem gives us a method to convert quasiequations into equations.

[31] We point out that in this article by the Local Deduction Theorem we mean its finitary version. Indeed, we notice that for the particular case $n = 1$ we get an equivalent statement.

REMARK A.5 (Modus Ponens as a unique rule)
The third statement in the previous proposition says in particular that if $\mathbb{Q}(\mathbf{A})$ is a variety then the finitary deductions of $\mathbf{\Lambda}(\mathbf{A})$ can be obtained using only one inference rule: Modus Ponens. The previous example shows that in general it is not true that $\mathbf{\Lambda}(\mathbf{A})$ can be axiomatized using Modus Ponens as a unique rule; maybe more rules have to be added.                                                    ⊣

Is there some way to characterize those $\mathbf{A}$'s such that $\mathbb{Q}(\mathbf{A})$ is a variety? Thanks to the congruence extension property [25, Lemma 3.57] (it implies $\mathbb{SH} = \mathbb{HS}$) together with Jónsson's Lemma [42] it is obvious that[32] $\mathbb{Q}(\mathbf{A})$ is a variety iff $\mathbb{HP}_U(\mathbf{A}) \subseteq \mathbb{Q}(\mathbf{A})$. In the particular case that $\mathbf{A}$ is finite (cf. [51, Theorem 1.1.1]) we get that $\mathbb{Q}(\mathbf{A})$ is a variety iff $\mathbb{H}(\mathbf{A}) \subseteq \mathbb{Q}(\mathbf{A})$. Therefore, for every finite residuate lattice $\mathbf{A}$, if $\mathbf{A}$ is simple then $\mathbb{Q}(\mathbf{A})$ is a variety. Next we state other cases where we know for sure that the generated quasivariety is a variety.

PROPOSITION A.6
Let $\mathbf{A}$ be a finite BL chain. Then, $\mathbb{Q}(\mathbf{A})$ is a variety. And $\mathbf{\Lambda}(\mathbf{A})$ can be axiomatized using Modus Ponens as a unique rule.

PROOF. It is enough to see that if $F$ is a congruence filter then $\mathbf{A}/F$ is embeddable into $\mathbf{A}$. The reader can easily check that

$$
\begin{array}{rccl}
i: & \mathbf{A}/F & \longrightarrow & \mathbf{A} \\
& x/F & \longmapsto & \left\{ \begin{array}{ll} 1, & \text{if } x \in F \\ x, & \text{if } x \notin F \end{array} \right.
\end{array}
$$

is well defined and it is an embedding. The reason why this map is well defined is because if $a \notin F$, then $a/F$ is a singleton.                                                                 ∎

EXAMPLE A.7
The previous statement fails in the case of arbitrary BL chains. Let us consider for instance $\mathbf{A}$ as the (simple) subalgebra of standard Łukasiewicz algebra $[\mathbf{0}, \mathbf{1}]_{\text{Ł}}$ with domain $A = \{x \in [0,1] : x \text{ is a rational number with denominator odd}\}$. Then, $p \leftrightarrow \neg p \vdash_{\mathbf{A}} q$ while for every $m \in \omega$, it holds that $\nvdash_{\mathbf{A}} (p \leftrightarrow \neg p)^m \to q$.                                                    ⊣

A natural question that appears here is whether there are complete BL chains $\mathbf{A}$ such that the quasivariety generated by $\mathbf{A}$ is not a variety. In [18, Corollary 12] this question is negatively answered for the case of standard BL algebras. Indeed, a similar proof can be used to answer the first question.

PROPOSITION A.8 (cf. [18, Corollary 7])
Let $\mathbf{A}$ be a complete BL chain. Then, $\mathbb{Q}(\mathbf{A})$ is a variety, and so $\mathbf{\Lambda}(\mathbf{A})$ can be axiomatized using Modus Ponens as a unique rule.

Throughout this article we have assumed that we already know an axiomatization of $\mathbf{\Lambda}(\mathbf{A})$ in order to see how to expand it for the modal logic. Unfortunately, in the literature there is no general result stating how to axiomatize an arbitrary logic $\mathbf{\Lambda}(\mathbf{A})$. Nevertheless, it is worth pointing out that for a lot of logics of the form $\mathbf{\Lambda}(\mathbf{A})$ we already know axiomatizations:

• If $\mathbf{A}$ is a finite BL algebra, then $\mathbf{\Lambda}(\mathbf{A})$ has been axiomatized, among other logics, in [1, Theorem 5.1] using as a rule only Modus Ponens. And the axiomatization has a finite number of axioms.

• If $\mathbf{A}$ is a standard BL algebra, then the finitary deductions[33] of $\mathbf{\Lambda}(\mathbf{A})$ are axiomatized in [20] (see also [38]). The axiomatization only uses one rule: Modus Ponens. And we notice that this axiomatization has a finite number of axioms.

REMARK A.9 (Recursive and finite axiomatizations)
Since it is known that the first order logic over any finite algebra is recursively axiomatizable [52, 2] it is clear that $\mathbf{\Lambda}(\mathbf{A})$ is recursively axiomatizable for every finite residuated lattice $\mathbf{A}$. It is worth pointing out that in general it is not true that for every finite residuated lattice **A** there is a finite[34] axiomatization of **Λ(A)**. In [17] there are counterexamples for finite Heyting algebras. However it is known that if **A** is a finite subdirectly irreducible residuated lattice[35] then **Λ(A)** is finitely axiomatizable (see [6] together with [43, Proposition 1.11]).   ⊣

---

[32] In the particular case that $\mathbf{A}$ is a BL chain, then [15, Theorem 3.8] is an strenghtening of this fact. This theorem states, with different terminology, that $\mathbb{Q}(\mathbf{A})$ is a variety iff every algebra in $\mathbb{HP}_U(\mathbf{A})$ is partially embeddable into $\mathbf{A}$.

[33] We point out that if $\mathbf{A}$ is a standard BL algebra, then $\mathbf{\Lambda}(\mathbf{A})$ is finitary iff $\mathbf{A}$ is the standard Gödel algebra $[\mathbf{0}, \mathbf{1}]_{\mathbf{G}}$ (see [18, p. 603]).

## A.2 Adding canonical constants to the non-modal logic

In this appendix we analyze the non-modal logic **Λ(A^c)** obtained by adding canonical constants. We remind that this logic is defined by the equivalence

$$\Gamma \vdash_{\mathbf{A^c}} \varphi \quad \Longleftrightarrow \quad \forall h \in \mathrm{Hom}(\mathbf{Fm}, \mathbf{A^c}), \text{ if } h[\Gamma] \subseteq \{1\} \text{ then } h(\varphi) = 1. \tag{2.8}$$

This logic is introduced following the same pattern that is used for **Λ(A)**, and hence **Λ(A^c)** is clearly a conservative expansion of **Λ(A)**. Noticing this pattern it is not surprising that Proposition A.1[36], Proposition A.2 and the two first statements of Proposition A.3 also hold if we replace **A** and **B** with **A^c** and **B^c**.

Unfortunately, the addition of canonical constants also has some undesired consequences, and so we must be really careful. These disadvantages are mainly caused by the fact that $\bar{a} \vdash_{\mathbf{\Lambda(A^c)}} 0$ whenever $a \in A \setminus \{1\}$. Among the properties that may fail are the Local Deduction Theorem and proofs by cases. If the set of idempotent elements of **A** is diferent than $\{0, 1\}$, then automatically the Local Deduction Theorem fails for **Λ(A^c)**. In particular we have that the addition of canonical constants to (even finite) Gödel algebras is destroying the Deduction Theorem. To avoid mistakes it is also worth pointing out that a consequence of the failure of the Local Deduction Theorem is that these logics are not core fuzzy logics [14, 34] and hence we cannot apply the machinery developed for core fuzzy logics as in [15].

**Proposition A.10**
Each one of the following conditions implies the next one.

1. **Λ(A^c)** satisfies the Local Deduction Theorem.

2. $\mathbb{Q}(\mathbf{A^c})$ is a variety.

3. **A** is simple.

4. The set of idempotents elements of **A** is $\{0, 1\}$.
In case **A** is finite all these conditions are equivalent.

**Proof.** $1 \Rightarrow 2$): It is trivial because the Local Deduction Theorem gives us a way to convert quasiequations into equations.

$2 \Rightarrow 3$): Let us assume that there is a non-trivial congruence filter $F$ of **A**. Let $a$ be an element of $F$ such that $a \neq 1$. Then, **A^c** satisfies the quasiequation $\bar{a} \approx 1 \Rightarrow 0 \approx 1$ while the quotient $\mathbf{A^c}/F$ does not satisfy this quasiequation. Therefore, $\mathbb{Q}(\mathbf{A^c})$ is not a variety.

$3 \Rightarrow 4$): It is trivial.

$4 \Rightarrow 1$): Assume that $A$ is finite. Let $n$ be the cardinal of $A$. Due to the assumption we know that for every element $a \neq 1$ it holds that $a^n = 0$. Using this it is obvious that the Local Deduction Theorem holds. ∎

As a consequence of this last proposition it is clear that if **A** is a finite and simple residuated lattice then the logic **Λ(A^c)** can be axiomatized using Modus Ponens as a unique rule.

**Proposition A.11**
**Λ(A^c)** admits proofs by cases     iff     $1$ is a join irreducible element of **A**.

---

PROOF. This easily follows from the fact that if $1 = a \vee b$ with $a, b \in A \setminus \{1\}$, then $\overline{a} \vdash_{\mathbf{A^c}} 0$ and $\overline{b} \vdash_{\mathbf{A^c}} 0$, while $\overline{a} \vee \overline{b} \nvdash_{\mathbf{A^c}} 0$. ∎

In the rest of this appendix we will explain a method to expand an axiomatization of $\mathbf{\Lambda(A)}$, when $\mathbf{A}$ is finite, into one of $\mathbf{\Lambda(A^c)}$. This method is not universal, but it can be applied to quite a lot of cases: the ones from Proposition A.11. As far as we know, this result about axiomatizing the consequence relation (not only its theorems) obtained by adding canonical constants is new in the literature (cf. [18]).

By the name *book-keeping axioms* we will refer as usual [30] to the equations

$$\overline{a} * \overline{b} \approx \overline{a * b} \quad \text{where } a, b \in A \text{ and } * \in \{\wedge, \vee, \odot, \rightarrow\},$$

or to the formulas

$$(\overline{a} * \overline{b}) \leftrightarrow \overline{a * b} \quad \text{where } a, b \in A \text{ and } * \in \{\wedge, \vee, \odot, \rightarrow\}.$$

The context will clarify whether we are thinking on equations or on formulas. Analogously, we will refer to

$$\bigvee_{a \in A}(x \leftrightarrow \overline{a}) \approx 1 \qquad \text{or} \qquad \bigvee_{a \in A}(p \leftrightarrow \overline{a})$$

as the *witnessing axiom*.

LEMMA A.12
Let $\mathbf{A}$ be a finite residuated lattice. The variety generated by $\mathbf{A^c}$ is axiomatized by the following list of equations:
- an equational presentation of $\mathbf{A}$,
- the witnessing and book-keeping axioms.

PROOF. It is obvious that the previous list of equations hold in $\mathbf{A^c}$. Thus, it is enough to prove that every algebra $\mathbf{B}$ (with the similarity type including the constants $\{\overline{a} : a \in A\}$) that satisfies the previous list of equations and that is subdirectly irreducible is a homomorphic image of $\mathbf{A^c}$. Let us consider $\mathbf{B}$ satisfying the previous properties. Then, it is obvious that the map

$$h: \quad \begin{array}{ccc} \mathbf{A^c} & \longrightarrow & \mathbf{B} \\ a & \longmapsto & \overline{a}^{\mathbf{B}} \end{array}$$

is a homomorphism due to the book-keeping axioms. And it is onto thanks to the witnessing axiom together with the fact that 1 is join irreducible in all subdirectly irreducible algebras (see [43, Proposition 1.4]). ∎

REMARK A.13
We point out that in this proof we have seen that all subdirectly irreducible algebras in the variety generated by $\mathbf{A^c}$ are of the form $\mathbf{A^c}/F$ for a certain congruence filter $F$. This can also be obtained as an easy consequence of Jónsson Lemma and the congruence extension property. ⊣

Next we consider the problem of axiomatizing $\mathbb{Q}(\mathbf{A^c})$ for the particular case that there is a unique coatom $k$. Although for this case it is possible to give an almost trivial proof of Theorem A.14 using Lemma A.12 together with the above displayed equivalence (where $k$ is the unique coatom)

$$\gamma_1, \ldots, \gamma_n \vdash_{\mathbf{A^c}} \varphi \qquad \text{iff} \qquad \vdash_{\mathbf{A^c}} (\gamma_1 \wedge \ldots \wedge \gamma_n) \rightarrow (\varphi \vee \overline{k}), \tag{4.4}$$

we have preferred to write down the following semantic proof since it could be easier to be generalized in a future to the general case with more than one coatom.

THEOREM A.14
Let $\mathbf{A}$ be a finite residuated lattice with a unique coatom $k \notin \{0, 1\}$. The quasivariety generated by $\mathbf{A^c}$ is axiomatized by the following list of quasiequations:
- a quasiequational presentation of $\mathbf{A}$,
- the witnessing and book-keeping axioms.
- $\overline{k} \vee x \approx 1 \Rightarrow x \approx 1$.

PROOF. It is obvious that the previous list of quasiequations hold in $\mathbf{A^c}$ (we notice that the last quasiequation is valid only because $\mathbf{A}$ has coatom $k$). Thus, it is enough to prove that every non trivial algebra $\mathbf{B}$ (with the similarity type including the constants $\{\overline{a} : a \in A\}$) that satisfies the previous list of quasiequations is a subalgebra of a direct product of copies of $\mathbf{A^c}$. Let us consider a non trivial $\mathbf{B}$ satisfying the list of quasiequations. By Lemma A.12 we know that $\mathbf{B}$ is in the variety generated by $\mathbf{A^c}$. And so, together with Remark A.13, it follows that $\mathbf{B} \subseteq \prod\{\mathbf{A^c}/F_i : i \in I\}$ where the $F_i$'s are congruence filters. By the last quasiequation in the list we know that $\overline{k}^{\mathbf{B}} \neq 1$. Therefore, $\{i \in I : F_i = \{1\}\} \neq \emptyset$. Let us consider the map

$$h : \quad \mathbf{B} \quad \longrightarrow \quad \prod\{\mathbf{A^c}/F_i : i \in I, F_i = \{1\}\}$$
$$b \quad \longmapsto \quad \langle \pi_i(b) : i \in I, F_i = \{1\}\rangle$$

We remark that $\pi_i$ refers to the $i$-th projection. By definition it is obvious that $h$ is a homomorphism. The proof is completed (because $\mathbf{A^c}/\{1\} \cong \mathbf{A^c}$) by showing that $h$ is an embedding. If this does not hold, then the kernel of the map $h$ is non trivial. So, there is an element $b \in B \setminus \{1\}$ such that $h(b) = 1$, i.e., for every $i \in I$, if $F_i = \{1\}$ then $\pi_i(b) = 1$. On the other hand, for every $i \in I$, if $F_i \neq \{1\}$ then $k \in F_i$, and hence $\pi_i(\overline{k}^{\mathbf{B}}) = 1$. Therefore, for every $i \in I$, it holds that $\pi_i(\overline{k}^{\mathbf{B}} \vee b) = 1$. Thus, $\overline{k}^{\mathbf{B}} \vee b = 1$ while $b \neq 1$. This contradicts the last quasiequation in the list. ∎

We stress that in the previous theorem we cannot replace the last quasiequation with $\overline{k} \vee 0 \approx 1 \Rightarrow 0 \approx 1$. The next example shows this.

EXAMPLE A.15
Let us consider the Gödel chain $\mathbf{G_3}$ with three elements, and let us denote its coatom by $k$ and its congruence filter $\{k, 1\}$ by $F$. It is obvious that in $\mathbf{G_3^c}/F \times \mathbf{G_3^c}$ fails the last quasiequation in Theorem A.14 (because $\overline{k} \vee (0, 1) = 1$ while $(0, 1) \neq 1$). On the other hand, $\mathbf{G_3^c}/F \times \mathbf{G_3^c}$ satisfies all quasiequations valid in $\mathbf{G_3}$, the witnessing and book-keeping axioms and the quasiequation $\overline{k} \vee 0 \approx 1 \Rightarrow 0 \approx 1$.

Next we state two trivial consequences of the previous results. We stress that the axiomatization given in Corollary A.17 only uses one rule (Modus Ponens).

COROLLARY A.16
Let $\mathbf{A}$ be a finite residuated lattice with a unique coatom $k \notin \{0, 1\}$. The logic $\mathbf{\Lambda(A^c)}$ is axiomatized by

- rules and axioms axiomatizing $\mathbf{\Lambda(A)}$,

- the witnessing and book-keeping axioms,

- the rule $\overline{k} \vee p \vdash p$.

COROLLARY A.17
Let $\mathbf{A}$ be a finite and simple residuated lattice. The logic $\mathbf{\Lambda(A^c)}$ is axiomatized by

- axioms axiomatizing $\mathbf{\Lambda(A)}$,

- the witnessing and book-keeping axioms,

- the Modus Ponens rule.

To finish the appendix we point out that it is an open problem to give, for every finite residuated lattice $\mathbf{A}$, a presentation of $\mathbb{Q}(\mathbf{A^c})$ (cf. Theorem A.14). In our opinion one of the difficulty to search this axiomatization is the failure of proofs by cases (cf. Proposition A.11).

# B   The non-modal companion method

In this appendix we introduce the non-modal companion method, which provides a really simple method to discard the validity of some modal formulas. This method does not work for all non valid formulas (but it works for (K)), but when it works it gives us a very simple explanation why validity fails.

DEFINITION B.1

We consider new variables $r_0, r_1, r_2, \ldots$ The *non-modal companion formula* of a modal formula is the non-modal formula that results from replacing every subformula of the form $\Box \varphi$ with $r_n \to \varphi$ where $n$ is the modal degree of this $\Box$ occurrence. For example, the non-modal companion of $\Box(p \to \Box q) \to (\Box p \to \Box q)$ is the formula $(r_1 \to (p \to (r_0 \to q))) \to ((r_0 \to p) \to (r_0 \to q))$. The non-modal companion of $\varphi$ is denoted $\pi(\varphi)$. ⊣

Before next proposition we remind the reader that $\mathbf{\Lambda}(\mathbf{A^c})$ is a conservative expansion of $\mathbf{\Lambda}(\mathbf{A})$, and so the statements in next proposition also apply to $\vdash_{\mathbf{A}}$.

PROPOSITION B.2

1. If $\mathsf{Fr} \models^1 \varphi$, then $\vdash_{\mathbf{A^c}} \pi(\varphi)$.

2. If $\mathsf{IFr} \models^1 \varphi$, then $\{r_n \to (r_n \odot r_n) : n \in \omega\} \vdash_{\mathbf{A^c}} \pi(\varphi)$.

3. If $\mathsf{CFr} \models^1 \varphi$, then[37] $\vdash_{\mathbf{A^c}} s(\pi(\varphi))$ where $s$ is any substitution such that $s[\{r_n : n \in \omega\}] \subseteq \{0, 1\}$.

PROOF. Let us prove the first statement. Let $h$ be a homomorphism from the algebra of formulas (including $\{r_n : n \in \omega\}$ among its variables and also including canonical constants) into $\mathbf{A^c}$. We consider the Kripke model $\langle \mathbb{N}, R, V \rangle$ where

$$R(n, n') = \begin{cases} h(r_n), & \text{if } n' = n + 1 \\ 0, & \text{otherwise.} \end{cases}$$

and $V(p, n) = h(p)$ for every variable $p$. A moment of reflection shows that for every modal formula $\varphi$ (without variables in $\{r_n : n \in \omega\}$) it holds that $V(\varphi, 0) = h(\pi(\varphi))$. Using that $\mathsf{Fr} \models^1 \varphi$ we get that $h(\pi(\varphi)) = 1$.

The proofs of the other two statements are analogous since the same Kripke model construction works for these classes of frames. ∎

EXAMPLE B.3

It is worth pointing out that the other direction is not true in general. For example, if $\mathbf{A}$ is an MTL algebra and $\varphi := \Box(p \lor q) \leftrightarrow (\Box p \lor \Box q)$ then $\vdash_{\mathbf{A}} \pi(\varphi)$. On the other hand, it is known that $\varphi$ is not valid in the class $\mathsf{CFr}(\mathbf{2})$ of crisp frames over the two element Boolean algebra (just consider a frame with two different successors); and so for every residuated lattice $\mathbf{A}$, it holds that $\mathsf{CFr}(\mathbf{A}) \not\models^1 \varphi$ (cf. [12]). ⊣

In the case of the normality axiom (K) the non-modal companion gives a very simple explanation about why in general this formula is not valid. Let us assume that $\mathsf{Fr} \models^1$ (K); then $\vdash_{\mathbf{A}} \pi(\mathsf{K})$ by Proposition B.2, i.e., $\vdash_{\mathbf{A}} (r_0 \to (p \to q)) \to ((r_0 \to p) \to (r_0 \to q))$. Replacing $p$ by $a$, $q$ by $a \odot a$ and $r_0$ by $a$ it clearly follows that $a = a \odot a$ for every element $a \in A$. That is, we have just proved that if $\mathsf{Fr} \models^1$ (K) then $\mathbf{A}$ is a Heyting algebra; and so we could have used the non-modal companion to give an alternative proof of the first item in Corollary 3.13.

Although in Example B.3 we saw that the method of non-modal companions is in general not enough to characterize all non valid formulas next we state in Corollary B.5 a particular case where it is indeed enough.

PROPOSITION B.4

Let $\delta(p_1, \ldots, p_n)$ and $\varepsilon(p)$ be two non-modal formulas such that $\delta$ is non decreasing over $\mathbf{A^c}$ in any of its arguments and $\varepsilon$ is an expanding formula over $\mathbf{A^c}$ (i.e., $\vdash_{\mathbf{A^c}} p \to \varepsilon(p)$). And let $\varphi_1, \ldots, \varphi_n, \varphi$ be non-modal formulas. If $\vdash_{\mathbf{A^c}} \delta(r \to \varphi_1, \ldots, r \to \varphi_n) \to \varepsilon(r \to \varphi)$, $r$ being a variable not appearing in $\{\varphi_1, \ldots, \varphi_n, \varphi\}$, then $\mathsf{Fr} \models^1 \delta(\Box \varphi_1, \ldots, \Box \varphi_n) \to \Box \varepsilon(\varphi)$.

PROOF. Take any frame $\mathfrak{F}$, a world $w$ in $\mathfrak{F}$, and a valuation $V$ on $\mathfrak{F}$. Then, for every world $w'$ it holds

---

[37] If $\mathbf{A}$ has the property that 1 is join irreducible (e.g., subdirectly irreducible residuated lattices, chains, etc.) then we can rewrite the conclusion as $\{r_n \lor \neg r_n : n \in \omega\} \vdash_{\mathbf{A^c}} \pi(\varphi)$.

that

$$V(\delta(\Box\varphi_1,\ldots,\Box\varphi_n),w) =$$
$$\delta^{\mathbf{A^c}}(V(\Box\varphi_1,w),\ldots,V(\Box\varphi_n,w)) \leqslant$$
$$\delta^{\mathbf{A^c}}(R(w,w') \to V(\varphi_1,w'),\ldots,R(w,w') \to V(\varphi_n,w')) \leqslant$$
$$\varepsilon^{\mathbf{A^c}}(R(w,w') \to V(\varphi,w')) =$$
$$\varepsilon^{\mathbf{A^c}}(R(w,w')) \to \varepsilon^{\mathbf{A^c}}(V(\varphi,w')) \leqslant$$
$$R(w,w') \to \varepsilon^{\mathbf{A^c}}(V(\varphi,w')) =$$
$$R(w,w') \to V(\varepsilon(\varphi),w').$$

Therefore, $V(\delta(\Box\varphi_1,\ldots,\Box\varphi_n),w) \leqslant \bigwedge\{R(w,w') \to V(\varepsilon(\varphi),w') : w' \in W\}$. Hence, $V(\delta(\Box\varphi_1,\ldots,\Box\varphi_n),w) \leqslant V(\Box\varepsilon(\varphi),w)$. ∎

For the particular case that we apply the previous proposition to $\delta := p_1 \wedge p_2$, $\varepsilon := p$ and $\varphi := \varphi_1 \wedge \varphi_2$ what we get is the validity of a particular case of (MD). Next we highlight another particular case.

Corollary B.5
Let $\varphi_1,\ldots,\varphi_n,\varphi$ be non-modal formulas. If $\vdash_{\mathbf{A^c}} ((r \to \varphi_1) \odot \ldots \odot (r \to \varphi_n)) \to (r \to \varphi)$, $r$ being a variable not appearing in $\{\varphi_1,\ldots,\varphi_n,\varphi\}$, then $\mathsf{Fr} \models^1 (\Box\varphi_1 \odot \ldots \odot \Box\varphi_n) \to \Box\varphi$.

Proof. This is a consequence of Proposition B.4; simply consider $\delta(p_1,\ldots,p_n) := p_1 \odot \ldots \odot p_n$ and $\varepsilon(p) := p$. ∎

To finish this appendix we point out that the notion of non-modal companion is also interesting from the point of view of the modal consequence relation (cf. Proposition B.2).

Proposition B.6
1. If $\Gamma \vdash^l_{\mathsf{Fr}(\mathbf{A^c})} \varphi$, then $\pi[\Gamma] \vdash_{\mathbf{A^c}} \pi(\varphi)$.

2. If $\Gamma \vdash^l_{\mathsf{IFr}(\mathbf{A^c})} \varphi$, then $\{r_n \leftrightarrow (r_n \odot r_n) : n \in \omega\}, \pi[\Gamma] \vdash_{\mathbf{A^c}} \pi(\varphi)$.

3. If $\Gamma \vdash^l_{\mathsf{CFr}(\mathbf{A^c})} \varphi$, then $s[\pi[\Gamma]] \vdash_{\mathbf{A^c}} s(\pi(\varphi))$ where $s$ is any substitution such that $s[\{r_n : n \in \omega\}] \subseteq \{0,1\}$.

Proof. The proof uses the same Kripke model given in Proposition B.2. ∎